\documentclass[11pt]{article}

\usepackage{latexsym}
\usepackage{epsfig}
\usepackage{subfig}
\usepackage{amsmath}
\usepackage{amssymb}
\usepackage{graphicx}
\usepackage{color}
\usepackage{multirow}
\usepackage{endnotes}
\usepackage{enumerate}
\usepackage{bm}
\usepackage{booktabs}
\usepackage[compact]{titlesec}
\usepackage{lscape}
\usepackage{natbib}

\usepackage{mathtools}
\allowdisplaybreaks

\titlespacing{\section}{0pt}{10pt}{*0}
\titlespacing{\subsection}{0pt}{5pt}{*0}
\titleformat{\section}{\large\bfseries}{\thesection}{1em}{}
\titleformat{\subsection}{\normalsize\bfseries}{\thesubsection}{1em}{}

\usepackage[a4paper, total={6.6in,9in}]{geometry}

\usepackage[colorlinks=true,breaklinks=true,bookmarks=true,urlcolor=blue,
     citecolor=blue,linkcolor=blue,bookmarksopen=false,draft=false]{hyperref}

\newcommand{\ol}[1]{\overline{#1}}
\newcommand{\oh}[1]{\widehat{#1}}
\newcommand{\ot}[1]{\widetilde{#1}}
\newcommand{\br}[1]{\breve{#1}}
\newcommand{\ru}[1]{\lceil #1 \rceil}
\newcommand{\rd}[1]{\lfloor #1 \rfloor}

\newcommand{\pr}{\mathbf{P}}
\newcommand{\E}{\mathbf{E}}
\newcommand{\Var}{\mathbf{Var}}

\newcommand{\T}{\mathrm{T}}

\newcommand{\dr}{\mathrm{d}}
\newcommand{\e}{\mathrm{e}}

\newcommand{\Cov}{\mathrm{Cov}}

\newcommand{\I}{\mathbb{I}}
\newcommand{\R}{\mathbb{R}}

\newcommand{\N}{\mathbb{N}}

\newcommand{\D}{\mathbb{D}}

\newcommand{\mE}{\mathcal{E}}
\newcommand{\mF}{\mathcal{F}}
\newcommand{\mA}{\mathcal{A}}
\newcommand{\mK}{\mathcal{K}}
\newcommand{\mJ}{\mathcal{J}}

\newcommand{\mB}{\mathcal{B}}

\newcommand{\mP}{\mathcal{P}}
\newcommand{\mG}{\mathcal{G}}
\newcommand{\mS}{\mathcal{S}}
\newcommand{\mZ}{\mathcal{Z}}
\newcommand{\mH}{\mathcal{H}}

\newtheorem{theorem}{Theorem}
\newtheorem{lemma}{Lemma}
\newtheorem{corollary}{Corollary}

\newtheorem{remark}{Remark}
\newtheorem{proposition}{Proposition}
\newtheorem{definition}{Definition}
\newtheorem{assumption}{Assumption}
\newenvironment{proof}{\hspace{0.20in}{\textbf{Pf:}}}  
{\hfill \rule{1.0ex}{1.0ex}}

\begin{document}


\begin{center}

\Large{{\bf On the Control of Fork-Join Networks}}

\vspace{0.25in} \normalsize

Erhun \"{O}zkan\footnote{Marshall School of Business, University of Southern California, Email: Erhun.Ozkan.2018@marshall.usc.edu.} and Amy R. Ward\footnote{Marshall School of Business, University of Southern California, Email: amyward@marshall.usc.edu.} \\

\end{center}

\begin{abstract}

Networks in which the processing of jobs occurs both sequentially and in parallel are prevalent in many application domains, such as computer systems, healthcare, manufacturing, and project management. The parallel processing of jobs gives rise to synchronization constraints that can be a main reason for job delay. In comparison with feedforward queueing networks that have only sequential processing of jobs, the approximation and control of networks that have synchronization constraints is less understood. One well-known modeling framework in which synchronization constraints are prominent is the fork-join processing network. 

Our objective is to find scheduling rules for fork-join processing networks with multiple job types in which there is first a fork operation, then activities that can be performed in parallel, and then a join operation. The difficulty is that some of the activities that can be performed in parallel require a shared resource. We solve the scheduling problem for that shared server (that is, which type of job to prioritize each time the server becomes available) when that server is in heavy traffic and prove an asymptotic optimality result.

\noindent \small{{\bf Keywords}: Fork-join processing network; Scheduling Control; Diffusion Approximation; Asymptotic Optimality. }

\end{abstract}

\section{Introduction}\label{introduction}

Networks in which processing of jobs occurs both sequentially and in parallel are prevalent in many application domains, such as computer systems (\citet{xia07}), healthcare (\citet{arm14}), manufacturing (\citet{dal92}), project management (\citet{adl95}), and the justice system (\citet{lar93}). The parallel processing of jobs gives rise to synchronization constraints that can be a main reason for job delay. Although delays in fork-join networks can be approximated under the common first-come-first-served (FCFS) scheduling discipline (\citet{ngu93,ngu94}), there is no reason to believe FCFS scheduling minimizes delay.  

\begin{figure}[htb]
\begin{center}
\includegraphics[width=0.95\textwidth]{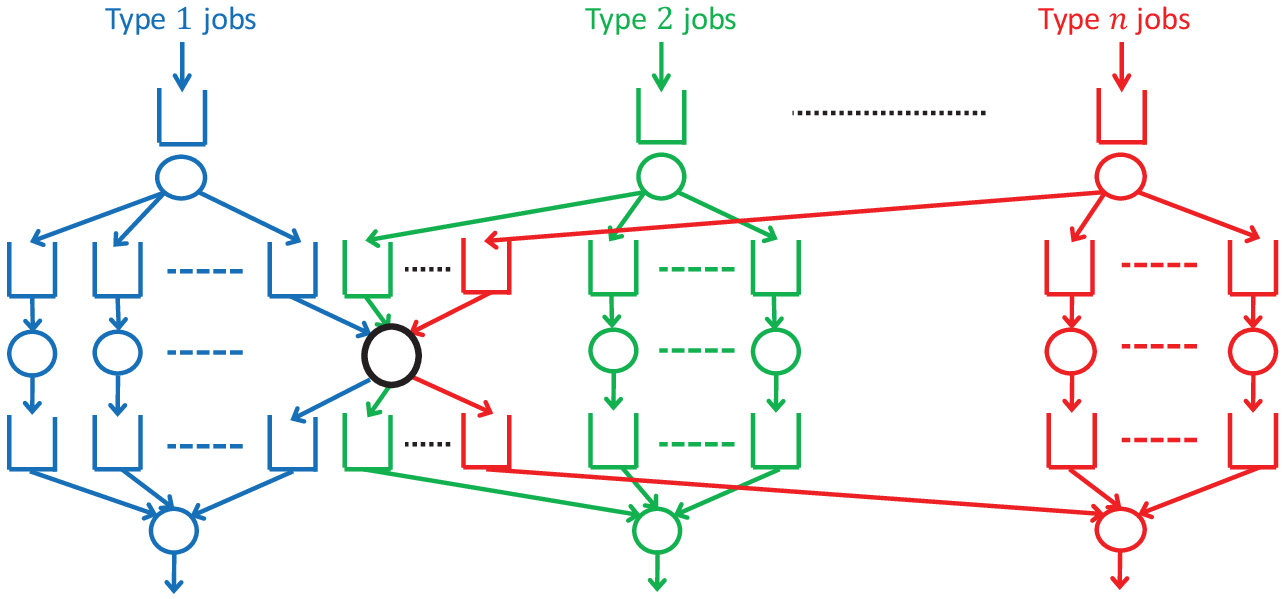}
\caption{(Color online) An example of a fork-join processing network with arbitrary number of job types, arbitrary number of forks associated with each job type, and a single shared server.}\label{fj_network_0}
\end{center}
\end{figure}

Our objective in this paper is to devise control policies that minimize delay (or, more generally, holding costs) in fork-join networks with multiple customer classes that share processing resources. For a concrete motivating example, consider the patient-flow process associated with the emergency department at Saintemarie University Hospital (cf. \citet{hub11}). An arriving patient is first triaged to determine condition severity, and then (after some potential waiting) moves to the patient management phase before being discharged. The patient management phase begins with the vital signs being taken and a first evaluation. Then, depending on the condition, there may be laboratory tests and radiology exams. Simple laboratory tests on the patient's blood and urine can be performed in parallel with the patient receiving a radiology exam, such as a CT scan. The discharge decision - whether the patient can return home or should be admitted to the hospital - cannot be made until all test results are received. Roughly speaking, we can imagine a process flow diagram such as that in Figure \ref{fj_network_0}, where the patient type corresponds to the condition severity determined at triage, the isolated operations correspond to the laboratory tests (which are necessarily associated with each individual patient), and the shared operation corresponds to the use of the CT scanner. The CT scanner is an expensive machine, and, as can be seen from the case teaching note (cf. \citet{hub11}), has a large impact on patient wait time. This motivates us to study the problem of how to schedule a shared resource that is used in parallel with other resources.

\begin{figure}[htb]
\begin{center}
\includegraphics[width=0.5\textwidth]{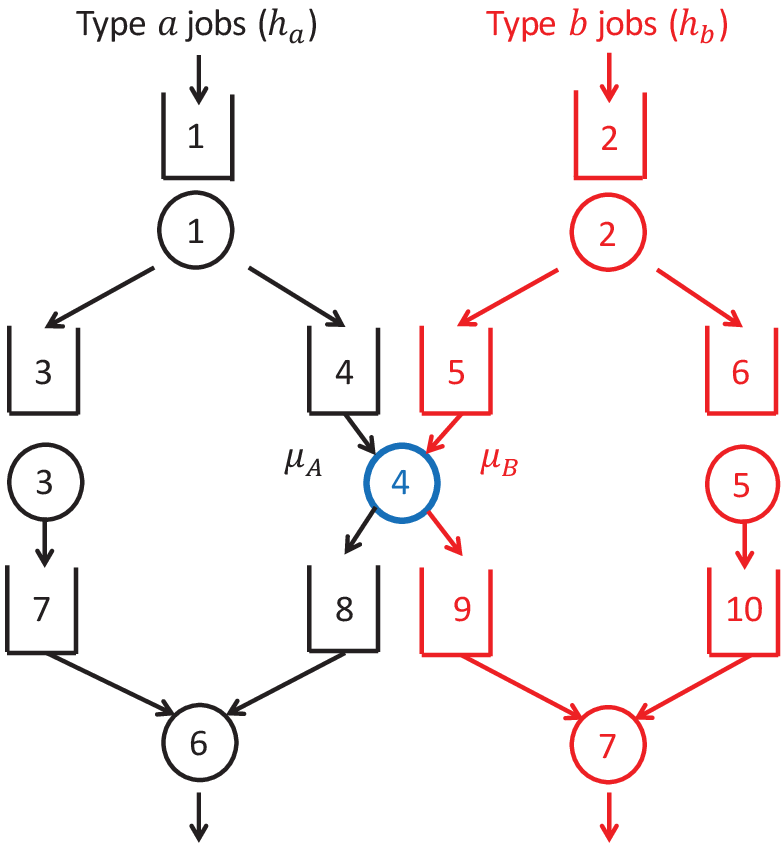}
\caption{(Color online) A fork-join processing network with two job types and a single shared server.}\label{fj_network}
\end{center}
\end{figure}

The simpler fork-join network shown in Figure \ref{fj_network} serves to illustrate why fork-join network control is difficult. In that network, there are two arriving job types ($a$ and $b$), seven servers (numbered 1 to 7), and ten buffers (numbered 1 to 10). We assume $h_a$ is the cost per unit time to hold a type $a$ job, and $h_b$ to hold a type $b$ job. Type $a$ ($b$) jobs are first processed at server 1 (2), then ``\textit{fork}'' into two jobs, one that must be processed at server 3 (5) and the other at server 4, and finally ``\textit{join}'' together to complete their processing at server 6 (7). There is synchronization because the processing at server 6 (7) cannot begin until there is at least one job in both buffers 7 and 8 (9 and 10). Server 4 processes both job types, but can only serve one job at a time. The control decision is to decide which job type server 4 should prioritize. This decision could be essentially ignored by serving the jobs in the order of their arrival regardless of type (that is, implementing FCFS policy). Another option is to always prioritizing the more expensive job type, in accordance with the well-known $c\mu$-rule. Then, if $h_a\mu_{A} \geq h_b\mu_{B}$ where $\mu_{A}$ ($\mu_{B}$) is the rate at which server 4 processes type $a$ ($b$) jobs, server 4 always prefers to work on a type $a$ job over a type $b$ job. However, when there are multiple jobs waiting at buffers 8 and 10, and no jobs waiting at buffer 9, it may be preferable to have server 4 work on a type $b$ job instead of a type $a$ job (and especially if also no jobs are waiting at buffer 7). This is because server 4 can prevent the ``join'' server 7 from idling without being the cause of server 6's forced idling. (Server 3 is the cause.)

The fork-join network control problem is too difficult to solve exactly, and so we search for an asymptotic solution. We do this under the assumption that server 4 is in heavy traffic. Otherwise, the scheduling control in server 4 has negligible impact on the delay of type $a$ and type $b$ jobs. We further assume the servers 6 and 7 are in light traffic, which focuses attention on when the required simultaneous processing of jobs at those servers forces idling. The servers 1, 2, 3, and 5 can all be in either light or heavy traffic.  

In the aforementioned heavy traffic regime, we formulate and solve an approximating diffusion control problem (DCP).  The DCP solution matches the number of jobs in buffer 4 to that in buffer 3, except when the total number of jobs waiting for processing by server 4 is too small for that to be possible. The implication is that when server 3 is in light traffic, so that buffer 3 is empty, buffer 4 is empty and all jobs waiting to be processed by server 4 are type $b$ jobs.  Otherwise, when server 3 is in heavy traffic, the control at server 4 must carefully pace its processing of type $a$ jobs to prevent ``getting ahead'' of server 3.

Our proposed policy is motivated by the observation from the DCP solution that there is no reason to have fewer type $a$ jobs in buffer 4 than in buffer 3. If server 3 can process jobs at least as quickly as server 4 can process type $a$ jobs, then the control under which server 4 gives static priority to type $a$ jobs performs well. Otherwise, we introduce a \textit{\textbf{slow departure pacing}} (SDP) control in which server 4 slows its processing of type $a$ jobs to match the departure process of type $a$ jobs from buffer 4 to the one from buffer 3.

SDP is a robust idea that is applicable to the more general network topology shown in Figure \ref{fj_network_0}. To see this, we formulate and solve approximating DCPs for the fork-join network in Figure \ref{fj_network} but with task dependent holding costs (cf. Section \ref{task_dep}), a fork-join network with an arbitrary number of ``forks'' (cf. Figure \ref{fj_network_2} in Section \ref{ex_fork}), and fork-join networks with more than two job types (cf. Figure \ref{fj_network_5} in Section \ref{ex_type_1}, Figure \ref{fj_network_0} and Remark \ref{most_general_network} in Section \ref{ex_type_1}, and Figure \ref{fj_network_3} in Section \ref{ex_type_2}). In each case, the DCP solutions suggest that, depending on the processing capacities of the servers and the network state, the servers in the network that process more than one job type should either give static priority to the more expensive job types or slow the departure process of these more expensive jobs in order to sometimes prioritize the less expensive jobs. This prevents unnecessary forced idling of the downstream ``join'' servers that process the less expensive job types, without sacrificing the speed at which the more expensive job types depart the network.

We prove that our proposed policy is asymptotically optimal in heavy traffic for the fork-join network in Figure \ref{fj_network}. To do this, we first show that the DCP solution provides a stochastic lower bound on the holding cost under any policy at every time instant. This is a strong objective, in line with that in \citet{ata05}, which implies asymptotic optimality for minimizing the expected total discounted cost or average cost over a finite time horizon (but does not guarantee asymptotic optimality for minimizing average cost over an infinite time horizon without a limit interchange result). Then, we prove a weak convergence result that implies the aforementioned lower bound is achieved under our proposed policy. That rigorous analysis suggests that our translations of the DCP solutions relevant to the more general topologies shown in Figures \ref{fj_network_0} and \ref{fj_network_3} should also perform well. 

The weak convergence result when the network operates under the SDP control is a major technical challenge for the paper. This is because the SDP control is a dynamic control that depends on the network state. In order to obtain the weak convergence, we must carefully construct random intervals on which we know the job type server 4 is prioritizing. Although this idea is similar in spirit to the random interval construction in \citet{bel01,gha13}, the proof to show convergence on the intervals is much different, due to the desired matching of the type $a$ job departure processes from the servers 3 and 4. More specifically, the interval construction is determined by tracking and comparing the job counts in buffers 3 and 4, because the SDP policy prioritizes type $a$ jobs when the number of jobs in buffer 4 exceeds that in buffer 3, and prioritizes type $b$ jobs otherwise. Then, the keys to obtaining the desired weak convergence result are as follows.
\begin{itemize}
\item On the intervals on which type $a$ jobs are prioritized, we must bound the difference between two renewal processes having different rates.
\item On the intervals on which type $b$ jobs are prioritized, we require a rate of convergence result on a light traffic $GI/GI/1$ queue that is different than any result we find in the literature (for example, the one in \citet{dai02}).
\end{itemize}
Finally, when we piece those intervals together, we see the DCP solution arise.

The remainder of this paper is organized as follows. We conclude this section with a literature review and a summary of our mathematical notation. Section \ref{model_description} specifies our model and Section \ref{asym_analysis} provides our asymptotic framework. We construct and solve an approximating DCP in Section \ref{DCP}. We introduce the SDP control in Section \ref{proposed_policy}, and specify when the proposed policy is SDP and when it is static priority. Section \ref{asym_opt_def} proves that the DCP solution provides a lower bound on the performance of any control, and Sections \ref{weak_convergence_proofs} and \ref{ud_proof} prove weak convergence under the proposed policy. Section \ref{numerical_analyses} provides extensive simulation results. In Section \ref{extensions}, we construct and solve approximating DCPs for a broader class of fork-join networks. Section \ref{conclusion} makes concluding remarks and proposes a future research direction. We separate out our rate of convergence result for a light traffic $GI/GI/1$ queue in the appendix, as that is a result of interest in its own right. We also provide the proofs of the results that use more standard methodology as well as more detailed simulation results in the appendix.

\subsection{Literature Review}\label{literature}

The inspiration for this work came from the papers \citet{ngu93,ngu94}. \citet{ngu93} establishes that a feedforward FCFS fork-join network with one job type and single-server stations in heavy traffic can be approximated by a reflected Brownian motion (RBM), and \citet{ngu94} extends this result to include multiple job types. The dimension of the RBM equals the number of stations, and its state space is a polyhedral region. In contrast to the RBM approximation for feedforward queueing networks (\citet{har96,har97,har98,har06}), the effect of the synchronization constraints in fork-join networks is to increase the number of faces defining the state space. Also in contrast to feedforward queueing networks, that number is increased further when moving from the single job type to multiple job type scenario. Although delay estimates for fork-join networks follow from the results of \citet{ngu93,ngu94}, they leave open the question of whether and how much delays can be shortened by scheduling jobs in a non-FCFS order. 

To solve the scheduling problem, we follow the ``standard Brownian machinery'' proposed in \citet{har96}. This is typically done by first introducing a heavy-traffic asymptotic regime in which resources are almost fully utilized and the buffer content processes can be approximated by a function of a Brownian motion, and second formulating an approximating Brownian control problem. Often, the dimension of the approximating Brownian control problem can be reduced by showing its equivalence with a so-called workload formulation (\citet{har97,har06}). The intriguing difference when the underlying network is a fork-join network is that the join servers must be in light traffic to arrive at an equivalent workload formulation. The issue is that otherwise the approximating problem is non-linear. This light traffic assumption is asymptotically equivalent to the assumption that processing times are instantaneous. Our simulation results suggest that our proposed control that is asymptotically optimal when the join servers are in light traffic also performs very well when the join servers are in heavy traffic. 

The papers \citet{ped14a,ped14b} are some of the few studies we find that consider the control of fork-join processing networks. In both papers, there are multiple job classes, but in \citet{ped14a} the servers can cooperate on the processing of jobs and in \citet{ped14b} they cannot. Their focus is on finding robust polices in the discrete-time setting that do not depend on system parameters and are rate stable. They do not determine whether or not their proposed policies minimize delay, which is our focus. The paper \citet{gur14} seeks to minimize delay, but in the context of a matching queue network that has only ``joins'' and no ``forks''.

An essential question to answer when thinking about controls for multiclass fork-join networks, as can be seen from the papers \citet{lu15a,lu15b,ata12}, is whether or not the jobs being joined are exchangeable; that is, whether or not a task forked from one job can be later joined with a task forked from a different job. Exchangeability is generally true in the manufacturing setting, because there is no difference between parts with the same specifications, and generally false in healthcare settings, because all paperwork and test results associated with one patient cannot be joined with another patient. The papers \citet{lu15a,lu15b} develop heavy traffic approximations for a non-exchangeable fork-join network with one arrival stream that forks into arrival streams to multiple many-server queues, and then must be joined together to produce one departure stream. The heavy-traffic approximation for the non-exchangeable network is different than for the exchangeable network, and the non-exchangeability assumption increases the problem difficulty. Their focus, different than ours, is on the effect of correlation in the service times, and there is no control. The paper \citet{ata12} looks at a fork-join network in which there is no control decision if jobs are exchangeable, and shows that the performance of the exchangeable network lower bounds the performance of the non-exchangeable network. Then, they propose a control for the non-exchangeable network that achieves performance very close to the exchangeable network. In comparison to the aforementioned papers, the exchangeability assumption is irrelevant in our case. This is because we assume head-of-line processing for each job type, so that the exact same type $a$ ($b$) jobs forked from server 1 (2) are the ones joined at server 6 (7).


\subsection{Notation and Terminology}

The set of nonnegative integers is denoted by $\N$ and the set of strictly positive integers are denoted by $\N_+$. The $k$ dimensional ($k\in\N_+$) Euclidean space is denoted by $\R^k$, $\R_+$ denotes $[0,+\infty)$, and $0_{k}$ is the zero vector in $\R^{k}$. For $x,y\in \R$, $x\vee y:= \max\{x,y\}$, $x\wedge y:= \min\{x,y\}$, and $(x)^+:=x\vee 0$. For any $x\in \R$, $\rd{x}$ ($\ru{x}$) denotes the greatest (smallest) integer which is smaller (greater) than or equal to $x$. The superscript $'$ denotes the transpose of a matrix or vector. 

For each $k\in \N_+$, $\D^k$ denotes the the space of all $\omega: \R_+\rightarrow \R^k$ which are right continuous with left limits. Let $\textbf{0}\in \D$ be such that $\textbf{0}(t)=0$ for all $t\in\R_+$. For $\omega\in \D$ and $T\in \R_+$, we let $\Vert \omega \Vert_T := \sup_{0\leq t \leq T} |\omega(t)|$. We consider $\D^k$ endowed with the usual Skorokhod $J_1$ topology (cf. Chapter 3 of \citet{bil99}). Let $\mB(\D^k)$ denote the Borel $\sigma$-algebra on $\D^k$ associated with Skorokhod $J_1$ topology. By Theorem 11.5.2 of \citet{whi02}, $\mB(\D^k)$ coincides with the Kolmogorov $\sigma$-algebra generated by the coordinate projections. For stochastic processes $W^r$, $r\in\R_+$, and $W$ whose sample paths are in $\D^k$ for some $k\in\N_+$, ``$W^r \Longrightarrow W$'' means that the probability measures induced by $W^r$ on $(\D^k,\mB(\D^k))$ weakly converge to the one induced by $W$ on $(\D^k,\mB(\D^k))$ as $r\rightarrow \infty$. For $x,y\in \D$, $x\vee y$, $x\wedge y$, and $(x)^+$ are processes in $\D$ such that $(x\vee y)(t):=x(t)\vee y(t)$, $(x\wedge y)(t):= x(t)\wedge y(t)$, and $(x)^+(t):=(x(t))^+$ for all $t\in\R_+$. For $x\in \D$, we define the mappings $\Psi, \Phi:\D\rightarrow\D$ such that for all $t\geq 0$,
\begin{equation*}
\Psi(x)(t) := \sup_{0\leq s\leq t} (-x(s))^+,\quad\quad \Phi(x)(t) := x(t)+\Psi(x)(t),
\end{equation*}
where $\Phi$ is the one-sided, one-dimensional reflection map, which is introduced by \citet{sko61}.

Let $\mZ=\{1,2,\ldots,n\}$ and $X_i$ be a process in $\D$ for each $i\in \mZ$. Then $(X_i,i\in \mZ)$ denotes the process $(X_1,X_2,\ldots,X_n)$ in $\D^n$. We denote $e$ as the deterministic identity process in $\D$ such that $e(t)=t$ for all $t\geq 0$. We abbreviate the phrase ``uniformly on compact intervals'' by ``u.o.c.'', ``almost surely'' by ``a.s.''. We let $\xrightarrow{a.s.}$ denote almost sure convergence and $\overset{d}{=}$ denote ``equal in distribution''. $\I$ denotes the indicator function and $BM_q(\theta,\Sigma)$ denotes a Brownian motion with drift vector $\theta$ and covariance matrix $\Sigma$ which starts at point $q$. The big-$O$ notation is denoted by $O(\cdot)$, i.e., if $x:\R_+\rightarrow\R$ and $y:\R_+\rightarrow\R$ are two functions, then $x(t)=O(y(t))$ as $t\rightarrow\infty$ if and only if there exist constants $C$ and $t_0$ such that $|x(t)|\leq C|y(t)|$ for all $t\geq t_0$. Lastly, $o_p(\cdot)$ is the little-$o$ in probability notation, i.e., if $\{X_n,n\in\N\}$ and $\{Y_n,n\in\N\}$ are sequences of random variables, then $X_n=o_p(Y_n)$ if and only if $|X_n|/|Y_n|$ converges in probability to 0.

\section{Model Description}\label{model_description}

We consider the control of the fork-join processing network depicted in Figure \ref{fj_network}. In this network, there are 2 job types, 7 servers, 10 buffers, and 8 activities. The set of job types is denoted by $\mJ$, where $\mJ = \{a,b\}$ and $a$ and $b$ denote the type $a$ and type $b$ jobs, respectively. The set of servers is denoted by $\mS$, where $\mS=\{1,2,\ldots,7\}$. The set of buffers is denoted by $\mK$, where $\mK=\{1,2,\ldots,10\}$, and the set of activities is denoted by $\mA$ where $\mA=\{1,2,3,A,B,5,6,7\}$. Except for server 4, each server is associated with a single activity. Server 4 is associated with two activities, denoted by $A$ and $B$, which are processing type $a$ jobs from buffer 4 and type $b$ jobs from buffer 5, respectively. Both server 6 and server 7 deplete jobs from 2 different buffers. Note that these two servers are join servers and process jobs whenever both of the corresponding buffers are nonempty. Hence, both server 6 and server 7 are associated with a single activity, namely activities 6 and 7, respectively. Let $s:\mA\rightarrow\mS$ be a function such that $s(j)$ denotes the server in which activity $j$, $j\in\mA$ takes place. Let $f:\mK\backslash\{1,2\}\rightarrow\mA$ be a function such that $f(k)$ denotes the activity which feeds buffer $k$, $k\in\mK\backslash\{1,2\}$. Lastly, let $d:\mK\rightarrow\mA$ be a function such that $d(k)$ denotes the activity which depletes buffer $k$, $k\in\mK$. For example, $s(A)=s(B)=4$, $f(4)=1$ and $d(4)=A$.

\subsection{Stochastic Primitives}\label{primitives}
We assume that all the random variables and stochastic processes are defined in the same complete probability space $(\Omega, \mF, \pr)$, $\E$ denotes the expectation under $\pr$, and $\pr(A, B):=\pr(A\cap B)$.

We associate the external arrival time of each job and the process time of each job in the corresponding activities with the sequence of random variables $\{\bar{v}_j(i), j\in \mJ\cup\mA\}_{i=1}^\infty$ and the strictly positive constants $\{\lambda_j, j\in\mJ\}$ and  $\{\mu_j, j\in\mA\}$. We assume that for each $j\in \mJ\cup\mA$, $\{\bar{v}_j(i)\}_{i=1}^\infty$ is a strictly positive, independent and identically distributed (i.i.d.) sequence of random variables mutually independent of $\{\bar{v}_k(i)\}_{i=1}^\infty$ for all $k\in (\mJ\cup\mA) \backslash \{j\}$, $\E[\bar{v}_j(1)]=1$, and the variance of $\bar{v}_j(1)$, denoted by $\Var(\bar{v}_j(1))$, is $\sigma_j^2$. For $j\in \mJ$, let $v_j(i):=\bar{v}_j(i)/\lambda_j$ be the interarrival time between the $(i-1)$st and $i$th type $j$ job. Then, $\lambda_j$ and $\sigma_j$ are the arrival rate and the coefficient of variation of the interarrival times of the type $j$ jobs, where $j\in\mJ$. For $j\in\{1,3,A,6\}$ ($j\in\{2,B,5,7\}$), let $v_j(i):=\bar{v}_j(i)/\mu_j$ be the service time of the $i$th incoming type $a$ ($b$) job associated with the activity $j$. Then, $\mu_j$ and $\sigma_j$ are the service rate and the coefficient of variation of the service times related to activity $j$, $j\in\mA$. For each $j\in \mJ\cup\mA$, let $V_j(0):=0$ and
\begin{equation}\label{renewal_process_def}
V_j(n):=\sum_{i=1}^n v_j(i)\quad\forall n\in\N_+,\quad\quad S_j(t):=\sup\{n\in\N:V_j(n)\leq t\}.
\end{equation}
Then, $S_j$ is a renewal process for each $j\in \mJ\cup\mA$. If $j\in \mJ$, $S_j(t)$ counts the number of external type $j$ arrivals until time $t$; if $j\in \mA$, $S_j(t)$ counts the number of service completions associated with activity $j$ until time $t$ given that the corresponding server works continuously on this activity during $[0,t]$.

\subsection{Scheduling Control and Network Dynamics}\label{SCND}

Let $T_j(t)$, $j\in\mA$, denote the cumulative amount of time server $s(j)$ devotes to activity $j$ during $[0,t]$. Then, a scheduling control is defined by the two dimensional service time allocation process $(T_{A},T_{B})$. Although a scheduling control indirectly affects $(T_6,T_7)$, since we do not have any direct control on servers 6 and 7, we exclude $(T_6,T_7)$ from the definition of the scheduling control. Let, 
\begin{subequations}\label{idle_def}
\begin{align}
I_{s(j)}(t)& :=t-T_j(t),\quad j\in\mA\backslash\{A,B\}, \label{idle_def_1} \\
I_4(t)& :=t-T_{A}(t)-T_{B}(t), \label{idle_def_2} 
\end{align}
\end{subequations}
denote the cumulative idle time of the servers during the interval $[0,t]$. For any $j\in\mA$, $S_j(T_j(t))$ denotes the total number of service completions related to activity $j$ in server $s(j)$ up to time $t$. For any $k\in\mK$, let $Q_k(t)$ be the number of jobs waiting in buffer $k$ at time $t$, $t\geq 0$, including the jobs that are being served. Then, for all $t\geq 0$,
\begin{subequations}\label{queue_length}
\begin{align}
&Q_1(t):=S_a(t)-S_1(T_1(t))\geq 0,\hspace{2cm} Q_2(t):=S_b(t)-S_2(T_2(t))\geq 0,\label{queue_length_1_1}\\
&\hspace{1cm}Q_k(t):=S_{f(k)}(T_{f(k)}(t))-S_{d(k)}(T_{d(k)}(t))\geq 0,\quad k\in\mK\backslash\{1,2\}.\label{queue_length_1_2}
\end{align}
\end{subequations}
For simplicity, we assume that initially all buffers are empty, i.e., $Q_k(0)=0$ for all $k\in\mK$. Later, we relax this assumption in Remark \ref{initial_buffer_remark}. 

We have
\begin{equation}\label{hl_policy}
V_j(S_j(T_j(t)))\leq T_j(t) < V_j(S_j(T_j(t))+1),\qquad\text{for all $j\in\mA$ and $t\geq 0$},
\end{equation}
which implies that we consider only head-of-the-line (HL) policies, where jobs are processed in FCFS order within each buffer. In this network, a task associated with a specific job cannot join a task originating in another job under the HL policies; that is, recalling our literature review (cf. Section \ref{literature}), the notion of exchangeability is not present.

It is straightforward to see that work conserving policies are more efficient than the others in this network. Hence, we ensure that all of the servers work in a work-conserving fashion by the following constraints: For all $t\geq 0$,
\begin{subequations}\label{w_conserving}
\begin{align}
& I_j(\cdot) \text{ is nondecreasing and }I_j(0)=0 \text{ for all }j\in\mS,\\  
& I_{s(d(k))}(t)\text{ increases if and only if }Q_k(t)=0, \text{ for all }k\in\{1,2,3,6\},\label{w_conserving_1}\\
& I_{4}(t)\text{ increases if and only if }Q_4(t)\vee Q_5(t)=0,\label{w_conserving_4}\\
& I_6(t)\text{ increases if and only if }Q_7(t)\wedge Q_8(t)=0,\label{w_conserving_2}\\
& I_7(t)\text{ increases if and only if }Q_9(t)\wedge Q_{10}(t)=0.\label{w_conserving_3}
\end{align}
\end{subequations}

A scheduling policy $\T:=\left(T_{A},T_{B}\right)$ is \textit{admissible} if it satisfies the following conditions: For any $T_i$, $i\in\mA$, $I_j$, $j\in\mS$, and $Q_k$, $k\in\mK$ satisfying \eqref{idle_def}, \eqref{queue_length}, \eqref{hl_policy}, and \eqref{w_conserving},
\begin{subequations}\label{sp_condition}
\begin{align}
&T_j(t)\in \mF,\quad \forall t\geq 0\text{ and } j\in\{A,B\},\label{sp_condition_1}\\
&T_j(\cdot)\text{ is continuous and nondecreasing with }T_j(0)=0,\; \forall j\in\{A,B\},\label{sp_condition_2}\\
&I_4(\cdot)\text{ is continuous and nondecreasing with }I_4(0)=0.\label{sp_condition_3}
\end{align}
\end{subequations}
Conditions \eqref{sp_condition_1}--\eqref{sp_condition_3} imply that we consider a wide range of scheduling policies including the ones which can anticipate the future. 

\subsection{The Objective}\label{obj_section}

A natural objective is to minimize the discounted expected total holding cost. Let $h_a$ and $h_b$ denote the holding cost rate per job per unit time for a type $a$ and $b$ job, respectively; and $\delta>0$ be the discount parameter. Moreover, let 
\begin{equation*}
Z_a(t):= Q_3(t)+Q_4(t)+Q_7(t)+Q_8(t),\qquad Z_b(t):= Q_5(t)+Q_6(t)+Q_9(t)+Q_{10}(t).
\end{equation*}
Then $Z_a(t)$ and $Z_b(t)$ denote the total number of type $a$ and $b$ jobs in the system except jobs waiting in buffers 1 and 2. Since $Q_1(t)$ and $Q_2(t)$ are independent of the scheduling policy, we exclude these processes from the definitions of  $Z_a(t)$ and $Z_b(t)$. Then, the expected total discounted holding cost under admissible policy $\T$ is
\begin{equation}\label{obj_1}
J_{\T}=\E\left[\int_0^{\infty} \e^{-\delta t}\left(h_a Z_a(t)+h_b Z_b(t)\right) \dr t \right],
\end{equation} 
and our objective is to find an admissible policy which minimizes \eqref{obj_1}. Another natural objective is to minimize the expected total cost up to time $t$, $t\in \R_+$, which is
\begin{equation}\label{obj_2}
J_{\T}=\E\left[\int_0^t (h_a Z_a(s)+h_b Z_b(s)) \dr s \right].
\end{equation} 
Yet another possible objective is to minimize the long-run average cost per unit time,
\begin{equation}\label{obj_3}
J_{\T}=\limsup_{t\rightarrow \infty}\frac{1}{t}\E\left[\int_0^t \left(h_a Z_a(s)+h_b Z_b(s)\right) \dr s \right].
\end{equation} 
We focus on a more challenging objective which is minimizing
\begin{equation}\label{obj_4}
\pr\left(h_a Z_a(t)+h_b Z_b(t) >x \right),\quad\quad \text{for all $t\in\R_+$ and $x>0$ }.
\end{equation}
It is possible to see that any policy that minimizes \eqref{obj_4} also minimizes the objectives \eqref{obj_1}, \eqref{obj_2}, and \eqref{obj_3}.

In this specific network, for all $t\geq 0$
\begin{equation}\label{balance_eq}
Q_3(t)+Q_7(t)=Q_4(t)+Q_8(t),\qquad Q_5(t)+Q_9(t)=Q_6(t)+Q_{10}(t).
\end{equation}
By \eqref{balance_eq}, a policy is optimal under the objective \eqref{obj_4} if and only if it is optimal under the objective of minimizing
\begin{equation}\label{obj_5}
\pr\left(h_a \left(Q_3(t)+Q_7(t)\right)+h_b \left(Q_6(t)+Q_{10}(t)\right)>x \right),\quad\quad \text{for all $t\in\R_+$ and $x>0$ }.
\end{equation}
We will focus on the objective \eqref{obj_5} from this point forward.

\section{Asymptotic Framework}\label{asym_analysis}

It is very difficult to analyze the system described in Section \ref{model_description} exactly. Even if we can accomplish this very challenging task, it is even less likely that the optimal control policy will be simple enough to be expressed by a few parameters. Therefore, we focus on finding an asymptotically optimal control policy under diffusion scaling and the assumption that server 4 is in heavy traffic. We first introduce a sequence of fork-join systems and present the main assumptions done for this study in Section \ref{seq_of_systems}. Then we formally define the fluid and diffusion scaled processes and present convergence results for the diffusion scaled workload facing server 4 and the diffusion scaled queue length processes associated with servers 1, 2, 3, and 5 in Section \ref{f_d_processes}. The question left open is to determine what should be the relationship between the control and the number of type $a$ and $b$ jobs in the workload facing server 4.

\subsection{A Sequence of Fork-Join Systems}\label{seq_of_systems}

We consider a sequence of fork-join systems indexed by $r$ where $r\rightarrow \infty$ through a sequence of values in $\R_+$. Each queueing system has the same structure defined in Section \ref{model_description} except that the arrival and service rates depend on $r$. To be more precise, in the $r$th system, we associate the external arrival time of each job and the process time of each job in the corresponding activities with the sequence of random variables $\{\bar{v}_j(i), j\in \mJ\cup\mA\}_{i=1}^\infty$, which we have defined in Section \ref{primitives}, and the strictly positive constants $\{\lambda_j^r, j\in\mJ\}$ and  $\{\mu_j^r, j\in\mA\}$ such that $v_j^r(i):=\bar{v}_j(i)/\lambda_j^r$, $j\in \mJ$ is the interarrival time between the $(i-1)$st and $i$th type $j$ job  and $v_j^r(i):=\bar{v}_j(i)/\mu_j^r$, $j\in \{1,3,A,6\}$ ($j\in \{2,B,5,7\}$) is the service time of the $i$th incoming type $a$ ($b$) job associated with the activity $j$ in the $r$th system. Therefore, $\lambda_j^r$, $j\in \mJ$ and $\mu_j^r$, $j\in \mA$ are the arrival rates and service rates in the $r$th system whereas the coefficient of variations are the same with the original system defined in Section \ref{model_description}. From this point forward, we will use the superscript $r$ to show the dependence of the stochastic processes to the $r$th queueing system. 

Next, we present our assumptions related to the system parameters. The first one is related to cost parameters.

\begin{assumption}\label{assumption_cost}
Without loss of generality, we assume that $h_a\mu_{A}\geq h_b\mu_{B}$.
\end{assumption}
This assumption implies that it is more expensive to keep type $a$ jobs than type $b$ jobs in server 4. 

Second, we make the following assumptions related to the stochastic primitives of the network.

\begin{assumption}\label{assumption_moment}
There exists a non-empty open neighborhood, $\mathcal{O}$, around $0$ such that for all $\alpha\in\mathcal{O}$, 
\begin{equation*}
\E[\e^{\alpha {v}_j(1)}]<\infty,\quad \text{for all $j\in\mJ\cup\mA$}.
\end{equation*}
\end{assumption}
Assumption \ref{assumption_moment} is the exponential moment assumption for the interarrival and service time processes. This assumption is common in the queueing literature, cf. \citet{har98,bel01,mag03,mey03}. 
 
Our final assumption concerns the convergence of the arrival and service rates, and sets up heavy traffic asymptotic regime.

\begin{assumption}\label{assumption_rate}
\begin{enumerate}
\item []
\item \label{lambda_rate_assumption} $\lambda_j^r\rightarrow \lambda_j > 0$ for all $j\in\mJ$ as $r\rightarrow\infty$.
\item \label{mu_rate_assumption} $\mu_j^r\rightarrow \mu_j>0$ for all $j\in\mA$ as $r\rightarrow\infty$.
\item \label{ht_assumption_s4_1} $\lambda_a/\mu_{A}+\lambda_b/\mu_{B}=1$.
\item \label{ht_assumption_s4_2} $r\left(\lambda_a^r/\mu_{A}^r+\lambda_b^r/\mu_{B}^r-1\right)\rightarrow \theta_4/\mu_{A}$ as $r\rightarrow\infty$, where $\theta_4$ is a constant in $\R$.
\item \label{ht_assumption_s1235} As $r\rightarrow\infty$,
\begin{subequations}\label{traffic_assumptions}
\begin{align}
& r(\lambda_a^r-\mu_1^r) \rightarrow \theta_1\in\R\cup\{-\infty\},\hspace{1cm} r(\lambda_b^r-\mu_2^r) \rightarrow\theta_2\in \R\cup\{-\infty\},\label{traffic_server_1_2}\\
& r(\lambda_a^r-\mu_3^r) \rightarrow\theta_3\in \R\cup\{-\infty\},\hspace{1cm}  r(\lambda_b^r-\mu_5^r) \rightarrow\theta_5\in \R\cup\{-\infty\}.\label{traffic_server_3_5}
\end{align}
\end{subequations}
\item \label{rate_assumption_s6s7} $\lambda_a<\mu_6$ and $\lambda_b<\mu_7$.
\item \label{equality_assumption} If $h_a\mu_{A} =  h_b\mu_{B}$, there exists an $r_0\in\R_+$ such that $h_a\mu_{A}^r\geq h_b\mu_{B}^r$ for all $r\geq r_0$.
\end{enumerate}
\end{assumption}
Parts \ref{ht_assumption_s4_1} and \ref{ht_assumption_s4_2} of Assumption \ref{assumption_rate} are the heavy traffic assumptions for server 4. Note that, if server 4 was in light traffic, its processing capacity would be high, thus any work-conserving control policy would perform well. Moreover, as explained in the emergency department example in Section \ref{introduction}, server 4 represents expensive and limited resources (e.g. a CT scanner) which needs to process multiple job types, thus it is expected from server 4 to be in heavy traffic in such a setting. By Parts \ref{lambda_rate_assumption}, \ref{mu_rate_assumption}, and  \ref{ht_assumption_s4_1}, we have $\mu_{A}>\lambda_a$ and $\mu_{B}>\lambda_b$. Part \ref{ht_assumption_s1235} states that each of the servers 1, 2, 3, and 5 can be either in light or heavy traffic. On the one hand, if $\theta_i=-\infty$ for some $i\in\{1,2,3,5\}$, then server $i$ is in light traffic. On the other hand, if $\theta_i\in\R$, then server $i$ is in heavy traffic. Part \ref{rate_assumption_s6s7} states that the two join servers are in light traffic. Note that, \citet{ata12,gur14,lu15a,lu15b} assume that the service processes in the join servers are instantaneous. Hence, Part \ref{rate_assumption_s6s7} of Assumption \ref{assumption_rate} extends the assumptions on the join servers made in the literature. Lastly, Part \ref{equality_assumption} of Assumption \ref{assumption_rate} is done for technical reasons that occur only when $h_a\mu_{A} =  h_b\mu_{B}$ (cf. Section \ref{asym_opt_proofs_1}).

For simplicity, we assume that 
\begin{equation*}
Q_k^r(0)=0,\quad\text{for all $k\in\mK$ and $r\in\R_+$}.
\end{equation*}
Later, we relax this assumption in Remark \ref{initial_buffer_remark}. Assumptions \ref{assumption_cost}, \ref{assumption_moment}, and \ref{assumption_rate} are assumed throughout the paper.

\subsection{Fluid and Diffusion Scaled Processes}\label{f_d_processes}

In this section, we present the fluid and diffusion scaled processes. For all $t\geq 0$, the fluid scaled processes are defined as
\begin{subequations}\label{fluid_process}
\begin{align}
\ol{S}_j^r(t)&:=r^{-2} S_j^r(r^2t),\; j\in\mJ\cup\mA,&  \ol{T}_j^r(t)&:=r^{-2} T_j^r(r^2t),\;j\in\mA,\\
\ol{Q}_j^r(t)&:=r^{-2} Q_j^r(r^2t),\quad j\in\mK,& \ol{I}_j^r(t)&:=r^{-2} I_j^r(r^2t),\quad j\in\mS,
\end{align}
\end{subequations}
and the diffusion scaled processes are defined as
\begin{subequations}\label{diffusion_process}
\begin{align}
\oh{S}_j^r(t)&:=r^{-1}\left(S_j^r(r^2t)-\lambda_j^r r^2t\right),\quad j\in\mJ,& \oh{S}_j^r(t)&:=r^{-1}\left(S_j^r(r^2t)-\mu_j^r r^2t\right),\quad j\in\mA,\label{diffusion_process_1}\\
\oh{T}_j^r(t)&:=r^{-1} T_j^r(r^2t),\quad j\in\mA,&  \oh{I}_j^r(t)&:=r^{-1} I_j^r(r^2t),\quad j\in\mS,\label{diffusion_process_2}\\
\oh{Q}_j^r(t)&:=r^{-1} Q_j^r(r^2t),\quad j\in\mK.\label{diffusion_process_3}
\end{align}
\end{subequations}

By the functional central limit theorem (FCLT), cf. Theorem 4.3.2 of \citet{whi02}, we have the following weak convergence result.
\begin{equation}\label{FCLT_on_renewal}
\left(\oh{S}_j^r, j\in\mJ\cup\mA\right)\Longrightarrow \left(\ot{S}_j, j\in\mJ\cup\mA\right),
\end{equation}
where $\ot{S}_j$ is a one-dimensional Brownian Motion for each $j\in\mJ\cup\mA$ such that $\ot{S}_j\overset{d}{=}BM_0(0,\lambda_j\sigma_j^2)$ for $j\in\mJ$, $\ot{S}_j\overset{d}{=}BM_0(0,\mu_j\sigma_j^2)$ for $j\in\mA$ and each $\ot{S}_j$ is mutually independent of $\ot{S}_i$, $i\in(\mJ\cup\mA)\backslash\{j\}$. 

For $t\geq 0$, let us define the workload process (up to a constant scale factor)
\begin{equation}\label{workload_process}
W_4^r(t):=Q_4^r(t)+\frac{\mu_{A}^r}{\mu_{B}^r} Q_5^r(t).
\end{equation}
Then, $W_4^r(t)/\mu_{A}^r$ is the expected time to deplete buffers 4 and 5 at time $t$, if no more jobs arrive after time $t$ in the $r$th system. Let $ \ol{W}_4^r(t):=r^{-2} W_4^r(r^2t)$ and $\oh{W}_4^r(t):=r^{-1} W_4^r(r^2t)$ denote the fluid and diffusion scaled workload processes, respectively. 

Next, we present the convergence of the fluid scaled processes under any work-conserving policy.  

\begin{proposition}\label{fluid_result_prop}
Under any work-conserving policy (cf. \eqref{w_conserving}), as $r\rightarrow \infty$
\begin{equation*}
\left(\ol{Q}^r_k,k\in\mK,\;\ol{W}^r_4,\; \ol{T}^r_j,j\in\mA\right)\xrightarrow{a.s.}\left(\ol{Q}_k,k\in\mK,\;\ol{W}_4,\; \ol{T}_j,j\in\mA\right)\qquad\text{u.o.c.},
\end{equation*}
where $\ol{Q}_k=\textbf{0}$ for all $k\in\mK$, $\ol{W}_4=\textbf{0}$, and $\ol{T}_j(t)=(\lambda_a/\mu_j) t$ for $j\in\{1,3,A,6\}$ and $\ol{T}_j(t)=(\lambda_b/\mu_j) t$ for $j\in\{2,B,5,7\}$ for all $t\geq 0$.
\end{proposition}
The proof of Proposition \ref{fluid_result_prop} is presented in Appendix \ref{fluid_results}. We will use Proposition \ref{fluid_result_prop} to prove convergence results for the diffusion scaled processes.

Considering Assumption \ref{assumption_rate} Part \ref {ht_assumption_s1235}, let $\mH:=\left\{i\in\{1,2,3,5\}: \theta_i\in\R\right\}$, i.e., $\mH$ is the set of servers which are in heavy traffic among the servers 1, 2, 3, and 5. Let $|\mH|$ be the cardinality of the set $\mH$, and for each $i\in\{1,2,3,5\}$, let
\begin{equation*}
\chi_i :=\begin{cases}
1, & \mbox{if server $i$ is in heavy traffic, i.e., $i\in\mH$,} \\
0, & \mbox{if server $i$ is in light traffic, i.e., $i\notin\mH$.}
\end{cases}
\end{equation*}
For all $t\geq 0$, let
\begin{subequations}\label{brew_process}
\begin{align}
\br{S}_1^r(t) &:= \chi_1 \left( \oh{S}^r_1\left(\ol{T}_1^r(t)\right)-r(\lambda_a^r-\mu_1^r)t\right) + \left(1-\chi_1\right)\left(\oh{S}_a^r(t)-\oh{Q}_1^r(t)\right), \label{brew_1}\\
\br{S}_2^r(t) &:= \chi_2 \left( \oh{S}^r_2\left(\ol{T}_2^r(t)\right)-r(\lambda_b^r-\mu_2^r)t\right) + \left(1-\chi_2\right)\left(\oh{S}_b^r(t)-\oh{Q}_2^r(t)\right). \label{brew_2}
\end{align}
\end{subequations}
Then, by using \eqref{brew_process}, we define the so-called ``\textit{netput process}'' for each buffer as
\begin{subequations}\label{netput_eq}
\begin{align}
\oh{X}_k^r(t) &:= \oh{S}_l^r(t)-\oh{S}^r_{d(k)}\big(\ol{T}_{d(k)}^r(t)\big)+r(\lambda_l^r-\mu_{d(k)}^r)t,\quad \text{for $(k,l)\in\{(1,a),(2,b)\}$},\label{netput_eq_1}\\
\oh{X}_k^r(t) &:= \br{S}_l^r(t)-\oh{S}^r_{d(k)}\big(\ol{T}_{d(k)}^r(t)\big)+r(\lambda_i^r-\mu_{d(k)}^r)t,\quad\text{for $(k,l,i)\in\{(3,1,a),(6,2,b)\}$}, \label{netput_eq_2}\\
\oh{X}_k^r(t) &:= \br{S}_l^r(t)-\oh{S}^r_{d(k)}\big(\ol{T}_{d(k)}^r(t)\big)+r\mu_{d(k)}^r\left(\frac{\lambda_i^r}{\mu_{d(k)}^r}-\frac{\lambda_i}{\mu_{d(k)}}\right)t,\quad \text{for $(k,l,i)\in\{(4,1,a),(5,2,b)\}$},\label{netput_eq_4}\\
\oh{X}_k^r(t) &:= \oh{S}_l^r(t)-\oh{Q}_i^r(t)-\oh{Q}_j^r(t)-\oh{S}^r_{d(k)}\big(\ol{T}_{d(k)}^r(t)\big)+r(\lambda_l^r-\mu_{d(k)}^r)t,\nonumber\\
&\hspace{4cm} \text{for $(k,l,i,j)\in\{(7,a,1,3),(8,a,1,4),(9,b,2,5),(10,b,2,6)\}$}. \label{netput_eq_3}
\end{align}
\end{subequations}

Let 
\begin{align}
&\hspace{1.5cm} \bm{\oh{Q}^r} :=
 \begin{bmatrix}
 \oh{Q}_1^r  \\
 \oh{Q}_2^r  \\
 \oh{Q}_3^r  \\
 \oh{Q}_6^r  \\
 \oh{W}_4^r  
 \end{bmatrix},\quad
 \bm{\oh{X}^r} :=
 \begin{bmatrix}
 \oh{X}_1^r  \\
 \oh{X}_2^r  \\
 \oh{X}_3^r  \\
 \oh{X}_6^r  \\
 \oh{X}_4^r +\frac{\mu_{A}^r}{\mu_{B}^r} \oh{X}_5^r 
 \end{bmatrix},\quad
 \bm{\oh{I}^r} :=
 \begin{bmatrix}
 \oh{I}_1^r  \\
 \oh{I}_2^r  \\
 \oh{I}_3^r  \\
 \oh{I}_5^r  \\
 \oh{I}_4^r  
 \end{bmatrix},\label{vector_def_1}\\
 &\bm{\theta} :=
 \begin{bmatrix}
 \theta_1  \\
\theta_2 \\
\theta_3-\chi_1 \theta_1  \\
\theta_5-\chi_2 \theta_2  \\
\theta_4-\chi_1 \theta_1 -\chi_2 \frac{\mu_{A}}{\mu_{B}}\theta_2 
 \end{bmatrix},\quad
 \bm{R^r} :=
 \begin{bmatrix}
\mu_1^r & 0 & 0 & 0 & 0 \\
0 & \mu_2^r & 0 & 0 & 0 \\
-\chi_1 \mu_1^r & 0 & \mu_3^r & 0 & 0 \\
0 & -\chi_2 \mu_2^r & 0 & \mu_5^r & 0  \\
-\chi_1 \mu_1^r & -\chi_2 \frac{\mu_{A}^r}{\mu_{B}^r}\mu_2^r & 0 & 0 & \mu_{A}^r
 \end{bmatrix},\label{vector_def_2}
\end{align}
and let $\bm{R}$ be a $5\times 5$ matrix which is the component-wise limit of $\bm{R}^r$. Then, we have
\begin{align}
\bm{\oh{Q}^r} &= \bm{\oh{X}^r} + \bm{R^r} \bm{\oh{I}^r},\label{queue_matrix}\\
\oh{Q}_k^r &= \oh{X}_k^r + \mu_{d(k)}^r \oh{I}_{s(d(k))}^r,\qquad k\in\{7,8,9,10\}.\label{downstream_q}
\end{align}

Let us define
\begin{equation}\label{main_covariance}
\Sigma := \bordermatrix{ ~ & 1& 2  & 3  & 6  & 4 \cr
1 & \lambda_a(\sigma_a^2+\sigma_1^2) & 0 & - \lambda_a\sigma_1^2 & 0 & - \lambda_a\sigma_1^2  \cr
2 & 0 &  \lambda_b(\sigma_b^2+\sigma_2^2) & 0  & - \lambda_b\sigma_2^2 &  - \frac{\mu_{A}}{\mu_{B}}\lambda_b \sigma_2^2  \cr
3 & - \lambda_a\sigma_1^2 & 0 & \Cov(3,3) & 0 & \Cov(3,4)  \cr
6 & 0 &  - \lambda_b\sigma_2^2 & 0 &  \Cov(6,6) & \Cov(4,6)  \cr
4 &  - \lambda_a\sigma_1^2 &  - \frac{\mu_{A}}{\mu_{B}}\lambda_b \sigma_2^2 & \Cov(3,4) & \Cov(4,6) & \Cov(4,4)\cr}
\end{equation}
where
\begin{align}
\Cov(3,3)&:=\lambda_a \left( \chi_1\sigma_1^2+ (1-\chi_1) \sigma_a^2 +\sigma_3^2\right),\quad\; \Cov(3,4):=\lambda_a \left( \chi_1\sigma_1^2+ (1-\chi_1) \sigma_a^2 \right),\nonumber\\
\Cov(4,6)&:=\frac{\mu_{A}}{\mu_{B}}\lambda_b \left( \chi_2\sigma_2^2+ (1-\chi_2) \sigma_b^2 \right),\qquad \Cov(6,6):=\lambda_b \left( \chi_2\sigma_2^2+ (1-\chi_2) \sigma_b^2+\sigma_5^2 \right),\nonumber\\
\Cov(4,4)&:= \lambda_a\left(\chi_1\sigma_1^2+(1-\chi_1) \sigma_a^2 +\sigma_{A}^2\right)+\left(\frac{\mu_{A}}{\mu_{B}}\right)^2\lambda_b\left( \chi_2\sigma_2^2+(1-\chi_2) \sigma_b^2 +\sigma_{B}^2\right).\label{covariance_4}
\end{align}
Then, we have the following weak convergence result.

\begin{proposition}\label{general_conv}
Under any work-conserving policy (cf. \eqref{w_conserving}),
\begin{equation}\label{diffusion_limit_result_0}
\left(\oh{Q}^r_1,\oh{Q}^r_2,\oh{Q}^r_3,\oh{Q}^r_6,\oh{W}^r_4\right) \Longrightarrow   \left(\ot{Q}_1,\ot{Q}_2,\ot{Q}_3, \ot{Q}_6,\ot{W}_4  \right), 
\end{equation}
where $\ot{Q}_i=\bm{0}$ for each $i\notin\mH$ and $\big(\ot{Q}_i, i\in\mH,\; \ot{W}_4\big)$ is a semimartingale reflected Brownian motion (SRBM) associated with the data $\big(P_{|\mH |},\bm{\theta}_{\mH},\Sigma_{\mH},\bm{R}_{\mH},0_{|\mH |}\big)$. $P_{|\mH |}$ is the nonnegative orthant in $\R^{|\mH |}$; $\bm{\theta}_{\mH}$ is an $|\mH |$-dimensional vector derived from the vector $\bm{\theta}$ (cf. \eqref{vector_def_2}) by deleting the rows corresponding to each $i$, $i\notin\mH$; $\Sigma_{\mH}$  and $\bm{R}_{\mH}$ are $|\mH |\times |\mH |$-dimensional matrices derived from $\Sigma$ (cf. \eqref{main_covariance}) and $\bm{R}$ (cf. \eqref{vector_def_2}) by deleting the rows and columns corresponding to each $i$, $i\notin\mH$, respectively; and $0_{|\mH |}$ is the origin in $P_{|\mH |}$.

The state space of the SRBM is  $P_{|\mH |}$; $\bm{\theta}_{\mH}$ and $\Sigma_{\mH}$ are the drift vector and the covariance matrix of the underlying Brownian motion of the SRBM, respectively; $\bm{R}_{\mH}$ is the \textit{reflection matrix}; and $0_{|\mH |}$ is the starting point of the SRBM.  
\end{proposition}

The formal definition of an SRBM can be found in Definition 3.1 of \citet{wil98b}. The proof of Proposition \ref{general_conv} is presented in Appendix \ref{proof_general_conv}.

\section{The Approximating Diffusion Control Problem}\label{DCP}

In this section, we construct an approximating diffusion control problem (DCP) with non-rigorous mathematical arguments. However, the solution of the DCP will help us to guess a heuristic control policy. Then, we will prove the asymptotic optimality of this heuristic control policy rigorously in Sections \ref{proposed_policy}, \ref{asym_opt_def}, \ref{weak_convergence_proofs}, and \ref{ud_proof}.

Parallel with the objective \eqref{obj_5}, consider the objective of minimizing
\begin{equation}\label{obj_6}
\pr\left(h_a \left(\oh{Q}_3^r(t)+\oh{Q}_7^r(t)\right)+h_b \left(\oh{Q}_6^r(t)+\oh{Q}_{10}^r(t)\right)>x \right),\quad \forall t\geq 0,\; x>0,
\end{equation}
for some $r$. Motivated by Assumption \ref{assumption_rate} Part \ref{rate_assumption_s6s7}, let us pretend that the service processes at servers $6$ and $7$ happen instantaneously. Since the diffusion scaled queue length process weakly converges to $\textbf{0}$ in a light traffic queue, we believe that considering instantaneous service rates in the downstream servers will not change the behavior of the system in the limit. In this case, jobs can accumulate in buffer 7 (8) only at the times buffer 8 (7) is empty. Similarly, jobs can accumulate in buffer 10 (9) only at the times buffer 9 (10) is empty. By this fact and \eqref{balance_eq}, 
\begin{equation}\label{queue_inst}
Q_7^r=\left(Q_4^r-Q_3^r\right)^+,\quad Q_8^r=\left(Q_3^r-Q_4^r\right)^+,\quad Q_{9}^r=\left(Q_6^r-Q_5^r\right)^+, \quad Q_{10}^r=\left(Q_5^r-Q_6^r\right)^+.
\end{equation}
By \eqref{queue_inst}, objective \eqref{obj_6} is equivalent to minimizing
\begin{equation}\label{obj_7}
\pr\left(h_a \left(\oh{Q}_3^r(t)+ \left(\oh{Q}_4^r(t)-\oh{Q}_3^r(t)\right)^+\right)+h_b \left(\oh{Q}_6^r(t)+ \left(\oh{Q}_5^r(t)-\oh{Q}_6^r(t)\right)^+\right)>x \right),\quad \forall t\geq 0,\; x>0.
\end{equation}
From the objective \eqref{obj_7}, we need approximations for $\oh{Q}_3^r$, $\oh{Q}_4^r$, $\oh{Q}_5^r$, and $\oh{Q}_6^r$ and we will achieve this goal by letting $r\rightarrow\infty$.

By \eqref{diffusion_process_3}, \eqref{workload_process}, and Proposition \ref{general_conv}, we know that $\big(\oh{Q}_3^r, \oh{Q}_6^r,\big)$ weakly converges to $\big(\ot{Q}_3, \ot{Q}_6\big)$ and $\oh{Q}_4^r+(\mu^r_A/\mu^r_B)\oh{Q}_5^r$ weakly converges to $\ot{W}_4$. At this point, let us assume that
\begin{equation*}
\left(\oh{Q}_4^r, \oh{Q}_5^r\right) \Longrightarrow \left(\ot{Q}_4, \ot{Q}_5\right).
\end{equation*}
Then we construct the following DCP: For each $x>0$ and $t\geq 0$, 
\begin{align}
\min\quad &\pr\left(h_a \left(\ot{Q}_3(t)+\left(\ot{Q}_4(t)-\ot{Q}_3(t)\right)^+\right)+h_b \left(\ot{Q}_6(t)+\left(\ot{Q}_5(t)-\ot{Q}_6(t)\right)^+\right)>x \right),\nonumber\\
\text{s.t.}\quad & \ot{Q}_4(t) + \frac{\mu_{A}}{\mu_{B}} \ot{Q}_5(t) = \ot{W}_4(t), \label{DCP_4}\\
& \ot{Q}_k(t)\geq 0,\quad\text{for all }k\in\{4,5\}. \nonumber
\end{align}
Intuitively, we want to minimize the total cost by splitting the total workload in server 4 to buffers 4 and 5 in the DCP \eqref{DCP_4}. Now, we will consider DCP \eqref{DCP_4} path-wise. Let $\omega$ ($\omega\in \Omega$) denote a sample path of the processes in DCP \eqref{DCP_4} and for any $F:\Omega\rightarrow \D$, $F(\omega_t)$ denote the value of the process $F$ at time $t$ in the sample path $\omega$. Then, consider the following optimization problem for each $\omega\in \Omega$ and $t\geq 0$.
\begin{subequations}\label{opt_problem_5}
\begin{align}
\min\quad &h_a \left(\ot{Q}_4(\omega_t)-\ot{Q}_3(\omega_t)\right)^+ + h_b \left(\ot{Q}_5(\omega_t)-\ot{Q}_6(\omega_t)\right)^+,\label{opt_problem_5_1}\\
\text{s.t.}\quad & \ot{Q}_4(\omega_t) + \frac{\mu_{A}}{\mu_{B}} \ot{Q}_5(\omega_t) = \ot{W}_4(\omega_t), \label{opt_problem_5_2}\\
& \ot{Q}_k(\omega_t)\geq 0,\quad\text{for all }k\in\{4,5\}. \label{opt_problem_5_3}
\end{align}
\end{subequations}
Note that, we exclude the term $h_a\ot{Q}_3(\omega_t)+h_b\ot{Q}_6(\omega_t)$ from the objective function \eqref{opt_problem_5_1} because $\ot{Q}_3(\omega_t)$ and $\ot{Q}_6(\omega_t)$ are independent of the decision variables $\ot{Q}_4(\omega_t)$ and $\ot{Q}_5(\omega_t)$. Although the objective function \eqref{opt_problem_5_1} is nonlinear, the optimization problem \eqref{opt_problem_5} is easy to solve because it has linear constraints. Moreover, Lemma \ref{opt_solution}, provides a closed-form solution for \eqref{opt_problem_5}.
  
\begin{lemma}\label{opt_solution}
Consider the optimization problem
\begin{align*}
\min\quad &h_a \left(q_4-q_3\right)^+ + h_b \left(q_5-q_6\right)^+,\\
\text{s.t.}\quad & q_4 + \frac{\mu_{A}}{\mu_{B}} q_5 = w_4, \\
& q_4\geq 0,\;q_5\geq 0, 
\end{align*}
where $q_4$ and $q_5$ are the decision variables, all of the parameters are nonnegative, and $h_a\mu_{A} \geq h_b\mu_{B}$. Then $q_4=q_3\wedge w_4$ and $q_5=(\mu_{B}/\mu_{A})(w_4-q_3)^+$ is an optimal solution of this problem.
\end{lemma}

\begin{proof}
Replacing $q_4$ with $w_4-(\mu_{A}/\mu_{B}) q_5$ gives us the following equivalent optimization problem which has only one decision variable.
\begin{subequations}\label{opt_problem_7}
\begin{align}
\min\quad &h_a \left(w_4-\frac{\mu_{A}}{\mu_{B}} q_5-q_3\right)^+ + h_b \left(q_5-q_6\right)^+,\label{opt_problem_7_1}\\
\text{s.t.}\quad & 0\leq q_5 \leq \frac{\mu_{B}}{\mu_{A}}w_4. \label{opt_problem_7_2}
\end{align}
\end{subequations}
We will solve the optimization problem \eqref{opt_problem_7} case by case. First, consider the case in which $w_4\leq q_3$. If $q_5=0$, then objective function value becomes $0$, which is the lowest possible objective function value of the problem \eqref{opt_problem_7}. Therefore, $q_4=w_4$ and $q_5=0$ is an optimal solution when $w_4\leq q_3$.

Second, consider the case where $w_4>q_3$. In this case, first consider the following subcase: $(\mu_{B}/\mu_{A})(w_4-q_3)\leq q_6$. In this case, if $q_5=(\mu_{B}/\mu_{A})(w_4-q_3)$, then the objective function value becomes 0. Thus, $q_4=q_3$ and $q_5=(\mu_{B}/\mu_{A})(w_4-q_3)$ is an optimal solution when $w_4>q_3$ and $(\mu_{B}/\mu_{A})(w_4-q_3)\leq q_6$. Lastly, consider the subcase in which $(\mu_{B}/\mu_{A})(w_4-q_3)>q_6$. Figure \ref{opt_sln} illustrates the change of the objective function value as $q_5$ increases from $0$ to its upper bound, $(\mu_{B}/\mu_{A})w_4$.

\begin{figure}[htb]
\begin{center}
\includegraphics[width=1\textwidth]{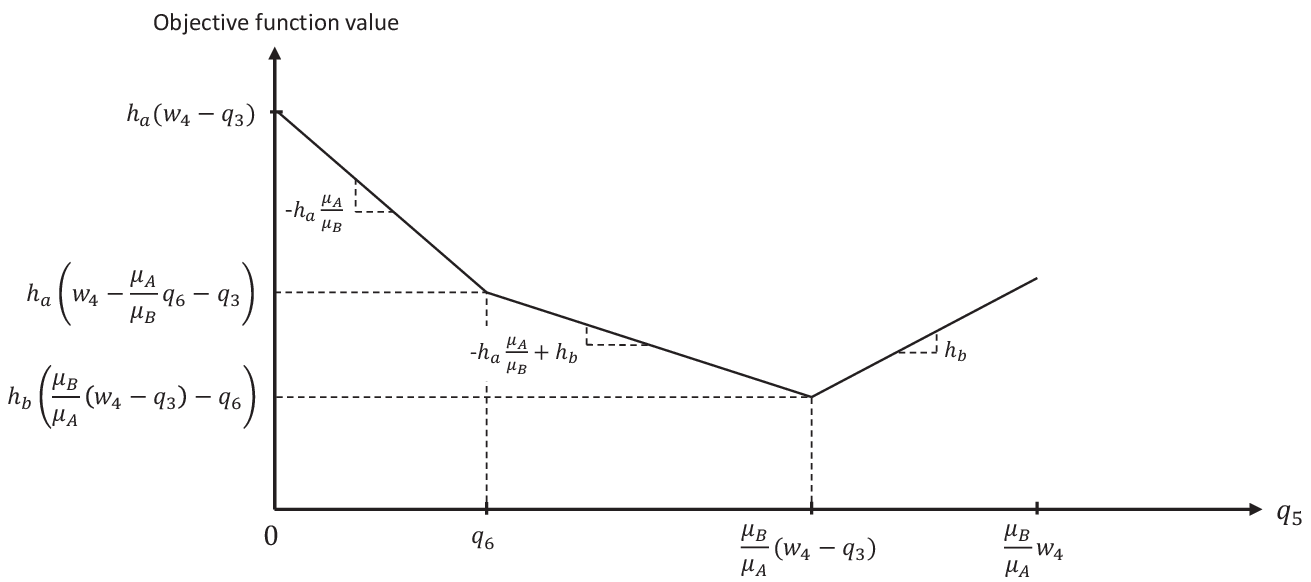}
\caption{Objective function value of the problem \eqref{opt_problem_7} with respect to different $q_5$ values.}\label{opt_sln}
\end{center}
\end{figure}
Note that, when $q_5=q_6$, the objective function value is $h_a \left(w_4-(\mu_{A}/\mu_{B}) q_6-q_3\right)$. At this point, increasing $q_5$ by a small amount, say $\epsilon$, will change the objective function value by $(-h_a \mu_{A}/\mu_{B}+h_b)\epsilon$, which is a negative number because $h_a\mu_{A} \geq h_b\mu_{B}$. When $q_5=(\mu_{B}/\mu_{A})(w_4-q_3)$, increasing $q_5$ will increase the objective function value with the rate $h_b$. Therefore, as seen in Figure \ref{opt_sln}, $q_4=q_3$ and $q_5=(\mu_{B}/\mu_{A})(w_4-q_3)$ is an optimal solution when $w_4>q_3$ and $(\mu_{B}/\mu_{A})(w_4-q_3)> q_6$.
\end{proof}

Therefore, by Lemma \ref{opt_solution}, we see that an optimal solution of the optimization problem \eqref{opt_problem_5} is $\ot{Q}_4(\omega_t)=\ot{W}_4(\omega_t)\wedge\ot{Q}_3(\omega_t)$ and $\ot{Q}_5(\omega_t)=(\mu_{B}/\mu_{A})\left(\ot{W}_4(\omega_t)-\ot{Q}_3(\omega_t)\right)^+$ for all $\omega$ and $t\geq 0$. This result and \eqref{queue_inst} imply the following proposition.

\begin{proposition}\label{lb_DCP}
\begin{align}
&\left(\ot{Q}_4,\ot{Q}_5,\ot{Q}_7,\ot{Q}_8,\ot{Q}_9,\ot{Q}_{10}\right)  \nonumber\\
&\hspace{2cm}=  \Bigg( \ot{Q}_3\wedge\ot{W}_4,\; \frac{\mu_{B}}{\mu_{A}}\left(\ot{W}_4-\ot{Q}_3\right)^+,\; \textbf{0},\; \left(\ot{Q}_3-\ot{W}_4\right)^+,\nonumber\\
&\hspace{3.5cm} \left(\ot{Q}_6-\frac{\mu_{B}}{\mu_{A}}\left(\ot{W}_4-\ot{Q}_3\right)^+\right)^+,\;\left(\frac{\mu_{B}}{\mu_{A}}\left(\ot{W}_4-\ot{Q}_3\right)^+ -\ot{Q}_6\right)^+ \Bigg)\label{DCP_sln}
\end{align}
is an optimal solution of the DCP \eqref{DCP_4}.
\end{proposition}
Note that the right hand side of \eqref{DCP_sln} is independent of the scheduling policies by Proposition \ref{general_conv}. Therefore, a control policy in which the corresponding processes weakly converge to the right hand side of \eqref{DCP_sln} is a good candidate for an asymptotically optimal policy. 

The DCP solution in Proposition \ref{lb_DCP} matches the content level of buffer 4 to that of buffer 3, except when the buffer 3 content level exceeds the total work facing server 4 (that is, the combined contents of buffers 4 and 5). This ensures that server 4 never causes server 6 to idle because of the join operation, as is evidenced by the fact that buffer 7 is always empty whereas buffer 8 sometimes has a positive content level. At the same time, server 4 prevents any unnecessary idling of server 7 by devoting its remaining effort to processing the contents of buffer 5. Sometimes, that effort is sufficient to ``keep up'' with server 5 and prevent the contents of buffer 10 from building and sometimes it is not. That is why sometimes buffer 9 has a positive content in the DCP solution and sometimes buffer 10 does (but never both simultaneously). In the next section, we formally introduce the proposed policy.

\section{Proposed Policy}\label{proposed_policy}

Our objective is to propose a policy under which the diffusion scaled queue-length processes track the DCP solution given in Proposition \ref{lb_DCP}. This is because the DCP solution provides a lower bound on the asymptotic performance of any admissible policy, as we will prove in Section \ref{asym_opt_def} (see Theorem \ref{asym_opt_1}).  

The key observation from the DCP solution in Proposition \ref{lb_DCP}  is that there is no reason for the departure process of the more expensive type $a$ jobs from server 4 to exceed that of server 3. Instead of ever letting server 4 ``get ahead'', it is preferable to have server 4 work on type $b$ jobs, so as to prevent as much forced idling at server 7 due to the join operation as possible. The only time there should have been more cumulative type $a$ job departures from server 4 than from server 3 is when the total number of jobs facing server 4 is less than that facing server 3. In that case, server 4 can outpace server 3 without forcing additional idling at server 7.

The intuition in the preceding paragraph motivates the following departure pacing policy, in which server 4 gives priority to type $a$ jobs when the number of type $a$ jobs in buffer 4 exceeds that in buffer 3 and gives priority to the type $b$ jobs in buffer 5 otherwise.

\begin{definition}\label{SDP_definition}
\textbf{Slow Departure Pacing (SDP) Policy}. The allocation process $\left(T_{A},T_{B}\right)$ satisfies
\begin{subequations}\label{formal_policy_1}
\begin{align}
& \int_0^{\infty} \I \left( Q_3(t) < Q_4(t)\right) \dr \left(t- T_A(t)\right) = 0, \label{formal_policy_1_1} \\
& \int_0^{\infty} \I \left( Q_3(t) \geq Q_4(t), Q_5(t)>0 \right) \dr \left(t- T_B(t)\right) = 0, \label{formal_policy_1_2} \\
& \int_0^{\infty} \I \left( Q_3(t) \geq Q_4(t), Q_4(t) + Q_5(t)>0 \right) \dr I_4(t) = 0, \label{formal_policy_1_3} 
\end{align}
\end{subequations}
together with \eqref{idle_def}, \eqref{queue_length}, \eqref{hl_policy}, and \eqref{w_conserving}. It is possible to see that $\left(T_{A},T_{B}\right)$ that satisfies \eqref{formal_policy_1} also satisfies \eqref{sp_condition}, and so is admissible. If $Q_3(t) <Q_4(t)$, \eqref{formal_policy_1_1} ensures that server 4 gives priority to buffer 4. If $Q_3(t) \geq Q_4(t)$ and $Q_5(t)>0$, \eqref{formal_policy_1_2} ensures that server 4 gives priority to buffer 5. \eqref{formal_policy_1_3} ensures a work-conserving control policy in server 4 when $Q_3(t) \geq Q_4(t)$.
\end{definition}

When $\mu_3 < \mu_{A}$, so that server 3 is the slower server, we use the slow departure pacing policy to determine when server 4 can allocate effort to processing type $b$ jobs without increasing type $a$ job delay. Otherwise, when $\mu_3 \geq \mu_{A}$, there is almost never extra processing power to allocate to type $b$ jobs, and so a static priority policy will perform similarly to the slow departure pacing policy (see Remark \ref{sp_optimality} regarding our numerical results).  

\begin{definition}\label{sp_definition}
\textbf{Static Priority Policy}. The allocation process $\left(T_{A},T_{B}\right)$ satisfies
\begin{subequations}\label{formal_policy_2}
\begin{align}
& \int_0^{\infty} \I \left( Q_4(t) >0\right) \dr \left(t- T_A(t)\right) = 0, \label{formal_policy_2_1} \\
& \int_0^{\infty} \I \left( Q_4(t) + Q_5(t) >0\right) \dr I_4(t) = 0, \label{formal_policy_2_2} 
\end{align}
\end{subequations}
together with \eqref{idle_def}, \eqref{queue_length}, \eqref{hl_policy}, and \eqref{w_conserving}. It is possible to see that $\left(T_{A},T_{B}\right)$ that satisfies \eqref{formal_policy_2} also satisfies \eqref{sp_condition}, and so is admissible. \eqref{formal_policy_2_1} ensures that server 4 gives static priority to buffer 4 and \eqref{formal_policy_2_2} ensures a work-conserving control policy in server 4. 
\end{definition}

The proposed policy is the SDP policy in Definition \ref{SDP_definition} when $\mu_3 < \mu_{A}$ and is the static priority policy in Definition \ref{sp_definition} when $\mu_3 \geq \mu_{A}$. We have the following weak convergence result associated with the proposed policy.

\begin{theorem}\label{diffusion_limit_result}
Under the proposed policy,
\begin{equation*}
 \left(\oh{Q}^r_k, k\in\mK,\;\oh{W}^r_4\right) \Longrightarrow  \left(\ot{Q}_k,k\in\mK\;,\ot{W}_4 \right),
\end{equation*}
where $\left(\ot{Q}_1,\ot{Q}_2,\ot{Q}_3, \ot{Q}_6,\ot{W}_4  \right)$ is defined in Proposition \ref{general_conv} and
$\left(\ot{Q}_4,\ot{Q}_5,\ot{Q}_7, \ot{Q}_8,\ot{Q}_9,\ot{Q}_{10} \right)$ is defined in Proposition \ref{lb_DCP}.
\end{theorem}

The proof of Theorem \ref{diffusion_limit_result} is presented in Section \ref{weak_convergence_proofs}. Theorem \ref{diffusion_limit_result}, the continuous mapping theorem (see, for example Theorem 3.4.3 of \citet{whi02}), and Theorem 11.6.6 of \citet{whi02} establish the asymptotic behavior of the objective function \eqref{obj_5} under the proposed policy.

\begin{corollary}\label{asym_opt_2}
Under the proposed policy, for all $t\geq 0$ and $x>0$, we have
\begin{multline}\label{asym_opt_2_1}
\left. \lim_{r\rightarrow \infty} \pr\left(h_a \left(\oh{Q}_3^r(t)+\oh{Q}_7^r(t)\right)+h_b \left(\oh{Q}_6^r(t)+\oh{Q}_{10}^r(t)\right)>x \right)\right.\\
\left.  = \pr\left(h_a \ot{Q}_3(t)+h_b \left(\ot{Q}_6(t)+\left(\frac{\mu_{B}}{\mu_{A}}\left(\ot{W}_4(t)-\ot{Q}_3(t)\right)^+- \ot{Q}_6(t)\right)^+\right)>x \right).\right.
\end{multline}
\end{corollary}

\begin{remark}\label{nonpreemtive_policy_remark}
The proposed policy is a preemptive policy. However, it is often preferred to use a non-preemptive policy. To specify a non-preemptive policy, we must specify which type of job server 4 chooses to process each time server 4 becomes free, and there are both type $a$ and type $b$ jobs waiting in buffers 4 and 5. The non-preemptive version of the SDP policy has server 4 choose to serve a type $a$ job when the number of jobs in buffer 4 exceeds that in buffer 3, and to serve a type $b$ job otherwise. The non-preemptive version of the static priority policy has server 4 always choose a type $a$ job. We expect the performance of the non-preemptive version of our proposed policy to be indistinguishable from our proposed policy in our asymptotic regime,  and we verify that the former policy performs very well by our numeric results in Section \ref{numerical_analyses}.
\end{remark}

\begin{remark}\label{implementation}
The non-preemptive version of the proposed policy (cf. Remark \ref{nonpreemtive_policy_remark}) is easy to implement in practice. The system controller does not need know any $\lambda_j$, $j\in\mJ$ or $\mu_j$, $j\in\mA$ but he needs to know whether $h_a\mu_A\geq h_b\mu_B$ or not in order to determine the job type that needs priority. Without loss of generality, suppose that the system controller knows that $h_a\mu_A\geq h_b\mu_B$. Then, he needs to know whether $\mu_3<\mu_A$ or not in order to determine which policy to implement among the non-preemptive versions of the SDP and static priority policies. Note that the simulation experiments (cf. Section \ref{numerical_results}) show that the nonpreemptive SDP policy performs well even when $\mu_3\geq\mu_A$. Thus, the nonpreemptive SDP policy is a safe choice when it is not clear whether $\mu_3<\mu_A$ or not. In order to implement this policy, the system controller only needs to keep track of the number of jobs in buffers 3 and 4 at the service completion epochs in server 4. The implementation of the nonpreemptive static priority policy is even more trivial. The system controller also needs to know whether servers 6 and 7 are in light traffic or not (cf. Part \ref{rate_assumption_s6s7} of Assumption \ref{assumption_rate}). However, since the simulation experiments (cf. Section \ref{numerical_results}) show that the nonpreemptive version of the proposed policy performs well even when at least one of the downstream servers are in heavy traffic, the latter required information is not very crucial in the implementation of the this policy. Lastly, whether server 4 is in light or heavy traffic does not affect the performance of the non-preemptive version of the proposed policy, because any work-conserving policy performs good when server 4 is in light traffic and the proposed policy is asymptotically optimal when server 4 is in heavy traffic (cf. Section \ref{asym_opt_def}).
\end{remark}

\begin{remark}
In the classical open processing networks, if a server is in light traffic, then the corresponding diffusion scaled buffer length process converges to $\textbf{0}$ (see Theorem 6.1 of \citet{che91}). However, we see in Theorem \ref{diffusion_limit_result} that although servers 6 and 7 are in light traffic, $\ot{Q}_9$ and $\ot{Q}_{10}$ are non-zero processes (moreover $\ot{Q}_8$ is a non-zero process when server 3 is in heavy traffic). Therefore, even though a join server has more than enough processing capacity, significant amount of jobs can accumulate in the corresponding buffers due to the synchronization requirements between the jobs in different buffers. This makes the control of fork-join networks more challenging than the one of classical open-processing networks.
\end{remark}

\begin{remark}\label{multiple_policy}
We will prove that the SDP and the static priority policies are asymptotically optimal when $\lambda_a\leq\mu_3<\mu_{A}$ and  $\mu_3>\lambda_a$, respectively (cf. Sections \ref{diffusion_results_1} and \ref{diffusion_results_2} and Remark \ref{sp_optimality_range}). Hence, both SDP and static priority policies are asymptotically optimal when $\lambda_a<\mu_3<\mu_{A}$ and this implies that there are many other asymptotically optimal control policies in addition to the one that we propose. However, the simulation experiments show that the static priority policy performs poorly when $\lambda_a\approx \mu_3$ but the SDP policy performs very well even when $\mu_3>\mu_{A}$ (cf. Section \ref{numerical_results}). Therefore, the performance of the proposed policy that we construct is robust with respect to the parameters $\lambda_a$, $\mu_3$, and $\mu_{A}$, which has more practical appeal.
\end{remark}

In the following section, we formally define asymptotic optimality and prove that the proposed policy is asymptotically optimal given Theorem \ref{diffusion_limit_result} and Corollary \ref{asym_opt_2}.

\section{Asymptotic Optimality} \label{asym_opt_def}

In this section, we prove that the DCP solution (cf. Proposition \ref{lb_DCP}) is a lower bound for all admissible policies with respect to the objective function \eqref{obj_5}, i.e., we have the following result.

\begin{theorem}\label{asym_opt_1}
Let $\{\T^r, r\geq 0\}$ be an arbitrary sequence of admissible policies. Then for all $t\geq 0$ and $x>0$, we have
\begin{multline}\label{asym_opt_1_1}
\left. \liminf_{r\rightarrow \infty} \pr\left(h_a \left(\oh{Q}_3^r(t)+\oh{Q}_7^r(t)\right)+h_b \left(\oh{Q}_6^r(t)+\oh{Q}_{10}^r(t)\right)>x \right)\right.\\
\left.  \geq \pr\left(h_a \ot{Q}_3(t)+h_b \left(\ot{Q}_6(t)+\left(\frac{\mu_{B}}{\mu_{A}}\left(\ot{W}_4(t)-\ot{Q}_3(t)\right)^+- \ot{Q}_6(t)\right)^+\right)>x \right).\right.
\end{multline}
\end{theorem}
Theorem \ref{asym_opt_1}  together with Corollary \ref{asym_opt_2} state that the proposed policy is asymptotically optimal.

\begin{remark}\label{obj_optimality}
In Section \ref{SCND}, we state that the objective \eqref{obj_5} implies the objectives \eqref{obj_1}, \eqref{obj_2}, and \eqref{obj_3}. Although \eqref{asym_opt_2_1} and \eqref{asym_opt_1_1} imply asymptotic optimality of the sequences of control policies with respect to the objectives \eqref{obj_1} and \eqref{obj_2}, they do not necessarily imply the same result related to the objective \eqref{obj_3}. That is because we need to change the order of the limits with respect to $r$ and $t$ to get the desired result. However, for that, we need additional results such as uniform convergence of the related processes.
\end{remark}

\subsection{Proof of Theorem \ref{asym_opt_1}}\label{asym_opt_proofs_1}

Let us consider the term in the left hand side of \eqref{asym_opt_1_1}. By \eqref{balance_eq} and the fact that $Q_k^r\geq \textbf{0}$ for all $k\in\mK$, we have
\begin{equation*}
Q_7^r (t) \geq \left(Q_4^r(t)-Q_3^r(t)\right)^+,\quad Q_{10}^r(t) \geq \left(Q_5^r(t) - Q_6^r(t)\right)^+,\qquad\forall t\geq 0.
\end{equation*}
Therefore, it is enough to prove for all $t\geq 0$ and $x>0$,
\begin{multline}\label{asym_opt_proofs_1_1}
\left. \liminf_{r\rightarrow \infty} \pr\left(h_a \left(\oh{Q}_3^r(t)+\left(\oh{Q}_4^r(t)-\oh{Q}_3^r(t)\right)^+\right)+h_b \left(\oh{Q}_6^r(t)+\left(\oh{Q}_5^r(t)-\oh{Q}_6^r(t)\right)^+\right)>x \right)\right.\\
\left.  \geq \pr\left(h_a \ot{Q}_3(t)+h_b \left(\ot{Q}_6(t)+\left(\frac{\mu_{B}}{\mu_{A}}\left(\ot{W}_4(t)-\ot{Q}_3(t)\right)^+- \ot{Q}_6(t)\right)^+\right)>x \right).\right.
\end{multline}

By \eqref{FCLT_on_renewal}, Proposition \ref{fluid_result_prop}, and Theorem 11.4.5 of \citet{whi02} (joint convergence when one limit is deterministic), we have
\begin{equation}\label{Skorokhod_1_0}
\left(\oh{S}_j^r, j\in\mJ\cup\mA,\; \ol{T}_i^r,i\in\mA\right)\Longrightarrow \left(\ot{S}_j, j\in\mJ\cup\mA,\; \ol{T}_i,i\in\mA\right).
\end{equation}
Now, we use Skorokhod's representation theorem (cf. Theorem 3.2.2 of \citet{whi02}) to obtain the equivalent distributional representations of the processes in \eqref{Skorokhod_1_0} (for which we use the same symbols and call ``Skorokhod represented versions'') such that all Skorokhod represented versions of the processes are defined in the same probability space and the weak convergence in \eqref{Skorokhod_1_0} is replaced by almost sure convergence. Then we have
\begin{equation}\label{Skorokhod_1_1}
\left(\oh{S}_j^r, j\in\mJ\cup\mA,\; \ol{T}^r_i,i\in\mA\right)\xrightarrow{a.s.} \left(\ot{S}_j, j\in\mJ\cup\mA,\; \ol{T}_i,i\in\mA\right),\quad\text{u.o.c.}\quad\text{as $r\rightarrow\infty$}.
\end{equation}
We will consider the Skorokhod represented versions of these processes from this point forward and prove \eqref{asym_opt_proofs_1_1} with respect to these processes.

By \eqref{queue_length}, \eqref{diffusion_process}, \eqref{workload_process}, \eqref{netput_eq}, \eqref{queue_matrix}, \eqref{Skorokhod_1_1}, and Proposition \ref{general_conv}, we have
\begin{equation}\label{asym_opt_proofs_1_3}
(\oh{Q}^r_3,\oh{W}^r_4,\oh{Q}^r_6)\xrightarrow{a.s.} (\ot{Q}_3,\ot{W}_4,\ot{Q}_6),\quad\text{u.o.c.},
\end{equation}
where 
\begin{equation}\label{asym_opt_proofs_1_6}
\oh{W}_4^r=\oh{Q}_4^r+\frac{\mu_{A}^r}{\mu_{B}^r}\oh{Q}_5^r,\quad\text{a.s.},
\end{equation}
and all of the processes in \eqref{asym_opt_proofs_1_3} and \eqref{asym_opt_proofs_1_6} have the same distribution with the original ones. Then by Fatou's lemma, the term in the left hand side of \eqref{asym_opt_proofs_1_1} is greater than or equal to
\begin{equation}\label{asym_opt_proofs_1_7}
\pr\left( \liminf_{r\rightarrow \infty} \left[h_a \left(\oh{Q}_3^r(t)+\left(\oh{Q}_4^r(t)-\oh{Q}_3^r(t)\right)^+\right)+h_b \left(\oh{Q}_6^r(t)+\left(\oh{Q}_5^r(t)-\oh{Q}_6^r(t)\right)^+\right)\right]>x \right).
\end{equation}
For each $t\geq 0$ and sufficiently large $r$ such that $h_a\mu_A^r\geq h_b \mu_B^r$ (note that such an $r$ exists by Assumption \ref{assumption_cost} and Parts \ref{mu_rate_assumption} and \ref{equality_assumption} of Assumption \ref{assumption_rate}), we will find a path-wise lower bound for the term
\begin{equation}\label{asym_opt_proofs_1_5}
\left[h_a \left(\oh{Q}_3^r(t)+\left(\oh{Q}_4^r(t)-\oh{Q}_3^r(t)\right)^+\right)+h_b \left(\oh{Q}_6^r(t)+\left(\oh{Q}_5^r(t)-\oh{Q}_6^r(t)\right)^+\right)\right].
\end{equation}
From this point forward, we will consider the sample paths in which $\oh{Q}_k^r(t)$, $k\in\{3,4,5,6\}$ are finite for all $r$ and $t$. By \eqref{asym_opt_proofs_1_3} and \eqref{asym_opt_proofs_1_6}, these sample paths occur with probability one. Let $\omega$ be a sample path and $\omega_t$ be defined as in Section \ref{DCP}. By \eqref{asym_opt_proofs_1_6}, \eqref{asym_opt_proofs_1_5}, and the fact that $\oh{Q}_3^r(t)$ and $\oh{Q}_6^r(t)$ are independent of the control policy, we construct the following optimization problem. For each $\omega$ in $\Omega$ except a null set and $t\geq 0$
\begin{subequations}\label{asym_opt_proofs_1_8}
\begin{align}
\min\;& h_a \left(\oh{Q}_4^r(\omega_t)-\oh{Q}_3^r(\omega_t)\right)^+ +h_b \left(\oh{Q}_5^r(\omega_t)-\oh{Q}_6^r(\omega_t)\right)^+,\label{asym_opt_proofs_1_8_1}\\
&\qquad\text{s.t.}\; \oh{Q}_4^r(\omega_t)+\frac{\mu_{A}^r}{\mu_{B}^r}\oh{Q}_5^r(\omega_t)=\oh{W}_4^r(\omega_t),\label{asym_opt_proofs_1_8_2}\\
&\qquad\quad\;\; \oh{Q}_4^r(\omega_t)\geq 0,\quad \oh{Q}_5^r(\omega_t)\geq 0,\label{asym_opt_proofs_1_8_3}
\end{align}
\end{subequations} 
where $\oh{Q}_4^r(\omega_t)$ and $\oh{Q}_5^r(\omega_t)$ are the decision variables. The optimization problem \eqref{asym_opt_proofs_1_8} has the same structure with the one presented in Lemma \ref{opt_solution}. Therefore, we can use Lemma \ref{opt_solution} to solve \eqref{asym_opt_proofs_1_8} and in the optimal solution
\begin{equation}\label{asym_opt_proofs_1_9} 
\oh{Q}_4^r(\omega_t)=\oh{Q}_3^r(\omega_t)\wedge\oh{W}_4^r(\omega_t),\quad \oh{Q}_5^r(\omega_t)=\frac{\mu_{B}^r}{\mu_{A}^r}\left(\oh{W}_4^r(\omega_t)-\oh{Q}_3^r(\omega_t)\right)^+.
\end{equation}
Therefore, by \eqref{asym_opt_proofs_1_9}, a path-wise lower bound on \eqref{asym_opt_proofs_1_5} under the admissible policy $\T^r$ is
\begin{equation}\label{asym_opt_proofs_1_10} 
h_a \oh{Q}_3^r(t)+ h_b \left(\oh{Q}_6^r(t)+\left(\frac{\mu_{B}^r}{\mu_{A}^r}\left(\oh{W}_4^r(t)-\oh{Q}_3^r(t)\right)^+ -\oh{Q}_6^r(t)\right)^+\right).
\end{equation}
When we take the $\liminf_{r\rightarrow \infty}$ of the term in \eqref{asym_opt_proofs_1_10}, by \eqref{asym_opt_proofs_1_3} and the continuous mapping theorem (specifically we use the continuity of the mapping $(\cdot)^+$ and Theorem 4.1 of \citet{whi80}, which shows the continuity of addition), \eqref{asym_opt_proofs_1_7} is greater than or equal to 
\begin{equation}\label{asym_opt_proofs_1_11} 
\pr\left(h_a \ot{Q}_3(t)+h_b \left(\ot{Q}_6(t)+\left(\frac{\mu_{B}}{\mu_{A}}\left(\ot{W}_4(t)-\ot{Q}_3(t)\right)^+- \ot{Q}_6(t)\right)^+\right)>x \right).
\end{equation}
Note that the lower bound in \eqref{asym_opt_proofs_1_11} is independent of control by Proposition \ref{general_conv} and \eqref{asym_opt_proofs_1_3}. Therefore, \eqref{asym_opt_proofs_1_11} proves \eqref{asym_opt_1_1} for the Skorokhod represented versions of the processes. Since these processes have the same joint distribution with the original ones, \eqref{asym_opt_proofs_1_11} also proves Theorem \ref{asym_opt_1}.

\section{Weak Convergence Proof}\label{weak_convergence_proofs}

In this section, we prove Theorem \ref{diffusion_limit_result}. We consider the cases $\lambda_a\leq \mu_3<\mu_A$ and $\mu_3\geq \mu_A$ separately. Note that the proposed policy is the SDP policy in \eqref{formal_policy_1} in the first case and the static priority policy in \eqref{formal_policy_2} in the second case. The proof of the second case is straightforward, but the proof of the first case is complicated because the SDP policy is a continuous-review and state-dependent policy.

\subsection{Case I:  \texorpdfstring{$\bm{\lambda_a\leq \mu_3<\mu_A}$} {$\lambda_a\leq \mu_3$} (Slow Departure Pacing Policy)}\label{diffusion_results_1}

The following result plays a crucial role in the weak convergence of $\oh{Q}^r_4$ and $ \oh{Q}^r_5$ under the proposed policy.

\begin{proposition}\label{ud_result}
Under the proposed policy, for all $\epsilon >0$ and $T>0$,
\begin{equation}\label{ud_result_eq}
\lim_r \pr \left( \big\Vert \oh{Q}^r_4- \oh{Q}^r_3\wedge\oh{W}^r_4\big\Vert_T >\epsilon \right)=0.
\end{equation}
\end{proposition}

The proof of Proposition \ref{ud_result} is presented in Section \ref{ud_proof}. By \eqref{diffusion_limit_result_0}, \eqref{ud_result_eq}, and Theorem 11.4.7 of \citet{whi02} (convergence-together theorem), we have the following joint convergence result.
\begin{equation}\label{diffusion_result_3}
\left(\oh{Q}^r_1,\oh{Q}^r_2,\oh{Q}^r_3,\oh{Q}^r_4,\oh{Q}^r_6,\oh{W}^r_4\right) \Longrightarrow \left(\ot{Q}_1,\ot{Q}_2,\ot{Q}_3,\ot{Q}_3\wedge\ot{W}_4,\ot{Q}_6,\ot{W}_4\right).
\end{equation}
By \eqref{workload_process} and \eqref{diffusion_result_3}, we have
\begin{equation}\label{s4_diffusion_1}
\oh{Q}^r_5 \Longrightarrow \frac{\mu_{B}}{\mu_{A}}(\ot{W}_4-\ot{Q}_3)^+.
\end{equation}

At this point, we invoke the Skorokhod representation theorem again for all the processes in \eqref{Skorokhod_1_0}, \eqref{diffusion_result_3}, and \eqref{s4_diffusion_1}, and we will use the same symbols again. Then, we can replace the weak convergence in \eqref{Skorokhod_1_0}, \eqref{diffusion_result_3}, and \eqref{s4_diffusion_1} with almost sure convergence in u.o.c. for the Skorokhod represented versions of the processes. Next, we will consider server 6. Let $\oh{Z}^r_6:=\oh{Q}^r_7\wedge \oh{Q}^r_8$. Then, by \eqref{netput_eq_3} and \eqref{downstream_q}, 
\begin{align}
\oh{Z}^r_6(t) &= \oh{U}^r_6(t)+\mu_6^r\oh{I}^r_6(t),\nonumber\\
\oh{U}^r_6(t) &:= \oh{S}_a^r(t)-\oh{Q}^r_1(t)-\oh{Q}^r_3(t)\vee\oh{Q}^r_4(t)-\oh{S}^r_6(\ol{T}_6^r(t))+r(\lambda_a^r-\mu_6^r)t,\label{s6_diffusion_1}
\end{align}
where $\oh{U}^r_6$ is defined in \eqref{s6_diffusion_1}. Since $\oh{I}^r_6(t)$ increases only if $\oh{Z}^r_6(t)=0$, we have a Skorokhod problem with respect to $\oh{U}^r_6$. By the uniqueness of the solution of the Skorokhod problem (cf. Proposition B.1 of \citet{bel01}) $\mu_6^r\oh{I}^r_6=\Psi(\oh{U}^r_6)$ and $\oh{Z}^r_6=\Phi(\oh{U}^r_6)$ for each $r$. By \eqref{Skorokhod_1_1} and the fact that $\big(\oh{Q}^r_1,\oh{Q}^r_3,\oh{Q}^r_4\big)\xrightarrow{a.s.}\big(\ot{Q}_1,\ot{Q}_3,\ot{Q}_3\wedge\ot{W}_4\big)$ 
\begin{equation}\label{s6_diffusion_2}
\left(\mu_6^r\oh{I}_6^r+r(\lambda_a^r-\mu_6^r)e,\; \oh{Z}_6^r\right)\xrightarrow{a.s.} \left(-\ot{S}_a+\ot{Q}_1+\ot{Q}_3+\ot{S}_6(\ol{T}_6),\;\textbf{0}\right),\qquad\text{u.o.c.},
\end{equation}
by Lemma 6.4 (ii) of \citet{che01}. By \eqref{netput_eq_3}, \eqref{downstream_q}, and \eqref{s6_diffusion_2}, 
\begin{align*}
&\oh{Q}_7^r = \oh{S}_a^r-\oh{Q}_1^r-\oh{Q}_3^r-\oh{S}^r_{6}(\ol{T}_{6}^r(t))+r(\lambda_a^r-\mu_{6}^r)e + \mu_6^r\oh{I}_6^r \xrightarrow{a.s.}  \\
& \hspace{4.7cm} \ot{S}_a-\ot{Q}_1-\ot{Q}_3-\ot{S}_6(\ol{T}_6) -\ot{S}_a+\ot{Q}_1+\ot{Q}_3+\ot{S}_6(\ol{T}_6)=\bm{0}\quad\text{u.o.c.},\\
&\oh{Q}_8^r = \oh{S}_a^r-\oh{Q}_1^r-\oh{Q}_4^r-\oh{S}^r_{6}(\ol{T}_{6}^r(t))+r(\lambda_a^r-\mu_{6}^r)e + \mu_6^r\oh{I}_6^r \xrightarrow{a.s.}  \\
& \hspace{2cm} \ot{S}_a-\ot{Q}_1-\ot{Q}_3\wedge\ot{W}_4-\ot{S}_6(\ol{T}_6) -\ot{S}_a+\ot{Q}_1+\ot{Q}_3+\ot{S}_6(\ol{T}_6)=\left(\ot{Q}_3-\ot{W}_4\right)^+\quad\text{u.o.c.}
\end{align*}

By the same way, we can derive the following result for server 7.
\begin{equation*}
\left(\oh{Q}_9^r,\oh{Q}_{10}^r\right)\xrightarrow{a.s.} \left(\left(\ot{Q}_6-\frac{\mu_{B}}{\mu_{A}}\left(\ot{W}_4-\ot{Q}_3\right)^+\right)^+,\left(\frac{\mu_{B}}{\mu_{A}}\left(\ot{W}_4-\ot{Q}_3\right)^+ -\ot{Q}_6\right)^+ \right),\qquad\text{u.o.c.}
\end{equation*}

Since the Skorokhod represented versions of the processes have the same joint distribution with the original ones, when the Skorokhod represented versions of the processes converge almost surely u.o.c., then the original processes weakly converge and we get the desired result.

\subsection{Case II:  \texorpdfstring{$\bm{ \mu_3\geq \mu_A}$}{$\mu_3\geq \mu_A$}  (Static Priority Policy)}\label{diffusion_results_2}

When $\mu_3\geq \mu_{A}$, server 3 is in light traffic because Assumption \ref{assumption_rate} Parts \ref{lambda_rate_assumption}, \ref{mu_rate_assumption}, and  \ref{ht_assumption_s4_1} imply $\mu_A>\lambda_a$. Then $\oh{Q}^r_3 \Longrightarrow \ot{Q}_3 =\bm{0}$ by Proposition \ref{general_conv}. Since server 4 gives static priority to buffer 4 over buffer 5 when $\mu_3\geq \mu_{A}$ in the proposed policy, then buffer 4 acts like a light traffic queue and $\oh{Q}^r_4 \Longrightarrow \bm{0}=\ot{Q}_3\wedge \ot{W}_4$. This implies that all of the workload in server 4 accumulates in buffer 5 and $\oh{Q}^r_5 \Longrightarrow (\mu_{B}/\mu_{A})\ot{W}_4$ by \eqref{workload_process} and Proposition \ref{general_conv}. The convergence proof of all other processes is very similar to the one presented in Section \ref{diffusion_results_1}. 

\begin{remark}\label{sp_optimality_range}
It is straightforward to see that the proof presented above holds when server 3 is in light traffic ($\lambda_a<\mu_3$). Hence, Theorem \ref{diffusion_limit_result} and the Corollary \ref{asym_opt_2} holds under the static priority policy whenever server 3 is in light traffic. Therefore, as stated in Remark \ref{multiple_policy}, the static priority policy is asymptotically optimal whenever server 3 is in light traffic.
\end{remark}

\begin{remark}
When $\mu_3\geq \mu_{A}$, we do not need Proposition \ref{ud_result} to get the desired weak convergence result. However, the proof of Proposition \ref{ud_result} is straightforward in this case. By Theorem \ref{diffusion_limit_result} and the continuous mapping theorem, $\oh{Q}^r_4- \oh{Q}^r_3\wedge\oh{W}^r_4\Longrightarrow \textbf{0}$, and this clearly gives us the desired result. 
\end{remark}

From this point forward, we will only consider the case $\lambda_a\leq \mu_3 < \mu_{A}$ and we prove Proposition \ref{ud_result} under this case in the following section. The key is the construction of random intervals such that in any given interval, we know if server 4 is giving priority to buffer 4 or to buffer 5. We first define the intervals and then prove Proposition \ref{ud_result} in the following section

\section{Proof of Proposition \ref{ud_result}}\label{ud_proof}

The SDP policy is a dynamic policy that changes the relative priorities of type $a$ and $b$ jobs depending on the system state. The analysis of such a policy requires different arguments to show $Q_4^r$ is close enough to $Q^r_3\wedge W^r_4$ to satisfy \eqref{ud_result_eq}, depending on which class has priority. This motivates us to partition the interval $[0,r^2T]$ according to the aforementioned priority rules, so that we can break the proof of \eqref{ud_result_eq} into two different parts.

We begin with the observation that type $a$ jobs are given priority at all times $t\in [0,r^2T]$ for which $Q^r_4(t)>Q^r_3(t)$, and type $b$ jobs are given priority otherwise. Then, we define ``up'' intervals during which $Q^r_4(t)>Q^r_3(t)$ and ``down'' intervals during which $Q^r_4(t)\leq Q^r_3(t)$ as follows. In the $r$th system, let $\tau_n^r:\Omega\rightarrow \R_+\cup\{+\infty\}$ be such that $\tau_0^r:=0$ and
\begin{subequations}\label{tau}
\begin{align}
\tau_{2n-1}^r &:= \inf\{t>\tau_{2n-2}^r: Q_3^r(t)=Q_4^r(t)-1\},\quad\forall n\in\N_+,\label{tau_odd}\\
\tau_{2n}^r &:= \inf\{t>\tau_{2n-1}^r: Q_3^r(t)=Q_4^r(t)\},\quad\forall n\in\N_+.\label{tau_even}
\end{align}
\end{subequations}
For completeness, if $\tau_{n_0}^r=+\infty$ for some $n_0\in \N_+$, we define $\tau_n^r:=+\infty$ for all $n\geq n_0$. Finally, we bound \eqref{ud_result_eq} using the ``up'' and ``down'' intervals so that 
\begin{align}
&\pr \left( \big\Vert \oh{Q}^r_4- \oh{Q}^r_3\wedge\oh{W}^r_4 \big\Vert_T >\epsilon \right)=\pr \left( \sup_{0\leq t\leq r^2T} \left|Q^r_4(t)- Q^r_3(t)\wedge W^r_4(t) \right| >r \epsilon \right)\nonumber\\
&\hspace{3cm}\leq \pr \left( \sup_{t\in [0,r^2T]\cap \bigcup_{n=1}^\infty [\tau_{2n-1}^r,\tau_{2n}^r)} \left|Q^r_4(t)- Q^r_3(t)\wedge W^r_4(t) \right| >r \epsilon \right)\label{ud_1}\\
&\hspace{4cm} + \pr \left( \sup_{t\in [0,r^2T]\cap \bigcup_{n=0}^\infty [\tau_{2n}^r,\tau_{2n+1}^r)} \left|Q^r_4(t)- Q^r_3(t)\wedge W^r_4(t) \right| >r \epsilon \right),\label{ud_2}
\end{align}
where $\epsilon>0$ and $T>0$ are arbitrary. For the proof, it is sufficient to prove that both of the probabilities in \eqref{ud_1} and \eqref{ud_2} converge to $0$ as $r\rightarrow\infty$.

The reason why the probability in \eqref{ud_1} converges to $0$ as $r\rightarrow\infty$ relies on the up interval construction. During an up interval, $Q^r_4>Q^r_3$, so $Q^r_3\wedge W^r_4 = Q^r_3$ by \eqref{workload_process}, and server $4$ gives priority to type $a$ jobs. Since server 4 is faster at processing type $a$ jobs than server 3, $Q^r_4$ never becomes much larger than $Q^r_3$. We make this argument rigorous in Section \ref{ud_proof_1} below.

\begin{figure}[htb]
\begin{center}
\includegraphics[width=1\textwidth]{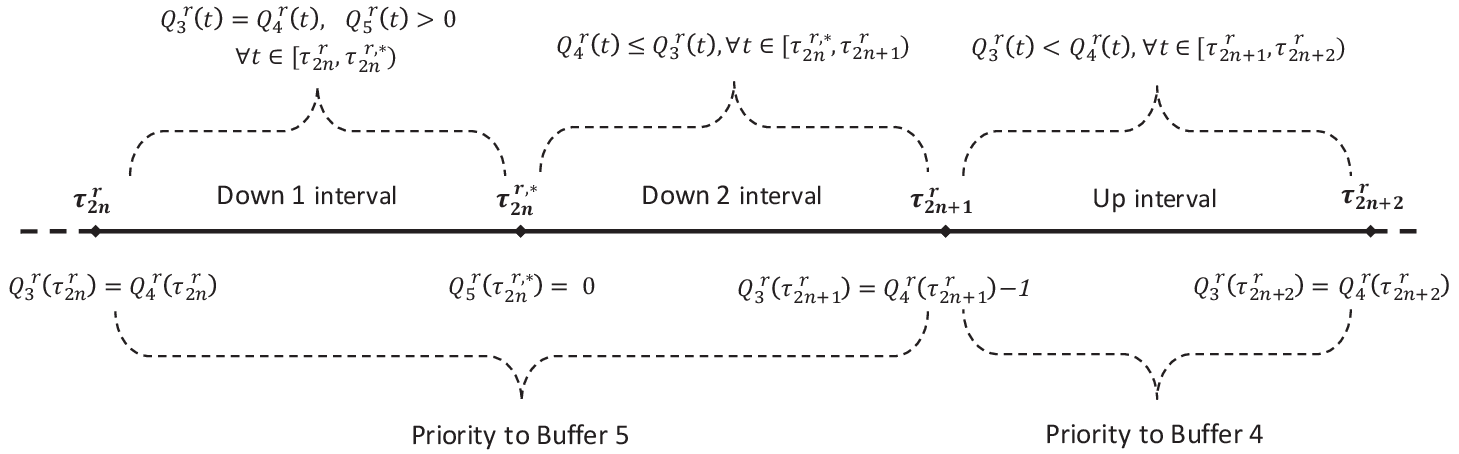}
\caption{Illustration of the down 1, down 2, and up intervals when $\lambda_a\leq\mu_3<\mu_A$.}\label{interval_figure}
\end{center}
\end{figure}

The argument to see the probability in \eqref{ud_2} converges to $0$ as $r\rightarrow\infty$ requires the splitting of the down intervals, as shown in Figure \ref{interval_figure}. To see this, first observe that at the beginning of a down interval, $Q^r_3(\tau_{2n}^r)= Q^r_4(\tau_{2n}^r)$. If $Q^r_5(\tau_{2n}^r)>0$, then server 4 works only on type $b$ jobs until the first time buffer 5 empties, defined as
\begin{equation}\label{tau_star}
\tau_{2n}^{r,*} := \begin{cases}
\inf\left\{t\in[\tau_{2n}^r,\tau_{2n+1}^r): Q_5^r(t)=0\right\}, &\parbox[t]{.4\textwidth}{if $ Q_5^r(t)=0$ for some $t\in[\tau_{2n}^r,\tau_{2n+1}^r)$ where $\tau_{2n}^r<\infty$,}\\
+\infty, & \mbox{otherwise}.
\end{cases}
\end{equation}
Note that it is possible that the next up interval starts before buffer 5 empties, in which case  $\tau_{2n}^{r,*}:=\infty$. For $t\in[\tau_{2n}^r,\tau_{2n}^{r,*})$, $Q^r_3(t)=Q^r_4(t)$. Since $Q^r_3(t)\wedge W^r_4(t) =Q^r_3(t)$ by \eqref{workload_process}, $Q^r_4(t) - Q^r_3(t)\wedge W^r_4(t) = 0$ trivially. During the interval $[\tau_{2n}^{r,*},\tau_{2n+1}^r)$ (when it exists), $Q^r_5$ behaves like a light traffic $GI/GI/1$ queue starting from the origin. Hence, $Q^r_5$ stays close to $0$, so that $Q^r_4\approx W^r_4$. Since also $Q^r_4 \leq Q^r_3$ by the down interval construction, $Q^r_3\wedge W^r_4 \approx Q^r_4$. If $Q^r_5(\tau_{2n}^r)=0$, then $\tau_{2n}^r=\tau_{2n}^{r,*}$, so that the second half of the above argument works for this case as well. We make this argument rigorous in Section \ref{ud_proof_2} below. The key is the following convergence result, that is of interest in its own right.

\paragraph{The light traffic $GI/GI/1$ Queue Convergence Rate Result}\label{light_traffic_conv_rate}

Let $\{u_k,k\in\N_+\}$ and $\{v_k,k\in\N_+\}$\footnote{We present a general result in this section, thus we use a general notation in this section and in Appendix \ref{light_traffic_conv_rate_proof}, where we present the proof of the result of this section.} be two independent sequences of strictly positive and i.i.d. random variables such that $\E[u_1]=1/\lambda$ and $\E[v_1]=1/\mu$. Consider a $GI/GI/1$ queue with infinite buffer capacity and FCFS service discipline such that the interarrival time between the $(k-1)$th and $k$th job arriving in the system after time $0$ is $u_k$ and service time of the $k$th job is $v_k$ for all $k\in\N_+$. We assume that there exists an open interval centered at zero denoted by $(-\ol{\alpha},\ol{\alpha})$ where $\ol{\alpha}>0$ such that $\E[e^{\alpha u_1}]<\infty$ and $\E[e^{\alpha v_1}]<\infty$ for all $\alpha\in(-\ol{\alpha},\ol{\alpha})$. Therefore, we assume exponential moment assumption for $u_1$ and $v_1$, i.e., all moments of $u_1$ and $v_1$ are finite. Let $Q(t)$ denote the total number of jobs in the system at time $t\geq 0$, and $Q(0)$ be a random variable which takes values in $\N$. Then, we have the following large deviations result.

\begin{proposition}\label{light_traffic_conv_rate_lemma}
Suppose that $\lambda<\mu$, i.e., the queue is in light traffic. Fix arbitrary $n\in\N_+$ and $\epsilon>0$, and suppose that there exists an $r_0\geq 1$ such that if $r\geq r_0$,
\begin{equation}\label{l_in_q}
\pr \left(  Q(0)> \frac{((\mu-\lambda)\wedge\epsilon)r}{2} \right)\leq C_0r^{2n-1}\e^{-C_1r},
\end{equation}
for some constants $C_0>0$ and $C_1>0$ which are independent of $r$. Then, there exists an $r_1\geq 1$ such that if $r\geq r_1$,
\begin{equation}\label{l_conv_rate}
\pr \left(  \sup_{0\leq t \leq r^n} Q(t)>r  \epsilon \right)\leq C_2r^{2n-1}\e^{-C_3r},
\end{equation}
such that $C_2>0$ and $C_3>0$ are constants which are independent of $r$.
\end{proposition}

The proof of Proposition \ref{light_traffic_conv_rate_lemma} is presented in Appendix \ref{light_traffic_conv_rate_proof}.

\subsection{Proof of Convergence of \eqref{ud_1} (Up Intervals)}\label{ud_proof_1}

Throughout Section \ref{ud_proof_1}, let $\gamma\in (0,1)$ be an arbitrary constant, $\alpha=\alpha(\gamma)>0$ be a constant such that $4/(1+\alpha)<\gamma$. Let $N^r:=\ru{(\mu_3^r+\epsilon_1)r^2T}$ where $\epsilon_1>0$ is an arbitrary constant. Since there is a service completion in server 3 at $\tau_{2n-1}^r$ for each $n\in \N_+$, then
\begin{equation}\label{ud_3}
\pr \left( \tau_{2N^r-1}^r \leq r^2T\right)\leq \pr \left(S_3^r(r^2T) \geq N^r\right)\leq \pr \left(S_3^r(r^2T) \geq (\mu_3^r+\epsilon_1)r^2T \right)\rightarrow 0\qquad\text{as $r\rightarrow\infty$},
\end{equation}
where \eqref{ud_3} is by functional strong law of large numbers (FSLLN) for renewal processes (cf. Theorem 5.10 of \citet{che01}).

Since, by construction, $Q^r_4(t)> Q^r_3(t)$ for all $t\in[\tau_{2n-1}^r,\tau_{2n}^r)$ and $n\in\N_+$, then $Q^r_3(t)\wedge W_4^r(t)=Q^r_3(t)$  for all $t\in[\tau_{2n-1}^r,\tau_{2n}^r)$ and $n\in\N_+$ by \eqref{workload_process}, and the  probability in \eqref{ud_1} is equal to
\begin{align}
& \pr \left( \sup_{t\in [0,r^2T]\cap \bigcup_{n=1}^\infty [\tau_{2n-1}^r,\tau_{2n}^r)} Q^r_4(t)- Q^r_3(t) >r \epsilon \right)\leq \nonumber\\
&\qquad\pr \left( \tau_{2N^r-1}^r \leq r^2T\right)+ \pr \left( \sup_{t\in [0,r^2T]\cap \bigcup_{n=1}^\infty [\tau_{2n-1}^r,\tau_{2n}^r)} Q^r_4(t)- Q^r_3(t) >r \epsilon ,\;  \tau_{2N^r-1}^r >r^2T \right)\label{ud_4}.
\end{align}
Note that the first probability in the right hand side of inequality \eqref{ud_4} converges to $0$ as $r\rightarrow\infty$ by \eqref{ud_3}. Hence, it is enough to consider the second probability in the right hand side of \eqref{ud_4} which is less than or equal to
\begin{align}
&\sum_{n=1} ^{N^r} \pr \left( \sup_{\tau_{2n-1}^r\leq t <\tau_{2n}^r} Q^r_4(t)- Q^r_3(t) >r \epsilon ,\;  \tau_{2n-1}^r \leq r^2T \right)\nonumber\\
&\qquad=\sum_{n=1} ^{N^r} \pr \left( \sup_{\tau_{2n-1}^r\leq t <\tau_{2n}^r} S_{3}^r(T_{3}^r(t))-S_{A}^r(T_{A}^r(t)) >r \epsilon ,\;  \tau_{2n-1}^r \leq r^2T\right)\label{ud_5}\\
&\qquad=\sum_{n=1} ^{N^r} \pr \left( \sup_{0\leq t <\tau_{2n}^r-\tau_{2n-1}^r} S_{3}^r(T_{3}^r(\tau_{2n-1}^r+t))-S_{A}^r(T_{A}^r(\tau_{2n-1}^r+t)) >r \epsilon ,\;  \tau_{2n-1}^r \leq r^2T\right)\nonumber\\
&\qquad \leq \sum_{n=1} ^{N^r} \pr \left( \sup_{0\leq t <\tau_{2n}^r-\tau_{2n-1}^r} S_{3}^r(T_{3}^r(\tau_{2n-1}^r)+t)-S_{A}^r(T_{A}^r(\tau_{2n-1}^r)+t) >r \epsilon ,\;  \tau_{2n-1}^r \leq r^2T\right)\label{ud_6}
\end{align}
where \eqref{ud_5} is by \eqref{queue_length_1_2}. We obtain \eqref{ud_6} in the following way. Server 4 works on buffer 4 during the whole up interval by construction. However, server 3 can be idle during an up interval. Hence, we have for all $t\in [0,\tau_{2n}^r-\tau_{2n-1}^r)$
\begin{equation}\label{ud_6_1}
T_{A}^r(\tau_{2n-1}^r+ t) =T_{A}^r(\tau_{2n-1}^r)+ t,\quad T_{3}^r(\tau_{2n-1}^r+ t) \leq T_{3}^r(\tau_{2n-1}^r)+ t,
\end{equation}
which gives \eqref{ud_6}. 

We expect the sum in \eqref{ud_6} to converge to zero because each term is the difference between two renewal processes with different rates, and the faster renewal process is the one being subtracted. To formalize this intuition, it is helpful to bound those differences by using the processes
\begin{equation*}
E_{3}^{r,n}(t) :=  \sup\left\{k\in \N: \sum_{l=1}^{k} v_{3}^r(l+B_n^r) \leq t\right\}, \qquad E_{A}^{r,n}(t) := \sup\left\{k\in \N: \sum_{l=1}^{k} v_{A}^r(l+A_n^r) \leq t\right\},
\end{equation*}
where $\sum_{l=1}^{0} x_l:=0$ for any sequence $\{x_l,l\in \N_+\}$, $A_n^r:= S_{A}^r(T_{A}^r(\tau_{2n-1}^r))$, and $B_n^r:= S_{3}^r(T_{3}^r(\tau_{2n-1}^r))$. Then by \eqref{queue_length_1_2} and \eqref{tau_odd}, $A_n^r=B_n^r-1$.
We have the following result.

\begin{lemma}\label{up_lemma} For all $t\geq 0$, $S_{3}^r(T_{3}^r(\tau_{2n-1}^r)+t) - S_{A}^r(T_{A}^r(\tau_{2n-1}^r)+t) \leq E_{3}^{r,n}(t) - E_{A}^{r,n}(t) +1$.
\end{lemma}
The proof of Lemma \ref{up_lemma} is presented in Appendix \ref{up_lemma_proof}. By Lemma \ref{up_lemma}, \eqref{ud_6} is less than or equal to
\begin{align}
&\sum_{n=1} ^{N^r} \pr \left( \sup_{0\leq t <\tau_{2n}^r-\tau_{2n-1}^r} E_{3}^{r,n}(t) - E_{A}^{r,n}(t) >r \epsilon -1 ,\;  \tau_{2n-1}^r \leq r^2T\right).\label{ud_6_0}\\
&\hspace{1cm}\leq \sum_{n=1} ^{N^r} \pr \left( \sup_{0\leq t <(\tau_{2n}^r-\tau_{2n-1}^r)\vee (r^{\gamma}T)} E_{3}^{r,n}(t) - E_{A}^{r,n}(t) >r \epsilon -1 ,\;  \tau_{2n-1}^r \leq r^2T\right)\nonumber\\
&\hspace{1cm}=\sum_{n=1} ^{N^r} \pr \left( \sup_{0\leq t < \tau_{2n}^r-\tau_{2n-1}^r} E_{3}^{r,n}(t) - E_{A}^{r,n}(t) >r \epsilon -1 ,\; \tau_{2n}^r-\tau_{2n-1}^r > r^{\gamma}T ,\; \tau_{2n-1}^r \leq r^2T \right)\nonumber\\
&\hspace{2.5cm}+ \sum_{n=1} ^{N^r} \pr \left( \sup_{0\leq t < r^{\gamma}T} E_{3}^{r,n}(t) - E_{A}^{r,n}(t) >r \epsilon -1,\; \tau_{2n}^r-\tau_{2n-1}^r \leq r^{\gamma}T ,\; \tau_{2n-1}^r \leq r^2T \right)\nonumber\\
&\hspace{1cm}\leq\sum_{n=1} ^{N^r} \pr \left(\tau_{2n}^r-\tau_{2n-1}^r > r^{\gamma}T ,\; \tau_{2n-1}^r \leq r^2T \right)\label{ud_7}\\
&\hspace{5cm} +  \sum_{n=1} ^{N^r} \pr \left( \sup_{0\leq t < r^{\gamma}T} E_{3}^{r,n}(t) - E_{A}^{r,n}(t) >r \epsilon -1 ,\; \tau_{2n-1}^r \leq r^2T \right).\label{ud_8}
\end{align}
We will show that both of the terms in \eqref{ud_7} and \eqref{ud_8} converge to $0$ as $r\rightarrow\infty$.

\subsubsection{Proof of Convergence of \eqref{ud_7} (Length of Up Intervals)}\label{length_of_ip_intervals}

In this section, we show that the sum in \eqref{ud_7} converges to $0$ as $r\rightarrow\infty$. This implies that the length of the up intervals within the time interval $[0,r^2T]$ is $o_p(r^{\delta})$ for any $\delta>0$. Let $\mE_n^r$ denote the event inside the probability in \eqref{ud_7}, i.e., 
\begin{equation}\label{lemma_up_eq_1}
\mE_n^r := \left\{\tau_{2n}^r-\tau_{2n-1}^r > r^{\gamma}T,\;  \tau_{2n-1}^r \leq r^2T \right\}.
\end{equation}
By \eqref{queue_length_1_2} and the fact that $Q_4^r(t)>Q_3^r(t)$ for all $t\in[\tau_{2n-1}^r,\tau_{2n}^r)$, \eqref{ud_7} is equal to 
\begin{align}
& \sum_{n=1} ^{N^r} \pr \left( \sup_{\tau_{2n-1}^r\leq t <\tau_{2n-1}^r+r^{\gamma}T} S_{A}^r(T_{A}^r(t))-S_{3}^r(T_{3}^r(t)) < 0,\; \mE_n^r \right)\nonumber\\
&\hspace{2cm}= \sum_{n=1} ^{N^r} \pr \left( \sup_{0\leq t<r^{\gamma}T} S_{A}^r(T_{A}^r(\tau_{2n-1}^r+t))-S_{3}^r(T_{3}^r(\tau_{2n-1}^r+t)) < 0,\;  \mE_n^r \right)\nonumber\\
&\hspace{2cm}\leq \sum_{n=1} ^{N^r} \pr \left( \sup_{0\leq t<r^{\gamma}T} S_{A}^r(T_{A}^r(\tau_{2n-1}^r)+t)-S_{3}^r(T_{3}^r(\tau_{2n-1}^r)+t) < 0,\; \mE_n^r \right),\label{lemma_up_eq_2} \\
&\hspace{2cm} \leq \sum_{n=1} ^{N^r} \pr \left( \sup_{0\leq t<r^{\gamma}T} E_{A}^{r,n}(t) - E_{3}^{r,n}(t) < 1,\;  \tau_{2n-1}^r \leq r^2T \right),\label{lemma_up_eq_7}
\end{align}
where \eqref{lemma_up_eq_2} is by \eqref{ud_6_1}, and \eqref{lemma_up_eq_7} is by Lemma \ref{up_lemma}.

Similar to \eqref{ud_6}, we consider the difference of two renewal processes in \eqref{lemma_up_eq_7}. We want to show that the probability that the number of renewals associated with the renewal process with higher renewal rate is always less than the number of renewals of the one with smaller renewal rate within a time interval of length $r^{\gamma}T$ converges to zero as $r\rightarrow\infty$. Let
\begin{equation}\label{lemma_up_eq_10_1}
\mE_{n,1}^r := \left\{ \sup_{0\leq t<r^{\gamma}T} \left|E_{3}^{r,n}(t) - \mu_{3}^r t\right| \leq \epsilon_2r^{\gamma}T,\;\sup_{0\leq t<r^{\gamma}T} \left|E_{A}^{r,n}(t) - \mu_{A}^r t\right| \leq \epsilon_2r^{\gamma}T  \right\},
\end{equation}
where $\epsilon_2$ is an arbitrary constant such that $0<\epsilon_2 < \mu_3 \wedge (( \mu_{A}- \mu_{3})/2)$. Then, the sum in \eqref{lemma_up_eq_7} is equal to
\begin{align}
&\sum_{n=1} ^{N^r} \pr \left( \sup_{0\leq t<r^{\gamma}T} E_{A}^{r,n}(t) - E_{3}^{r,n}(t) < 1,\;  \tau_{2n-1}^r \leq r^2T,\; \mE_{n,1}^r \right) \label{lemma_up_eq_11}\\
&\hspace{3cm} + \sum_{n=1} ^{N^r} \pr \left( \sup_{0\leq t<r^{\gamma}T} E_{A}^{r,n}(t) - E_{3}^{r,n}(t) < 1,\;  \tau_{2n-1}^r \leq r^2T,\; (\mE_{n,1}^r)^c \right),\label{lemma_up_eq_11_0}
\end{align}
where superscript $c$ denote the complement of a set. Note that, on the set $\mE_{n,1}^r$, we have for all $t\in [0,r^{\gamma}T)$,
\begin{equation}\label{lemma_up_eq_11_1}
\mu_{3}^r t-\epsilon_2r^{\gamma}T\leq E_{3}^{r,n}(t) \leq \mu_{3}^r t+\epsilon_2r^{\gamma}T,\qquad \mu_{A}^r t-\epsilon_2r^{\gamma}T\leq E_{A}^{r,n}(t) \leq \mu_{A}^r t+\epsilon_2r^{\gamma}T.
\end{equation}
This implies that the sum in \eqref{lemma_up_eq_11} is less than or equal to
\begin{align}
&\sum_{n=1} ^{N^r} \pr \left(  \sup_{0\leq t<r^{\gamma}T} \mu_{A}^r t-\epsilon_2r^{\gamma}T- (\mu_{3}^r t+\epsilon_2r^{\gamma}T) < 1\right) \nonumber\\
&\hspace{4cm}= N^r \pr \left( (\mu_{A}^r - \mu_{3}^r)r^{\gamma}T-2\epsilon_2r^{\gamma}T < 1\right) \label{lemma_up_eq_12_1}\\
&\hspace{4cm}= N^r \I \left( \frac{\mu_{A}^r - \mu_{3}^r}{2} -\frac{1}{2r^{\gamma}T} <\epsilon_2 \right) \rightarrow 0\qquad\text{as $r\rightarrow\infty$},\label{lemma_up_eq_12}
\end{align}
where \eqref{lemma_up_eq_12_1} is by the fact that $\mu_{A}^r > \mu_{3}^r$ for all $r\geq r_1$ for some $r_1\in \R_+$; and \eqref{lemma_up_eq_12} is by the fact that $\epsilon_2 < ( \mu_{A}- \mu_{3})/2$. Now, let us look at the sum in \eqref{lemma_up_eq_11_0}, which is less than or equal to
\begin{align}
&\sum_{n=1} ^{N^r} \pr \left( \tau_{2n-1}^r \leq r^2T,\; (\mE_{n,1}^r)^c\right)\label{lemma_up_eq_13}\\
&\hspace{2cm}\leq \sum_{n=1} ^{N^r} \pr \left( \sup_{0\leq t<r^{\gamma}T} \left|E_{3}^{r,n}(t) - \mu_{3}^r t\right| > \epsilon_2r^{\gamma}T,\;  \tau_{2n-1}^r \leq r^2T\right)\label{lemma_up_eq_13_10}\\
&\hspace{4cm}+ \sum_{n=1} ^{N^r} \pr \left( \sup_{0\leq t<r^{\gamma}T} \left|E_{A}^{r,n}(t) - \mu_{A}^r t\right| > \epsilon_2r^{\gamma}T,\;  \tau_{2n-1}^r \leq r^2T \right).\label{lemma_up_eq_13_11}
\end{align}
It is straightforward to see that the sums in \eqref{lemma_up_eq_13_10} and \eqref{lemma_up_eq_13_11} converges to 0 as $r\rightarrow\infty$ by the following result, whose proof is presented in Appendix \ref{proof_of_ata_modification}. 

\begin{lemma}\label{ata_modification}
For all $\gamma >0$, $\epsilon_2>0$ such that $\epsilon_2 < \mu_3 \wedge (( \mu_{A}- \mu_{3})/2)$, and $j\in\{3,A\}$,
\begin{equation}\label{ata_m_convergence}
\sum_{n=1} ^{N^r}  \pr \left( \sup_{0\leq t<r^{\gamma}T} \left|E_{j}^{r,n}(t) - \mu_{j}^r t\right| > \epsilon_2r^{\gamma}T,\;  \tau_{2n-1}^r \leq r^2T \right)\rightarrow 0,\qquad\text{as $r\rightarrow\infty$}.
\end{equation}
\end{lemma}
Note that Lemma \ref{ata_modification} extends Lemma 9 of \citet{ata05} to a renewal process that starts from a random time. Consequently, the sum in \eqref{ud_7} converges to $0$ as $r\rightarrow\infty$.

\begin{remark}
The sum in \eqref{ud_7} converges to $0$ for all $\gamma>0$ as $r\rightarrow\infty$. We need $\gamma<1$ in the next section (see \eqref{lemma_up_eq_26}).
\end{remark}

\subsubsection{Proof of Convergence of \eqref{ud_8}}

In this section, we will show that the sum in \eqref{ud_8} converges to $0$ as $r\rightarrow\infty$. The sum in \eqref{ud_8} is less than or equal to
\begin{equation} \label{lemma_up_eq_23}
\sum_{n=1} ^{N^r}\pr \left(  \tau_{2n-1}^r \leq r^2T,\; (\mE_{n,1}^r)^c\right)+\sum_{n=1} ^{N^r}\pr \left( \sup_{0\leq t < r^{\gamma}T} E_{3}^{r,n}(t) - E_{A}^{r,n}(t) >r \epsilon -1,\;  \tau_{2n-1}^r\leq r^2T  ,\; \mE_{n,1}^r \right),
\end{equation}
where $\mE_{n,1}^r$ is defined in \eqref{lemma_up_eq_10_1}. Note that the first sum in \eqref{lemma_up_eq_23} converges to zero with the same way \eqref{lemma_up_eq_13} does. The second sum in \eqref{lemma_up_eq_23} is less than or equal to
\begin{align}
&\sum_{n=1} ^{N^r} \pr \left(\sup_{0\leq t < r^{\gamma}T} E_{3}^{r,n}(t) - E_{A}^{r,n}(t) >r \epsilon -1 ,\; \mE_{n,1}^r\right)\nonumber\\
&\hspace{2cm} \leq \sum_{n=1} ^{N^r} \pr \left(  \sup_{0\leq t<r^{\gamma}T} \mu_{3}^r t+\epsilon_2r^{\gamma}T- (\mu_{A}^r t-\epsilon_2r^{\gamma}T) >r \epsilon -1\right)\label{lemma_up_eq_25_0}\\
&\hspace{2cm} = N^r \pr \left( 2\epsilon_2r^{\gamma}T >r \epsilon -1\right) \label{lemma_up_eq_25}\\
&\hspace{2cm} = N^r \I \left( 2\epsilon_2r^{\gamma}T >r \epsilon -1 \right)\rightarrow 0\qquad\text{as $r\rightarrow\infty$},\label{lemma_up_eq_26}
\end{align}
where \eqref{lemma_up_eq_25_0} is by \eqref{lemma_up_eq_11_1}; \eqref{lemma_up_eq_25} is by the fact that when $r$ is sufficiently large, we have $\mu_{A}^r > \mu_{3}^r$; and \eqref{lemma_up_eq_26} is by the fact that $\gamma<1$. 

Therefore, we prove that \eqref{ud_1} converges to 0 as $r\rightarrow\infty$. In the next section, we will prove that \eqref{ud_2} converges to 0 as $r\rightarrow\infty$.

\subsection{Proof of Convergence of \eqref{ud_2} (Down Intervals)}\label{ud_proof_2}

In this section, we consider the down intervals and prove the convergence of \eqref{ud_2}. Let $N_2^r:=\ru{(\mu_{A}^r+\epsilon_3)r^2T}$ where $\epsilon_3>0$ is an arbitrary constant. Since there is a service completion of a type $a$ job in server 4 at $\tau_{2n}^r$ for all $n\in \N_+$, then as $r\rightarrow\infty$
\begin{equation}\label{lemma_down_eq_1}
\pr \left( \tau_{2N_2^r}^r \leq r^2T\right)\leq \pr \left(S_{A}^r(r^2T) \geq N_2^r\right)\leq \pr \left(S_{A}^r(r^2T) \geq (\mu_{A}^r+\epsilon_3)r^2T \right)\rightarrow 0\qquad\text{as $r\rightarrow\infty$},
\end{equation}
where \eqref{lemma_down_eq_1} is by FSLLN. Let
\begin{equation}\label{lemma_down_eq_1_1}
\Gamma_n^r(t) := \left|Q^r_4(t)- Q^r_3(t)\wedge W^r_4(t) \right|.
\end{equation}
Then, the probability in \eqref{ud_2} is less than or equal to
\begin{equation}\label{lemma_down_eq_2}
\pr \left( \tau_{2N_2^r}^r \leq r^2T\right) + \pr \left( \sup_{t\in [0,r^2T]\cap \bigcup_{n=0}^\infty [\tau_{2n}^r,\tau_{2n+1}^r)} \Gamma_n^r(t) >r \epsilon,\;\tau_{2N_2^r}^r > r^2T \right).
\end{equation}
Note that the first probability in \eqref{lemma_down_eq_2} converges to 0 as $r\rightarrow\infty$ by \eqref{lemma_down_eq_1} and the second probability in \eqref{lemma_down_eq_2} is less than or equal to
\begin{align}
& \sum_{n=0}^{N_2^r} \pr \left( \sup_{\tau_{2n}^r\leq t <\tau_{2n+1}^r\wedge r^2(T+\epsilon_4)}  \Gamma_n^r(t) >r \epsilon,\;\tau_{2n}^r\leq r^2T  \right) \nonumber \\
& \hspace{1cm} = \sum_{n=0}^{N_2^r} \pr \Bigg( \left[ \left(\sup_{\tau_{2n}^r\leq t \leq \tau_{2n}^{r,*}\wedge r^2(T+\epsilon_4)}  \Gamma_n^r(t)\right) \vee \left(\sup_{\tau_{2n}^{r,*}\wedge r^2(T+\epsilon_4) < t <\tau_{2n+1}^r\wedge r^2(T+\epsilon_4)}  \Gamma_n^r(t)\right) \right] \nonumber \\
&\hspace{2cm}  \times \I(\tau_{2n}^{r,*}<\infty) + \left[\sup_{\tau_{2n}^r\leq t <\tau_{2n+1}^r\wedge r^2(T+\epsilon_4)}  \Gamma_n^r(t)\right] \I(\tau_{2n}^{r,*}=\infty) >r \epsilon,\;\tau_{2n}^r\leq r^2T \Bigg),\label{lemma_down_eq_3}
\end{align}
where $\epsilon_4>0$ is an arbitrary constant introduced to cover the time instant $r^2T$. For completeness, we define $\sup_{t\in\emptyset} X(t)=0$ for any $X\in\D$. By definition (cf. \eqref{tau_star}), $\tau_{2n}^{r,*}=\infty$ implies that down 2 interval does not exist within $[\tau_{2n}^r,\tau_{2n+1}^r)$, thus buffer 5 never becomes empty during the corresponding down interval. Hence, server 4 does not work on buffer 4 during the same interval and the down interval ends with the first service completion in server 3. Therefore, when $\tau_{2n}^{r,*}=\infty$, $Q^r_4(t)=Q^r_3(t)$ and $\Gamma_n^r(t)=0$ for all $t\in[\tau_{2n}^r,\tau_{2n+1}^r)$ by \eqref{lemma_down_eq_1_1} and the fact that $Q^r_4(t)\leq W^r_4(t)$ for all $t\geq 0$ (cf. \eqref{workload_process}). Moreover, when $\tau_{2n}^{r,*}<\infty$, by the same logic, we have $Q^r_4(t)=Q^r_3(t)$ and $\Gamma_n^r(t)=0$ for all $t\in[\tau_{2n}^r,\tau_{2n}^{r,*}]$. Therefore, the only nonzero term in \eqref{lemma_down_eq_3} is the one associated with the down 2 interval $(\tau_{2n}^{r,*}, \tau_{2n+1}^r)$, and \eqref{lemma_down_eq_3} is equal to
\begin{equation}\label{lemma_down_eq_3_0}
\sum_{n=0}^{N_2^r} \pr \left( \left[ \sup_{\tau_{2n}^{r,*}\wedge r^2(T+\epsilon_4) < t <\tau_{2n+1}^r\wedge r^2(T+\epsilon_4)}  \Gamma_n^r(t)\right] \I(\tau_{2n}^{r,*}<\infty) >r \epsilon,\;\tau_{2n}^r\leq r^2T \right).
\end{equation}

From the preceding argument, it is enough to show that the term in \eqref{lemma_down_eq_3_0} converges to 0 as $r\rightarrow\infty$. The term in \eqref{lemma_down_eq_3_0} is equal to
\begin{align}
&\sum_{n=0}^{N_2^r} \pr \left( \left[ \sup_{\tau_{2n}^{r,*} < t < \tau_{2n+1}^r\wedge r^2(T+\epsilon_4)}  \Gamma_n^r(t)\right] \I(\tau_{2n}^{r,*}<\infty) >r \epsilon,\;\tau_{2n}^r\leq r^2T,\; \tau_{2n}^{r,*}\leq r^2(T+\epsilon_4)\right)\nonumber\\
&\hspace{1cm} \leq \sum_{n=0}^{N_2^r} \pr \left(  \sup_{\tau_{2n}^{r,*}< t <\tau_{2n+1}^r\wedge r^2(T+\epsilon_4)}  \Gamma_n^r(t)  >r \epsilon,\; \tau_{2n}^{r,*}\leq r^2(T+\epsilon_4) \right)\nonumber\\
& \hspace{1cm}= \sum_{n=0}^{N_2^r} \pr \left(  \sup_{\tau_{2n}^{r,*} < t <\tau_{2n+1}^r\wedge r^2(T+\epsilon_4)} Q^r_3(t)\wedge W^r_4(t)- Q^r_4(t) >r \epsilon,\; \tau_{2n}^{r,*}\leq r^2(T+\epsilon_4) \right)\label{lemma_down_eq_4}\\
&\hspace{1cm} \leq \sum_{n=0}^{N_2^r} \pr \left(  \sup_{\tau_{2n}^{r,*}< t <\tau_{2n+1}^r\wedge r^2(T+\epsilon_4)} \frac{\mu_{A}^r}{\mu_{B}^r}Q^r_5(t)  >r \epsilon,\; \tau_{2n}^{r,*}\leq r^2(T+\epsilon_4) \right),\label{lemma_down_eq_5}
\end{align}
where \eqref{lemma_down_eq_4} is by \eqref{workload_process}, \eqref{lemma_down_eq_1_1}, and the fact that $Q^r_4(t)\leq Q^r_3(t)$ for all $t\in(\tau_{2n}^{r,*},\tau_{2n+1}^r)$ by construction; and \eqref{lemma_down_eq_5} is by \eqref{workload_process}. Let $\epsilon_5:=\inf_{r\geq r_1} (\mu_{B}^r/\mu_{A}^r) \epsilon $ for some $r_1>0$. By Assumption \ref{assumption_rate}, $\epsilon_5>0$ when $r_1$ is sufficiently large. Then, when $r\geq r_1$, the sum in \eqref{lemma_down_eq_5} is less than or equal to
\begin{equation}\label{lemma_down_eq_6}
\sum_{n=0}^{N_2^r} \pr \left(  \sup_{\tau_{2n}^{r,*}< t <\tau_{2n+1}^r\wedge r^2(T+\epsilon_4)} Q^r_5(t) >r \epsilon_5,\; \tau_{2n}^{r,*}\leq r^2(T+\epsilon_4) \right).
\end{equation}
Note that, when $\tau_{2n}^{r,*}<\infty$, $Q^r_5(\tau_{2n}^{r,*})=0$ for all $n\in\N$ by construction.

By Assumption \ref{assumption_rate} Parts \ref{lambda_rate_assumption}, \ref{mu_rate_assumption}, and  \ref{ht_assumption_s4_1}, $\mu_{B}^r>\lambda_b^r$ for r sufficiently large. Since server 4 gives priority to buffer 5 during a down 2 interval, each probability within the sum in \eqref{lemma_down_eq_6} resembles the probability that the supremum of the buffer length in a light traffic $GI/GI/1$ queue is at least $O(r)$ within a down 2 interval. Moreover, there are $O(r^2)$ probabilities in the sum in \eqref{lemma_down_eq_6}. Therefore, the key to the proof is to first construct a specific light traffic $GI/GI/1$ queue and then bound each probability in the sum in \eqref{lemma_down_eq_6} by the probability that the supremum of the buffer length process in the constructed $GI/GI/1$ queue is greater than $r\epsilon_5$ within a compact time interval with length $O(r^2)$. Then we obtain the desired convergence result by Proposition \ref{light_traffic_conv_rate_lemma}.

\subsubsection{Construction of the light traffic $GI/GI/1$ queue}\label{construction_of_light_traffic}

In this section, we construct a light traffic $GI/GI/1$ queue to bound the sum in \eqref{lemma_down_eq_6}. This is doable because by Assumption \ref{assumption_rate} Parts \ref{lambda_rate_assumption}, \ref{mu_rate_assumption}, and  \ref{ht_assumption_s4_1}, $\lambda_b<\mu_{B}$. However, the approach is different depending on if server 2 is in light traffic ($\lambda_b<\mu_2$) or heavy traffic ($\lambda_b=\mu_2$); see the case 1 and case 2 below.

When server 2 is in light traffic, the assumption that server 2 has instantaneous processing capacity allows us to bound $Q_5^r$ over $(\tau_{2n}^{r,*},\tau_{2n+1}^r\wedge r^2(T+\epsilon_4))$ for each $n$ by a light traffic $GI/GI/1$ queue with arrival rate $\lambda_b$, service rate $\mu_{B}$, and initial buffer length $Q_2^r(\tau_{2n}^{r,*})$, which is not too large. The size of the initial buffer length cannot be controlled when server 2 is in heavy traffic. In that case, since $\lambda_b=\mu_2<\mu_{B}$, the light traffic $GI/GI/1$ queue bound is obtained by pretending there is an infinite number of jobs in front of server 2.(The reason why this approach does not work when server 2 is in light traffic is that $\mu_2\geq\mu_{B}$ is possible.)

\paragraph{Case 1: Server 2 is in light traffic} In this section, we prove that the sum in \eqref{lemma_down_eq_6} converges to $0$ as $r\rightarrow\infty$ under the assumption that server 2 is in light traffic ($\lambda_b<\mu_2$). Let, $N_3^r:= \ru{(\lambda_b^r+\mu_{B}^r)r^2(T+\epsilon_4)}$, $\kappa>0$ be an arbitrary constant such that $\kappa<\mu_2\wedge\mu_B-\lambda_b$, and
\begin{equation}\label{lemma_down_eq_7}
\mE_{n,2}^r := \left\{ S_{B}^r(T_{B}^r(\tau_{2n}^{r,*}))\leq N_3^r, \;S_{b}^r(\tau_{2n}^{r,*}) \leq N_3^r, \;Q_2^r(\tau_{2n}^{r,*})\leq \left(\left( \mu_2\wedge\mu_B-\kappa-\lambda_b \right)\wedge \epsilon_5 \right)r/2 \right\}.
\end{equation}
Then the sum in \eqref{lemma_down_eq_6} is less than or equal to
\begin{align}
& \sum_{n=0}^{N_2^r} \pr \left(  \sup_{\tau_{2n}^{r,*} < t <\tau_{2n+1}^r\wedge r^2(T+\epsilon_4)} Q^r_5(t) >r  \epsilon_5,\;\tau_{2n}^{r,*} \leq r^2(T+\epsilon_4),\; \mE_{n,2}^r \right)\label{lemma_down_eq_8}\\
&\hspace{7cm} + \sum_{n=0}^{N_2^r} \pr \left(\tau_{2n}^{r,*} \leq r^2(T+\epsilon_4),\; (\mE_{n,2}^r)^c \right)\label{lemma_down_eq_9}.
\end{align}
The term in \eqref{lemma_down_eq_9} converges to $0$ because $Q_2^r$ is a light traffic $GI/GI/1$ queue, and the number of external type $b$ job arrivals and service completions associated with activity $B$ in an interval of length $O(r^2)$ can be bounded with high probability using the rate of the renewal processes. We have the following result whose proof is presented in Appendix \ref{down_lemma_proof}.

\begin{lemma}\label{down_lemma}
As $r\rightarrow\infty$, $\sum_{n=0}^{N_2^r} \pr \left(\tau_{2n}^{r,*} \leq r^2(T+\epsilon_4),\; (\mE_{n,2}^r)^c \right)\rightarrow 0$.
\end{lemma}

To complete the proof of Proposition \ref{ud_result} when server 2 is in light traffic, it remains only to show the sum in \eqref{lemma_down_eq_8} converges to 0 as $r\rightarrow\infty$. To do this, we first couple $Q_5^r$ with a hypothetical single server queue in which type $b$ jobs are given priority even beyond time $\tau_{2n+1}^r$, so that the random time at the end of the interval $(\tau_{2n}^{r,*},\tau_{2n+1}^r\wedge r^2(T+\epsilon_4))$ can be ignored. Next, we assume that the processing at server 2 is instantaneous, and that the first arrival which occurs after time $\tau_{2n}^{r,*}$ occurs instantaneously, in order to create an upper bound light traffic $GI/GI/1$ queue with arrival rate $\lambda_b^r$ and service rate $\mu_{B}^r$. Finally, we show that the upper bound light traffic $GI/GI/1$ queue converges to 0 as $r\rightarrow\infty$ by Proposition \ref{light_traffic_conv_rate_lemma}.

By sample-path arguments, we will construct  an upper bound (depending on $r$) to the sum in \eqref{lemma_down_eq_8} and show that this upper bound converges to 0 as $r\rightarrow\infty$. Let us define a new hypothetical queueing system in which we use the same control until $\tau_{2n+1}^r$ but after this time epoch, server $4$ always gives static preemptive priority to buffer 5 over buffer 4. We call this hypothetical system as ``system 1''. Let the queue length process for buffer $k$, $k\in\mK$ be $Q^{(1),r}_k$ in system 1. Then, clearly, for each $\omega\in\Omega$, $Q^{(1),r}_5(\omega_t)=Q^{r}_5(\omega_t)$ for all $t< \tau_{2n+1}^r(\omega)$ and the sum in \eqref{lemma_down_eq_8} is less than or equal to 
\begin{align}
&\sum_{n=0}^{N_2^r} \pr \left(  \sup_{\tau_{2n}^{r,*} < t < r^2(T+\epsilon_4)} Q^{(1),r}_5(t)>r  \epsilon_5,\;\tau_{2n}^{r,*}\leq r^2(T+\epsilon_4) ,\;\mE_{n,2}^r \right)\nonumber\\
& \hspace{1cm}=\sum_{n=0}^{N_2^r} \sum_{i=0}^{N_3^r} \sum_{j=0}^{N_3^r} \pr \Bigg(  \sup_{\tau_{2n}^{r,*} < t < r^2(T+\epsilon_4)} Q^{(1),r}_5(t)>r  \epsilon_5,\;\tau_{2n}^{r,*}\leq r^2(T+\epsilon_4), \nonumber\\
& \hspace{3cm} S_{B}^r(T_{B}^r(\tau_{2n}^{r,*}))=i,\; S_{b}^r(\tau_{2n}^{r,*})=j,\; Q_2^{(1),r}(\tau_{2n}^{r,*}) \leq \frac{\left(\left( \mu_2\wedge\mu_B-\kappa-\lambda_b \right)\wedge \epsilon_5 \right)r}{2} \Bigg),\label{lemma_down_eq_11}
\end{align}
where \eqref{lemma_down_eq_11} is by \eqref{lemma_down_eq_7}.

Next, let us construct another hypothetical queueing system, which we call ``system 2'', by modifying system 1 in the following way. Let the queue length process for buffer $k$, $k\in\mK$ be $Q^{(2),r}_k$ in system 2. For each $\omega\in\Omega$ such that $\tau_{2n}^{r,*}(\omega)\leq r^2(T+\epsilon_4)$, server 2 can process jobs instantaneously after $\tau_{2n}^{r,*}(\omega)$. Furthermore, assume that the next arrival after time $\tau_{2n}^{r,*}$, which occurs at time $V_b^r(S_{b}^r(\tau_{2n}^{r,*})+1)$ occurs instead at time $\tau_{2n}^{r,*}$. This implies that when $\tau_{2n}^{r,*}(\omega)< r^2(T+\epsilon_4)$, $Q^{(2),r}_5(\omega_t)=Q^{(1),r}_5(\omega_t)$ for all $t< \tau_{2n}^{r,*}(\omega)$; and $Q^{(2),r}_5(\omega_t) \geq Q^{(1),r}_5(\omega_t)$ for all $t\geq \tau_{2n}^{r,*}(\omega)$. Note that, $Q^{(1),r}_5(\tau_{2n}^{r,*}(\omega))=0$ on the set $\{\tau_{2n}^{r,*}<\infty\}$ and $Q^{(2),r}_5(\tau_{2n}^{r,*}(\omega))= Q^{(1),r}_2(\tau_{2n}^{r,*}(\omega))+1$ on the same set because server 2 depletes buffer 2 instantaneously starting from $\tau_{2n}^{r,*}(\omega)$ and the next arrival occurs immediately in the system 2.

In summary, after $\tau_{2n}^{r,*}$, buffer 5 behaves like a $GI/GI/1$ queue with initial buffer length less than or equal to $\left(\left( \mu_2\wedge\mu_B-\kappa-\lambda_b \right)\wedge \epsilon_5 \right)r/2+1$, service times
\begin{equation}\label{lemma_down_eq_13}
\left\{ v_B^r(S_{B}^r(T_{B}^r(\tau_{2n}^{r,*}))+k),\; k\in\N_+ \right\},
\end{equation}
and interarrival times
\begin{equation}\label{lemma_down_eq_14}
\left\{  v_b^r(S_{b}^r(\tau_{2n}^{r,*})+1+l),\; l\in\N_+ \right\}.
\end{equation}
Since there is a service completion related to activity $B$ at time $\tau_{2n}^{r,*}$, the first service time in \eqref{lemma_down_eq_13} is $v_B^r(S_{B}^r(T_{B}^r(\tau_{2n}^{r,*}))+1$. Then the sum in \eqref{lemma_down_eq_11} is less than or equal to 
\begin{equation}\label{lemma_down_eq_12}
\sum_{n=0}^{N_2^r} \sum_{i=0}^{N_3^r} \sum_{j=0}^{N_3^r} \pr \Bigg(  \sup_{\tau_{2n}^{r,*} < t < r^2(T+\epsilon_4)} Q^{(2),r}_5(t)>r  \epsilon_5,\;\tau_{2n}^{r,*}\leq r^2(T+\epsilon_4),S_{B}^r(T_{B}^r(\tau_{2n}^{r,*}))=i, S_{b}^r(\tau_{2n}^{r,*})=j \Bigg),
\end{equation}
given that $Q_5^{(2),r}(\tau_{2n}^{r,*}) \leq \left(\left( \mu_2\wedge\mu_B-\kappa-\lambda_b \right)\wedge \epsilon_5 \right)r/2+1$.

At this point we use the i.i.d. property of the service and interarrival times. For all $i,j\in \{0,1,\ldots,N_3^r\}$, the interarrival times are equal to $ \left\{ v_b^r(l),\; l\in\N_+ \right\}$ in distribution, and service times are equal to $\left\{ v_B^r(k),\; k\in\N_+ \right\}$  in distribution by \eqref{lemma_down_eq_13} and \eqref{lemma_down_eq_14} in \eqref{lemma_down_eq_12}. Then, let us construct a hypothetical $GI/GI/1$ queue, where the buffer length process is denoted by $Q^{(3),r}_5$, the interarrival and service time sequences are $ \left\{ v_b^r(l),\; l\in\N_+ \right\}$ and $\left\{ v_B^r(k),\; k\in\N_+ \right\}$, respectively, and $Q^{(3),r}_5(0)=\left(\left( \mu_2\wedge\mu_B-\kappa-\lambda_b \right)\wedge \epsilon_5 \right)r/2+1$. It is possible to see that the sum in \eqref{lemma_down_eq_12} is less than or equal to
\begin{equation}\label{lemma_down_eq_15}
N_2^r \left(N_3^r\right)^2  \pr \left(  \sup_{0\leq t \leq r^2(T+\epsilon_4)} Q^{(3),r}_5(t)>r  \epsilon_5 \right),
\end{equation}
given that $Q^{(3),r}_5(0)= \left(\left( \mu_2\wedge\mu_B-\kappa-\lambda_b \right)\wedge \epsilon_5 \right)r/2+1$. 

Lastly, let us construct a hypothetical $GI/GI/1$ queue, which we call ``system 4'' and  where the buffer length process is denoted by $Q^{(4)}_5$. In system 4, the interarrival and service time processes are $\left\{ v_b^r(l)(\lambda_b^r/(\lambda_b+\kappa/2)),\; l\in\N_+ \right\}$ and $\left\{ v^r_B(k)(\mu_B^r/(\mu_B-\kappa/2)),\; k\in\N_+ \right\}$, respectively\footnote{Since $v_b^r(l)(\lambda_b^r/(\lambda_b+\kappa/2))= \bar{v}_b(l)/(\lambda_b+\kappa/2)$ and $v_B^r(l)(\mu_B^r/(\mu_B-\kappa/2))= \bar{v}_B(l)/(\mu_B-\kappa/2)$ for all $l\in\N_+$ (cf. Section \ref{seq_of_systems}), the corresponding sequences are independent of $r$}, and $Q^{(4)}_5(0)=\left(\left( \mu_2\wedge\mu_B-\kappa-\lambda_b \right)\wedge \epsilon_5 \right)r/2+1$. Hence, the arrival and service rates in the hypothetical $GI/GI/1$ queue in system 4 are $\lambda_b+\kappa/2$ and $\mu_B-\kappa/2$, respectively. By Assumption \ref{assumption_rate} Parts \ref{lambda_rate_assumption} and \ref{mu_rate_assumption}, there exists an $r_2\geq 1$ such that if $r\geq r_2$, the term in \eqref{lemma_down_eq_15} is less than or equal to 
\begin{equation}\label{lemma_down_eq_15_1}
N_2^r \left(N_3^r\right)^2  \pr \left(  \sup_{0\leq t \leq r^2(T+\epsilon_4)} Q^{(4)}_5(t)>r  \epsilon_5\right),
\end{equation} 
given that $Q^{(4)}_5(0)= \left(\left( \mu_2\wedge\mu_B-\kappa-\lambda_b \right)\wedge \epsilon_5 \right)r/2+1$. The term in \eqref{lemma_down_eq_15_1} converges to 0 as $r\rightarrow\infty$ by the fact the corresponding server is in light traffic and Proposition \ref{light_traffic_conv_rate_lemma}.

\paragraph{Case 2: Server 2 is in heavy traffic} In this section, we prove that the sum in \eqref{lemma_down_eq_6} converges to $0$ as $r\rightarrow\infty$ under the assumption that server 2 is in heavy traffic ($\lambda_b=\mu_2$). Now, Lemma \ref{down_lemma} does not hold because we cannot show that the event $\{Q_2^r(\tau_{2n}^{r,*})\leq  \epsilon_6 r,\;\tau_{2n}^{r,*} \leq r^2(T+\epsilon_4) \}$ occur with high probability for any $\epsilon_6>0$ in this case. However, we can solve this problem by modifying the part of the proofs starting from \eqref{lemma_down_eq_7} in the following way. First, we modify \eqref{lemma_down_eq_7} by
\begin{equation}\label{lemma_down_eq_7_0}
\mE_{n,2}^r := \left\{ S_{B}^r(T_{B}^r(\tau_{2n}^{r,*}))\leq N_3^r \right\}.
\end{equation}
Then, again, the sum in \eqref{lemma_down_eq_6} is less than or equal to the sum of the terms in \eqref{lemma_down_eq_8} and \eqref{lemma_down_eq_9}, and \eqref{lemma_down_eq_9} converges to 0 as $r\rightarrow\infty$ by the convergence of the sum in \eqref{lemma_down_eq_9_1} to 0 (cf. proof of Lemma \ref{down_lemma} in Appendix \ref{down_lemma_proof}, specifically see \eqref{lemma_down_eq_10}).

We again use sample-path arguments to construct  an upper bound (depending on $r$) to the sum in \eqref{lemma_down_eq_8}. First, let us construct the hypothetical system 1 as before. By \eqref{queue_length_1_2} and the fact that $Q^{(1),r}_5(\tau_{2n}^{r,*})=0$ on the set $\{\tau_{2n}^{r,*}<\infty\}$,
\begin{equation}\label{lemma_down_eq_6_1}
S_{B}^r(T_{B}^r(\tau_{2n}^{r,*})) = S_2^r(T_2^r(\tau_{2n}^{r,*})),
\end{equation}
on the set $\{\tau_{2n}^{r,*}<\infty\}$. By \eqref{lemma_down_eq_7_0} and \eqref{lemma_down_eq_6_1}, the sum in \eqref{lemma_down_eq_8} is less than or equal to
\begin{equation}\label{lemma_down_eq_16}
\sum_{n=0}^{N_2^r} \sum_{i=0}^{N_3^r}  \pr \Bigg(  \sup_{\tau_{2n}^{r,*} < t < r^2(T+\epsilon_4)} Q^{(1),r}_5(t)>r  \epsilon_5,\;\tau_{2n}^{r,*}\leq r^2(T+\epsilon_4),\;S_{B}^r(T_{B}^r(\tau_{2n}^{r,*}))=S_{2}^r(T_{2}^r(\tau_{2n}^{r,*}))=i\Bigg),
\end{equation}
given that $Q^{(1),r}_5(\tau_{2n}^{r,*})=0$. Next, let us construct the hypothetical system 2 in the following way. We assume that at time $\tau_{2n}^{r,*}$, all remaining type $b$ jobs arrive immediately to buffer 2. Then, server 2 never becomes idle after time $\tau_{2n}^{r,*}$. Furthermore, assume that the next service completion in server 2 after time $\tau_{2n}^{r,*}$, which occurs at time $V_2^r(S_{2}^r(T_2^r(\tau_{2n}^{r,*}))+1)$ occurs instead at time $\tau_{2n}^{r,*}$. This implies that when $\tau_{2n}^{r,*}(\omega)< r^2(T+\epsilon_4)$, $Q^{(2),r}_5(\omega_t)=Q^{(1),r}_5(\omega_t)$ for all $t< \tau_{2n}^{r,*}(\omega)$; and $Q^{(2),r}_5(\omega_t) \geq Q^{(1),r}_5(\omega_t)$ for all $t\geq \tau_{2n}^{r,*}(\omega)$. Note that, $Q^{(1),r}_5(\tau_{2n}^{r,*}(\omega))=0$ on the set $\{\tau_{2n}^{r,*}<\infty\}$ and $Q^{(2),r}_5(\tau_{2n}^{r,*}(\omega))= 1$ on the same set because the next service completion in server 2 after $\tau_{2n}^{r,*}$ occurs immediately in the system 2.

In summary, after $\tau_{2n}^{r,*}$, buffer 5 behaves like a $GI/GI/1$ queue with initial buffer level $1$, service times
\begin{equation}\label{lemma_down_eq_17}
\left\{ v_B^r(S_{B}^r(T_{B}^r(\tau_{2n}^{r,*}))+k),\; k\in\N_+ \right\},
\end{equation}
and interarrival times
\begin{equation}\label{lemma_down_eq_18}
\left\{  v_2^r(S_{2}^r(T_2^r(\tau_{2n}^{r,*}))+1+l),\; l\in\N_+ \right\}.
\end{equation}
Then the sum in \eqref{lemma_down_eq_16} is less than or equal to
\begin{equation}\label{lemma_down_eq_19}
\sum_{n=0}^{N_2^r} \sum_{i=0}^{N_3^r}  \pr \Bigg(  \sup_{\tau_{2n}^{r,*} < t < r^2(T+\epsilon_4)} Q^{(2),r}_5(t)>r  \epsilon_5,\;\tau_{2n}^{r,*}\leq r^2(T+\epsilon_4),\;S_{B}^r(T_{B}^r(\tau_{2n}^{r,*}))=S_{2}^r(T_{2}^r(\tau_{2n}^{r,*}))=i \Bigg),
\end{equation}
given that $Q^{(2),r}_5(\tau_{2n}^{r,*})=1$.

At this point we use the i.i.d. property of the service and interarrival times. For all $i\in \{0,1,\ldots,N_3^r\}$ the interarrival times are equal to $\left\{ v_2^r(l),\; l\in\N_+ \right\}$ in distribution, and service times are equal to $\left\{ v_B^r(k),\; k\in\N_+ \right\}$ in distribution by \eqref{lemma_down_eq_17} and \eqref{lemma_down_eq_18} in \eqref{lemma_down_eq_19}. Then, let us construct a hypothetical $GI/GI/1$ queue, where the buffer length process is denoted by $Q^{(3),r}_5$, the interarrival and service time sequences are $ \left\{ v_2^r(l),\; l\in\N_+ \right\}$ and $\left\{ v_B^r(k),\; k\in\N_+ \right\}$, respectively, and $Q^{(3),r}_5(0)=1$. Then, it is possible to see that the sum in \eqref{lemma_down_eq_19} is less than or equal to
\begin{equation}\label{lemma_down_eq_20}
N_2^r N_3^r  \pr \left(  \sup_{0\leq t \leq r^2(T+\epsilon_4)} Q^{(3),r}_5(t)>r  \epsilon_5 \right),
\end{equation}
given that $Q^{(3),r}_5(0)= 1$. 

Lastly, let us construct a hypothetical $GI/GI/1$ queue, which we call ``system 4'' and  where the buffer length process is denoted by $Q^{(4)}_5$. Let us modify the definition of $\kappa>0$ such that $\kappa<\mu_B-\mu_2$. In system 4, the interarrival and service time processes are $ \left\{ v_2^r(l)(\mu_2^r/(\mu_2+\kappa/2)),\; l\in\N_+ \right\}$ and $\left\{ v^r_B(k)(\mu_B^r/(\mu_B-\kappa/2)),\; k\in\N_+ \right\}$, respectively, and $Q^{(4)}_5(0)=1$. Hence, the arrival and service rates in the hypothetical $GI/GI/1$ queue in system 4 are $\mu_2+\kappa/2$ and $\mu_B-\kappa/2$, respectively. By Assumption \ref{assumption_rate} Parts \ref{lambda_rate_assumption} and \ref{mu_rate_assumption}, there exists an $r_3\geq 1$ such that if $r\geq r_3$, the term in \eqref{lemma_down_eq_20} is less than or equal to 
\begin{equation}\label{lemma_down_eq_20_1}
N_2^r \left(N_3^r\right)^2  \pr \left(  \sup_{0\leq t \leq r^2(T+\epsilon_4)} Q^{(4)}_5(t)>r  \epsilon_5\right),
\end{equation} 
given that $Q^{(4)}_5(0)= 1$. The term in \eqref{lemma_down_eq_20_1} converges to 0 as $r\rightarrow\infty$ by the fact the corresponding server is in light traffic and Proposition \ref{light_traffic_conv_rate_lemma}.

\begin{remark}\label{initial_buffer_remark}
We assume that $Q_k^r(0)=0$ for all $k\in\mK$ up to now. However, it is straightforward to see that the results presented so far hold under the following weaker assumption.

\begin{assumption}\label{assumption_initial} For each $r$, $\left(Q_k^r(0), k\in\mK\right)$ is a nonnegative random vector which takes values in $\N^{10}$ such that 
\begin{enumerate}
\item \label{initial_queue_1} $Q_3^r(0)=Q_4^r(0)$, $Q_7^r(0)=Q_8^r(0)$, $Q_5^r(0)=Q_6^r(0)$, and $Q_9^r(0)=Q_{10}^r(0)$ for each $r$.
\item \label{initial_queue_2} $r^{-1} Q_k^r(0)\Longrightarrow 0$ and $r^{-2} Q_k^r(0)\xrightarrow{a.s.} 0$ for all $k\in\mK$ as $r\rightarrow\infty$,
\item \label{initial_buffer_2_rate} If server 2 is in light traffic, i.e., $\lambda_b<\mu_2$, then for each $\epsilon>0$, there exists a $\kappa\in(0,\mu_2\wedge\mu_B-\lambda_b)$ and $r_0\geq 1$ such that if $r\geq r_0$, 
\begin{equation*}
\pr \left(Q_2^r(0)\geq \left(\left( \mu_2\wedge\mu_B-\kappa-\lambda_b \right)\wedge \epsilon \right)r/2\right)\leq C_0r^3\e^{-C_1r},
\end{equation*}
for some constants $C_0>0$ and $C_1>0$ which are independent of $r$.
\end{enumerate}
\end{assumption}
Assumption \ref{assumption_initial} Part \eqref{initial_queue_1} guarantees that the initial buffer lengths do not violate \eqref{balance_eq} at time $t=0$. We need Assumption \ref{assumption_initial} Part \eqref{initial_queue_2} because of two main reasons. First, since we consider cases in which some of the servers are in light traffic, in order to obtain the weak convergence of the diffusion scaled queue length processes corresponding to these servers to $\textbf{0}$, we need the first convergence result in Assumption \ref{assumption_initial} Part \eqref{initial_queue_2}. Second, in order to obtain the a.s. convergence of the fluid scaled queue length processes in each buffer to $\textbf{0}$ u.o.c., we need the second convergence result in Assumption \ref{assumption_initial} Part \eqref{initial_queue_2}. We need Part \eqref{initial_buffer_2_rate} of this assumption in order to invoke Proposition \ref{light_traffic_conv_rate_lemma} in the proof of Lemma \ref{down_lemma}, specifically in \eqref{lemma_down_eq_10_1}.
\end{remark}

\begin{remark}\label{weak_moment_remark}
It is possible to weaken the exponential moment assumption (cf. Assumption \ref{assumption_moment}) for some of the interarrival and service time processes. First, in Section \ref{length_of_ip_intervals}, we use the exponential moment assumption only for the service time processes associated with activities $3$ and $A$ in the proof of Lemma \ref{ata_modification}. However, we can relax this assumption in the following way. We need $4/(1+\alpha)<\gamma<1$ (cf. proof of Lemma \ref{ata_modification}). Since $\gamma\in(0,1)$ can be chosen arbitrarily close to 1, any $\alpha$ strictly greater than $3$ satisfies $4/(1+\alpha)<\gamma<1$. Since we need $\E[(v_j^r(1))^{2+2\alpha}]<\infty$ for $j\in\{3,A\}$ in the proof of Lemma \ref{ata_modification} (cf. \eqref{lemma_up_eq_16} and \eqref{lemma_up_eq_19}), $\E[(v_j^r(1))^{8+\beta}]<\infty$, $j\in\{3,A\}$ for some $\beta>0$ is sufficient to get the convergence result associated with \eqref{ud_1}. Second, in Section \ref{construction_of_light_traffic}, we have used the exponential moment assumption for the interarrival times of type $b$ jobs (cf.  \eqref{lemma_down_eq_15_1} and \eqref{lemma_down_eq_10_1}, which is in the proof of Lemma \ref{down_lemma}), service times associated with activity $2$ (cf. \eqref{lemma_down_eq_20_1} and \eqref{lemma_down_eq_10_1}), and service times associated with activity $B$ (cf. \eqref{lemma_down_eq_15_1} and \eqref{lemma_down_eq_20_1}). For the interarrival times of type $a$ jobs and the service times associated with activities $1$, $5$, $6$, and $7$, we only need finite second moment assumption to use the FCLT in Sections \ref{diffusion_results_1} and \ref{diffusion_results_2}.
\end{remark}


\section{Simulation}\label{numerical_analyses}

In this Section, we use discrete-event simulation to test the performance of the non-preemptive version of the proposed policy which is described in Remark \ref{nonpreemtive_policy_remark}. We consider $36$ different test instances and at each instance, we compare the performance of the non-preemptive version of the proposed policy with the performances of $4$ other non-preemptive control policies. The instances are designed to consider various cases related to the processing capacities of the servers and variability in the service times. We describe the simulation setup in Section \ref{sim_setup}, then we present the results of the experiments in Section \ref{numerical_results}.

\subsection{Simulation Setup}\label{sim_setup}

At each instance, the arrival processes of type $a$ and $b$ jobs are independent Poisson processes with rate one, thus $\lambda_a=\lambda_b=1$\footnote{For convenience in notation, we have dropped the superscript $r$ on the parameters. The reader should understand that each set of parameters $(\lambda_a,\lambda_b,\mu_1,\mu_2,\ldots,\mu_7)$ is associated with a particular $r$, and is not the limit parameter that appear in Assumption \ref{assumption_rate}.}. At each instance, servers $5$, $6$, and $7$ are in either heavy or light traffic, where heavy and light traffic mean $95\%$ and $70\%$ long-run utilization rates of the corresponding server, respectively. For example, when server 6 is in heavy (light) traffic at an instance, then $\mu_6=1/0.95$ ($\mu_6=1/0.7$), which gives the desired long-run utilization rate. At each instance, servers $1$ and $2$ are in light traffic and server 4 is in heavy traffic such that $\mu_A=\mu_B=2/0.95$. At each instance, we assume that $\mu_3\in\{1/0.95,1/0.7125,1/0.35\}$ to test the performance of the proposed policy when $\lambda_a\approx \mu_3$, $\lambda_a<\mu_3<\mu_A$, and $\mu_3>\mu_A$, respectively. 

We use three different distribution types for the service time processes, namely Erlang-3, Exponential, and Gamma distributions. When the service time process associated with an activity is Gamma distributed, we choose the distribution parameters such that the squared coefficient of variation of the corresponding service time process is $3$. Note that the squared coefficient of variations in Erlang-3 and Exponential distributions are $1/3$ and $1$, respectively. Therefore, Erlang-3, Exponential, and Gamma distributions correspond to the low, moderate, and high level variability in the service time processes, respectively. At each instance, we use the same distribution type for all service time processes. 

Table \ref{test_beds} shows the parameter choices related to the service time processes at each instance. For example, at instance 7, $\mu_3=1/0.7125$, server 5 is in heavy traffic, and servers 6 and 7 are in light traffic. Moreover, the service time processes associated with each activity are Erlang-3 distributed at this instance. In the first 18 instances, the downstream servers, i.e. servers 6 and 7, are in light traffic. Note that we prove the asymptotic optimality of the proposed policy under this assumption (cf. Assumption \ref{assumption_rate} Part \ref{rate_assumption_s6s7}). However, in order to test the robustness of the proposed policy, we consider the cases in which server 6 or 7 is in heavy traffic in the last 18 instances.

\begin{table}[htbp]\scriptsize
\centering \caption{Parameter choices at each instance.}
\resizebox{0.9\textwidth}{!}{
\begin{tabular}{c c c c|c c|c c|c c|c c c}
\cline{2-13}
 & \multicolumn{3}{ c|}{$\bm{\mu_3}$} & \multicolumn{2}{|c|}{$\bm{\mu_5}$} & \multicolumn{2}{|c|}{$\bm{\mu_6}$} & \multicolumn{2}{|c|}{$\bm{\mu_7}$} & \multicolumn{3}{|c}{\textbf{Service Time Variability}}\\
\hline
\textbf{Instance} & $1/0.95$ & $1/0.7125$ & $1/0.35$ & $1/0.95$ & $1/0.7$ & $1/0.95$ & $1/0.7$ & $1/0.95$ & $1/0.7$ & Low & Moderate & High\\
\hline
1 & X &  &  & X &  &  & X &  & X & X &  & \\

2 & X &  &  & X &  &  & X &  & X &  & X & \\

3 & X &  &  & X &  &  & X &  & X &  &  & X\\

4 & X &  &  &  & X &  & X &  & X & X &  & \\

5 & X &  &  &  & X &  & X &  & X &  & X & \\

6 & X &  &  &  & X &  & X &  & X &  &  & X\\
\hline
7 &  & X &  & X &  &  & X &  & X & X &  & \\

8 &  & X &  & X &  &  & X &  & X &  & X & \\

9 &  & X &  & X &  &  & X &  & X &  &  & X\\

10 &  & X &  &  & X &  & X &  & X & X &  & \\

11 &  & X &  &  & X &  & X &  & X &  & X & \\

12 &  & X &  &  & X &  & X &  & X &  &  & X\\
\hline
13 &  &  & X & X &  &  & X &  & X & X &  & \\

14 &  &  & X & X &  &  & X &  & X &  & X & \\

15 &  &  & X & X &  &  & X &  & X &  &  & X\\

16 &  &  & X &  & X &  & X &  & X & X &  & \\

17 &  &  & X &  & X &  & X &  & X &  & X & \\

18 &  &  & X &  & X &  & X &  & X &  &  & X\\
\hline
19 & X &  &  & X &  & X &  &  & X &  & X & \\

20 &  & X &  & X &  & X &  &  & X &  & X & \\

21 &  &  & X & X &  & X &  &  & X &  & X & \\

22 & X &  &  &  & X & X &  &  & X &  & X & \\

23 &  & X &  &  & X & X &  &  & X &  & X & \\

24 &  &  & X &  & X & X &  &  & X &  & X & \\
\hline
25 & X &  &  & X &  & X &  & X &  &  & X & \\

26 &  & X &  & X &  & X &  & X &  &  & X & \\

27 &  &  & X & X &  & X &  & X &  &  & X & \\

28 & X &  &  &  & X & X &  & X &  &  & X & \\

29 &  & X &  &  & X & X &  & X &  &  & X & \\

30 &  &  & X &  & X & X &  & X &  &  & X & \\
\hline
31 & X &  &  & X &  &  & X & X &  &  & X & \\

32 &  & X &  & X &  &  & X & X &  &  & X & \\

33 &  &  & X & X &  &  & X & X &  &  & X & \\

34 & X &  &  &  & X &  & X & X &  &  & X & \\

35 &  & X &  &  & X &  & X & X &  &  & X & \\

36 &  &  & X &  & X &  & X & X &  &  & X & \\
\hline
\end{tabular}}\label{test_beds}
\end{table}

Even though $\oh{W}_4^r$ weakly converges to the same limit independent of the control (cf. Proposition \ref{general_conv}), $\oh{W}_4^r$ is policy dependent in the pre-limit. Since our class of admissible scheduling policies is very large, and includes the ones that can anticipate the future, we cannot construct a pre-limit lower bound based on the solution of DCP \eqref{DCP_4} given in Proposition \ref{lb_DCP}. Hence, there is no pre-limit lower bound on performance. Furthermore, we cannot use the DCP solution to develop an approximate pre-limit lower bound on the average holding cost, because that requires knowing the stationary distribution of the relevant SRBM and that is only straightforward to find under very specific conditions (cf. \citet{har87,dai92,die11,dai14}). The relevant SRBM does not have a product form stationary distribution in our case (cf. \citet{har87}), when it is at least two-dimensional. Therefore, we compare the performances of $5$ different non-preemptive control policies. The first one is the non-preemptive version of the proposed policy (see Remark \ref{nonpreemtive_policy_remark} for the definition). Since the proposed policy is the SDP policy and the static priority policy when $\lambda_a\leq\mu_3<\mu_A$ and $\mu_3\geq\mu_A$, respectively, the second and the third control policies that we consider are the non-preemptive versions of the SDP and the static priority policies, respectively. The fourth policy that we consider is the \textit{FCFS policy}, in which whenever server 4 is ready to process a new job, it chooses the job which has arrived the earliest to the buffers 4 or 5. The fifth policy that we consider is the \textit{randomized policy}, in which whenever server 4 is ready to process a job and if both of the buffers 4 and 5 are non-empty, server 4 chooses a type $a$ job with half probability. If only one of the buffers 4 and 5 is non-empty, then server 4 chooses the job from the non-empty one. 

We have used Omnet$++$ discrete-event simulation freeware in our experiments. At each instance associated with each control policy that we considered, we have done 30 replications. At each replication, we have created approximately $1$ million type $a$ jobs and $1$ million type $b$ jobs and we have considered the time interval in which the first $50,000$ type $a$ jobs and the first $50,000$ type $b$ jobs arrive as the warm-up period.


\subsection{Simulation Results}\label{numerical_results}

In this section, we present the results of the simulation experiments. We assume that $h_b=1$ in all experiments and consider various $h_a$ values such that $h_a\geq h_b$. In each experiment, we use the same $h_a$ value at all instances.

Let $\mP$ denote the set of the five control policies that we consider in the simulation experiments and $Q_k^{p,j}(i)$ denote the average length of buffer $k$, $k\in\mK$ in replication $j$, $j\in\{1,2,\ldots,30\}$ with respect to policy $p$, $p\in\mP$ at instance $i$, $i\in\{1,2,\ldots,36\}$. Then, $\ol{Q}_k^p(i):= (1/30) \sum_{j=1}^{30} Q_k^{p,j}(i)$ is the average length of buffer $k$ corresponding to policy $p$ at instance $i$. Note that, we have proved the asymptotic optimality of the proposed policy with respect to the objectives \eqref{obj_1}, \eqref{obj_2}, and \eqref{obj_5} but not with respect to the average cost objective \eqref{obj_3} (cf. Remark \ref{obj_optimality}). Still, our results suggest that our proposed policy performs well with respect to \eqref{obj_3}, which is a natural objective to consider in simulation experiments. Hence, we use the objective \eqref{obj_3}. Recalling the equality \eqref{balance_eq}, we only need to compute $J_p(i) := h_a\big(\ol{Q}_3^p(i)+\ol{Q}_7^p(i)\big)+h_b\big(\ol{Q}_6^p(i)+\ol{Q}_{10}^p(i)\big)$, which is the average total holding cost in buffers 3, 6, 7, and 10 per unit time corresponding to policy $p$ at instance $i$. Let $L(i):=\min_{p\in\mP} J_p(i)$ denote the lowest average cost among the policies in $\mP$ at instance $i$.

We present the detailed results of the simulation experiments in Table \ref{detailed_results}, which is in Appendix \ref{detailed_numerical_results}. This table contains $\ol{Q}_3^p(i)+\ol{Q}_7^p(i)$ and $\ol{Q}_6^p(i)+\ol{Q}_{10}^p(i)$ for each $i$ and $p$ with their $95\%$ confidence intervals. 

Table \ref{aggregate_results} shows the average and maximum deviations of the cost of the policies from the lowest realized average costs among the first and last 18 instances for different $h_a$ values. For example, for given $h_a$, the ``Avg.'' and ``Max.'' columns corresponding to policy $p$ and the first 18 instances are
\begin{equation*}
\frac{1}{18} \sum_{i=1}^{18} \left( 100\times\frac{J_p(i)-L(i)}{L(i)}\right),\qquad \max_{i\in\{1,2,\ldots,18\}} \left\{100\times\frac{J_p(i)-L(i)}{L(i)}\right\},
\end{equation*}
respectively. In the ``All'' row, ``Avg.'' (``Max.'') column corresponding to policy $p$ denotes the average (maximum) of the values in the  ``Avg.'' (``Max.'') column among all $h_a$ values.

\begin{table}[htbp]\scriptsize
\centering \caption{Average and maximum deviations of the cost of the policies from the lowest realized average cost.}
\resizebox{0.8\textwidth}{!}{
\begin{tabular}{ c c c|c c|c c|c c|c c }
\hline
 & \multicolumn{10}{ c }{\textbf{First 18 instances}}\\
\cline{2-11}
 & \multicolumn{2}{ c|}{\textbf{Proposed}} & \multicolumn{2}{ c |}{\textbf{SDP}} & \multicolumn{2}{c|}{\textbf{Static}} & \multicolumn{2}{c|}{\textbf{FCFS}} & \multicolumn{2}{c}{\textbf{Randomized}} \\
\hline
$\bm{h_a}$ & \textbf{Avg.} & \textbf{Max.} & \textbf{Avg.} & \textbf{Max.} & \textbf{Avg.} & \textbf{Max.} & \textbf{Avg.} & \textbf{Max.} & \textbf{Avg.} & \textbf{Max.}  \\
\hline
$\bm{1}$ & $1.6\%$ & $6.4\%$ & $1.8\%$ & $6.5\%$ & $9.1\%$ & $30.7\%$ & $3.7\%$ & $9.5\%$ & $4.1\%$ & $8.7\%$ \\

$\bm{1.25}$ & $0.0\%$ & $0.3\%$ & $0.3\%$ & $1.4\%$ & $6.2\%$ & $25.8\%$ & $6.2\%$ & $15.3\%$ & $6.7\%$ & $14.6\%$ \\

$\bm{1.5}$ & $0.1\%$ & $0.8\%$ & $0.5\%$ & $1.7\%$ & $5.3\%$ & $22.2\%$ & $10.0\%$ & $20.7\%$ & $10.6\%$ & $20.2\%$ \\

$\bm{1.75}$ & $0.2\%$ & $1.4\%$ & $0.7\%$ & $2.0\%$ & $4.6\%$ & $19.3\%$ & $13.6\%$ & $25.9\%$ & $14.2\%$ & $25.4\%$ \\

$\bm{2}$ & $0.3\%$ & $2.0\%$ & $0.8\%$ & $2.3\%$ & $4.0\%$ & $17.1\%$ & $16.9\%$ & $30.8\%$ & $17.5\%$ & $30.4\%$ \\

$\bm{3}$ & $0.8\%$ & $4.1\%$ & $1.6\%$ & $4.1\%$ & $2.5\%$ & $11.3\%$ & $28.5\%$ & $48.0\%$ & $29.2\%$ & $48.0\%$\\

$\bm{4}$ & $1.3\%$ & $5.8\%$ & $2.5\%$ & $5.8\%$ & $1.8\%$ & $8.4\%$ & $37.8\%$ & $65.3\%$ & $38.6\%$ & $65.5\%$ \\

$\bm{10}$ & $3.1\%$ & $11.6\%$ & $5.2\%$ & $11.6\%$ & $0.4\%$ & $2.8\%$ & $68.9\%$ & $125.1\%$ & $70.1\%$ & $123.3\%$ \\
\hline
\textbf{All} & $0.9\%$ & $11.6\%$ & $1.7\%$ & $11.6\%$ & $4.2\%$ & $30.7\%$ & $23.2\%$ & $125.1\%$ & $23.9\%$ & $123.3\%$ \\
\hline
 & \multicolumn{10}{ c }{\textbf{Last 18 instances}}\\
\hline
$\bm{1}$ & $4.1\%$ & $17.0\%$ & $4.1\%$ & $16.1\%$ & $9.5\%$ & $25.6\%$ & $0.8\%$ & $5.2\%$ & $1.6\%$ & $6.5\%$ \\

$\bm{1.25}$ & $2.3\%$ & $11.4\%$ & $2.4\%$ & $11.4\%$ & $6.9\%$ & $22.1\%$ & $1.6\%$ & $8.4\%$ & $2.4\%$ & $10.1\%$ \\

$\bm{1.5}$ & $1.4\%$ & $7.4\%$ & $1.5\%$ & $7.4\%$ & $5.3\%$ & $19.3\%$ & $2.9\%$ & $11.1\%$ & $3.8\%$ & $13.0\%$ \\

$\bm{1.75}$ & $0.7\%$ & $3.8\%$ & $0.8\%$ & $3.8\%$ & $4.0\%$ & $17.3\%$ & $4.1\%$ & $13.3\%$ & $5.0\%$ & $15.5\%$ \\

$\bm{2}$ & $0.2\%$ & $1.4\%$ & $0.4\%$ & $1.4\%$ & $3.2\%$ & $15.8\%$ & $5.3\%$ & $15.3\%$ & $6.2\%$ & $17.7\%$ \\

$\bm{3}$ & $0.4\%$ & $2.0\%$ & $0.7\%$ & $2.0\%$ & $2.1\%$ & $11.5\%$ & $10.9\%$ & $24.8\%$ & $11.8\%$ & $26.1\%$ \\

$\bm{4}$ & $0.6\%$ & $2.5\%$ & $1.1\%$ & $2.5\%$ & $1.6\%$ & $8.8\%$ & $15.2\%$ & $35.9\%$ & $16.2\%$ & $37.3\%$ \\

$\bm{10}$ & $1.4\%$ & $7.0\%$ & $2.3\%$ & $7.0\%$ & $0.4\%$ & $2.5\%$ & $29.2\%$ & $81.9\%$ & $30.3\%$ & $84.0\%$ \\
\hline
\textbf{All} & $1.4\%$ & $17.0\%$ & $1.7\%$ & $16.1\%$ & $4.1\%$ & $25.6\%$ & $8.7\%$ & $81.9\%$ & $9.7\%$ & $84.0\%$  \\
\hline
\end{tabular}}\label{aggregate_results}
\end{table}

According to the results, the proposed policy in general performs the best with respect to both the average and maximum deviations from the lowest realized average cost. The SDP policy performs very close to the proposed policy, which shows that it performs well even when $\mu_3>\mu_A$. As expected, the performance of the static priority policy becomes better as $h_a$ increases, and it performs the best when $h_a$ is much larger than $h_b$. Performances of the FCFS and the randomized policies become worse as $h_a$ increases, which is expected because these policies do not give more priority to type $a$ jobs than type $b$ jobs. On average, these two policies perform the worst with respect to both the average and maximum deviations from the lowest realized average cost. 

At the first 18 instances, when $h_a=10$, the static priority and the proposed policies perform the best and the second best, respectively. This is not surprising given that in all of the instances in which server 3 is in light traffic, both the static priority and the SDP policies are asymptotically optimal (cf. Remarks \ref{multiple_policy} and \ref{sp_optimality_range}). The superior performance of the static priority policy in the pre-limit can be attributed to the high holding cost of type $a$ jobs.

At the last 18 instances, the proposed policy still performs the best in average, which suggests that the performance of the proposed policy is robust with respect to the processing capacities of the downstream servers.

Another interesting result that we see from Table \ref{aggregate_results} is that the percentage deviation of the worst performing policy from the lowest realized average cost at the first 18 instances is much higher than the same deviation at the last 18 instances. This result is due to the fact that at least one of the downstream servers is in heavy traffic at the last 18 instances and this decreases the waiting time of the jobs for the ones that they are going to be matched in the join servers. In other words, synchronization requirements between the jobs in different buffers become less important in this case. Hence, the effect of the control policy on the system performance is less important when at least one of the downstream servers is in heavy traffic, which justifies Assumption \ref{assumption_rate} Part \ref{rate_assumption_s6s7}. Moreover, FCFS policy performs the best when $h_a=1$ at the last 18 instances because all jobs are equally expensive and this policy approximately matches the arrival time of the same type of jobs to the downstream buffers.

\begin{figure}[htbp]
  \centering
  \subfloat[Proposed, SDP, and static policies.]{\label{Figure_r1}\includegraphics[width=0.5\textwidth]{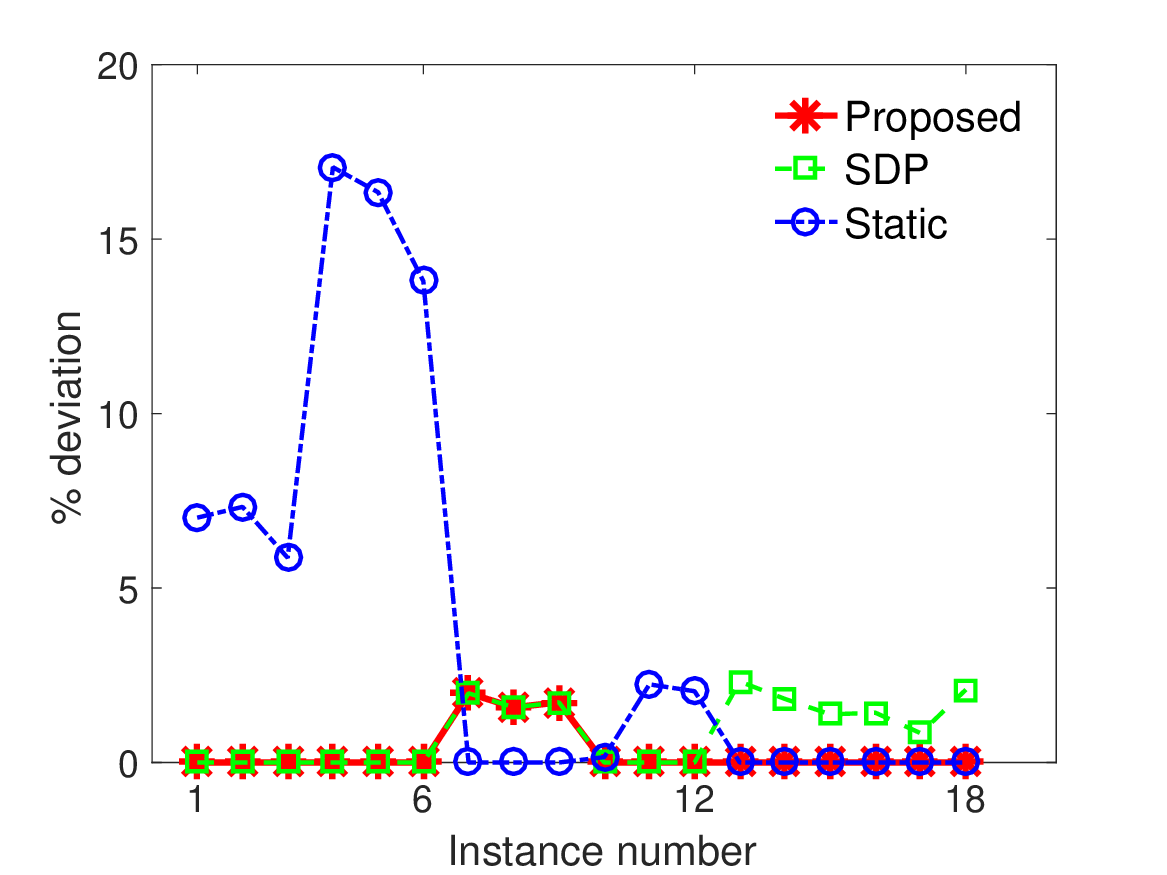}}
  \subfloat[Proposed, FCFS, rand., and rand.-$2/3$ policies.]{\label{Figure_r2}\includegraphics[width=0.5\textwidth]{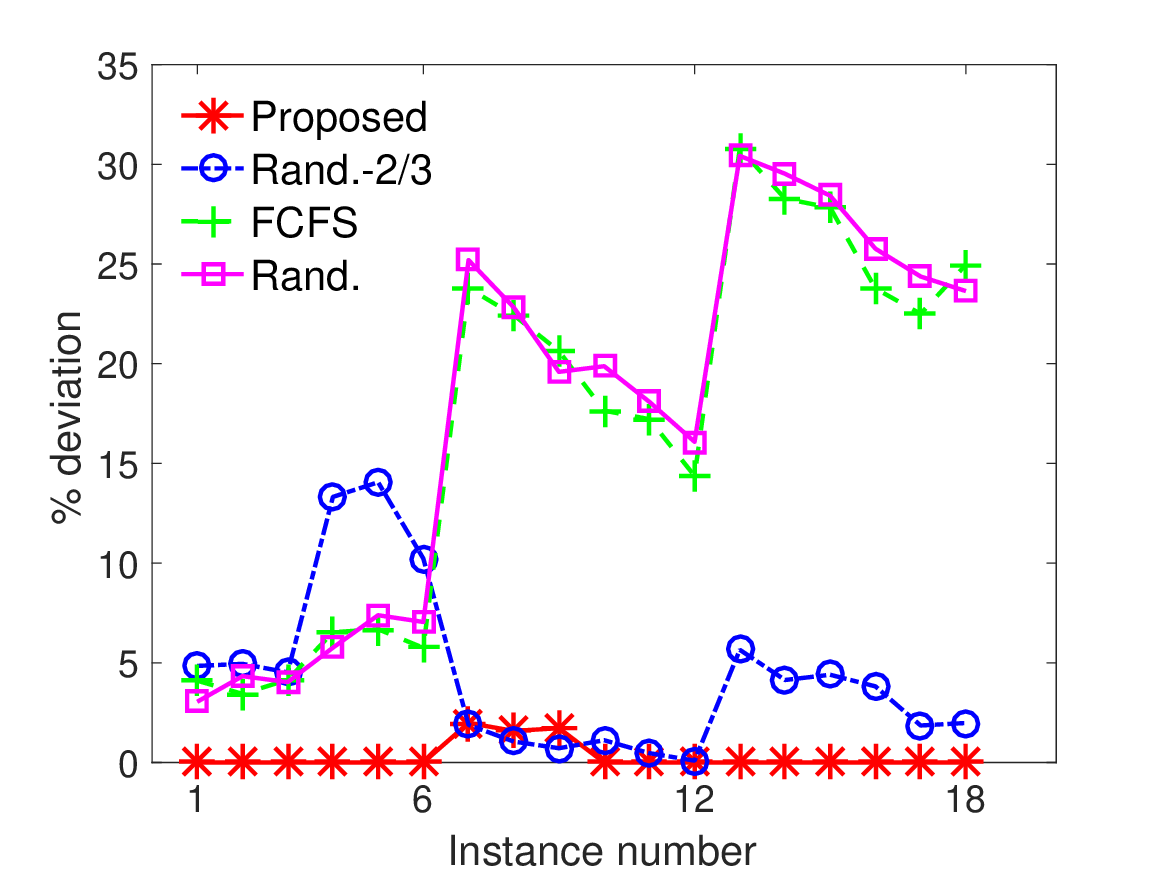}}
  \caption{(Color online) Percentage deviations of the costs of the policies from the lowest realized average cost ($L(i)$ for $i\in\{1,2,\ldots,18\}$) in the first 18 instances when $h_a=2$ and $h_b=1$.}
  \label{Figure_results}
\end{figure}

Figure \ref{Figure_results} shows $100\times(J_p(i)-L(i))/L(i)$ for each policy $p$, $p\in\mP$ at each $i$, $i\in\{1,2,\ldots,18\}$ when $h_a=2$ and $h_b=1$. Since $\mu_A=\mu_B$ and $h_a=2 h_b$, it makes more sense to give more priority to type $a$ jobs. This is why in addition to the policies in $\mP$, we also consider a \textit{randomized-$2/3$ policy}, in which whenever server 4 becomes available and if both of the buffers 4 and 5 are non-empty, server 4 chooses a type $a$ job with $2/3$ probability. If only one of the buffers 4 and 5 is non-empty, then server 4 chooses the job from the non-empty one. 

In Figure \ref{Figure_r1}, we see that the proposed and the SDP policies perform well at each instance. SDP policy performs well even at the instances $\{13,14,\ldots,18\}$, where $\mu_3>\mu_A$. However, the static priority policy does not perform well in the first 6 instances, where server 3 is in heavy traffic. This result is expected because when server 4 gives static priority to buffer 4 under the condition that server 3 is in heavy traffic, jobs in buffer 8 will wait for the jobs in buffer 7. In this case, it is more efficient to give priority to the type $b$ jobs in server 4 and this is exactly what the SDP policy does. The static priority policy performs the best at the instances $\{13,14,\ldots,18\}$, where $\mu_3>\mu_A$. Moreover, its performance is close to the one of the SDP policy at the instances $\{7,8,\ldots,12\}$, where $\lambda_a<\mu_3<\mu_A$. This result is not surprising because static priority policy is asymptotically optimal when server 3 is in light traffic (cf. Remark \ref{sp_optimality_range}). On the one hand, at the instances $\{7,8,9\}$, the static priority policy performs better than the SDP policy because server 5 is in heavy traffic at these instances, hence server 4 should not give priority to the type $b$ jobs at all. On the other hand, at the instances $\{10,11,12\}$, SDP policy performs better than the static priority policy, because server 5 is in light traffic at these instances; hence giving priority to type $b$ jobs in server 4 decreases the waiting time of the jobs in buffer 10.

\begin{figure}[htbp]
\begin{center}
\includegraphics[width=0.5\textwidth]{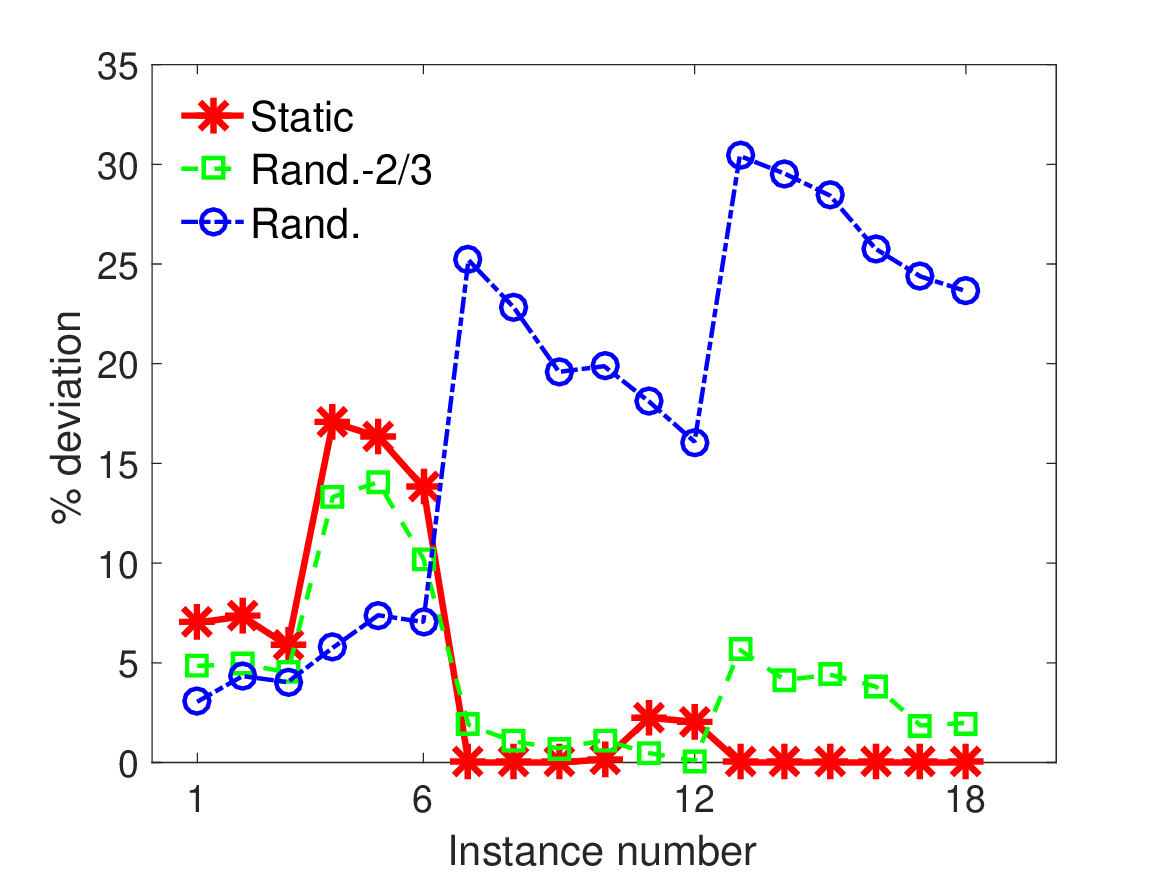}
\caption{(Color online) Percentage deviations of the costs of the static, randomized, and the randomized-$2/3$ policies from the lowest realized average cost ($L(i)$ for $i\in\{1,2,\ldots,18\}$) in the first 18 instances when $h_a=2$ and $h_b=1$.}\label{Figure_results_2}
\end{center}
\end{figure}

Figure \ref{Figure_r2} shows that the FCFS and the randomized policies perform worse than the proposed policy at each instance. The randomized-$2/3$ policy performs well only at the instances $\{7,8,\ldots,12\}$. Figure \ref{Figure_results_2}, where we compare the static priority, randomized, and randomized-$2/3$ policies, explains this result in the following way. On the one hand, at the first 6 instances, we see that the randomized and randomized-$2/3$ policies perform the best and the second best among these three policies, respectively. This implies that when server 3 is in heavy traffic, giving more priority to type $a$ jobs decreases the performance. On the other hand, at the instances $\{7,8,\ldots,18\}$, the static priority and randomized-$2/3$ policies perform the best and the second best among these three policies, respectively. This implies that when server 3 is in light traffic, giving more priority to type $a$ jobs increases the performance. In summary, the processing capacity of server 3 affects the performance of the control policy significantly. This is consistent with our theoretic results showing that an asymptotically optimal policy should not give static priority to type $a$ jobs when server 3 is in heavy traffic.

\begin{remark}\label{sp_optimality}
Figure \ref{Figure_r1} shows that the SDP policy performs very well at the instances $\{13,14,\ldots,18\}$ where $\mu_3>\mu_{A}$. Moreover, Table \ref{aggregate_results} shows that the SDP policy performs close to the proposed policy at all cases. Hence, we conjecture that the SDP policy is asymptotically optimal even when $\mu_3\geq\mu_{A}$. When $\mu_3\geq\mu_{A}$, the proofs associated with the up intervals (cf. Section \ref{ud_proof_1}) do not work. Note that we consider the difference of two renewal processes with rates $\mu_{A}$ and $\mu_3$, respectively in the proofs associated with the up intervals (cf. \eqref{ud_6}). We show that the renewal process with rate $\mu_{A}$ stays sufficiently close to the one with rate $\mu_3$ during the up intervals when $\mu_3<\mu_{A}$. However, this argument does not hold when $\mu_3\geq\mu_{A}$.
\end{remark}

\section{Extensions}\label{extensions}

It is possible to extend our results to more complex networks. We first consider a fork-join network in which there are task dependent holding costs in Section \ref{task_dep}. Next, we consider a fork-join network in which jobs fork into an arbitrary number of jobs in Section \ref{ex_fork}. Lastly, we consider fork-join networks with more than two job types in Section \ref{ex_type}. In each of these extensions, we first construct and solve an approximating DCP, then interpret a control policy from the solution. Since we omit the asymptotic optimality proofs corresponding to these extensions, the suggested control policies are all heuristics.

\subsection{Task Dependent Holding Costs}\label{task_dep}

So far, we assume that the holding cost rate of a type $a$ ($b$) job in buffers 1, 3, 4, 7, and 8 (2, 5, 6, 9, and 10) are the same and denoted by $h_a$ ($h_b$). In this section, we extend this assumption by considering task dependent holding cost rates. We denote the holding cost rate of a job in buffer $k$ as $h_k$ for all $k\in\mK$ and we assume $h_k\geq 0$ for all $k\in\mK$. In this case, our objective in the DCP is to minimize 
\begin{equation}\label{obj_t}
\pr\left( \sum_{k=1}^{10} h_k \ot{Q}_k(t) >x \right),\quad \forall t\geq 0,\; x>0.
\end{equation}
Since the control policy in server 4 has no effect on the costs in buffers 1, 2, 3, and 6, the objective \eqref{obj_t} can be simplified to minimizing
\begin{equation}\label{obj_t2}
\pr\left(h_4 \ot{Q}_4(t) + h_5 \ot{Q}_5(t) + h_7 \ot{Q}_7(t) + h_8 \ot{Q}_8(t) + h_9 \ot{Q}_9(t) + h_{10} \ot{Q}_{10}(t) >x \right),\quad \forall t\geq 0,\; x>0.
\end{equation}
By first simplifying objective \eqref{obj_t2} by \eqref{queue_inst}, which is the instantaneous service process assumption in the downstream servers, we can modify the DCP \eqref{DCP_4} in the following way: For each $x>0$ and $t\geq 0$, 
\begin{align}
\min\quad &\pr\bigg(h_4 \ot{Q}_4(t) + h_5 \ot{Q}_5(t) + h_7 \left(\ot{Q}_4(t)-\ot{Q}_3(t)\right)^+ + h_8 \left(\ot{Q}_3(t)-\ot{Q}_4(t)\right)^+   \nonumber\\
&\hspace{5cm}  + h_9 \left(\ot{Q}_6(t)-\ot{Q}_5(t)\right)^+ + h_{10} \left(\ot{Q}_5(t)-\ot{Q}_6(t)\right)^+  >x \bigg),\label{DCP_5}\\
&\hspace{2cm}\text{s.t. }  \ot{Q}_4(t) + \frac{\mu_{A}}{\mu_{B}} \ot{Q}_5(t) = \ot{W}_4(t), \nonumber\\
& \hspace{2cm}\qquad \ot{Q}_k(t)\geq 0,\quad\text{for all }k\in\{4,5\}. \nonumber
\end{align}
When we consider the DCP \eqref{DCP_5} path-wise, we have the following optimization problem for all $t\in\R_+$ and $\omega$ in $\Omega$ except a null set:
\begin{align}
\min\quad &h_4 \ot{Q}_4(\omega_t) + h_5 \ot{Q}_5(\omega_t) + h_7 \left(\ot{Q}_4(\omega_t)-\ot{Q}_3(\omega_t)\right)^+ + h_8 \left(\ot{Q}_3(\omega_t)-\ot{Q}_4(\omega_t)\right)^+   \nonumber\\
&\hspace{5cm}  + h_9 \left(\ot{Q}_6(\omega_t)-\ot{Q}_5(\omega_t)\right)^+ + h_{10} \left(\ot{Q}_5(\omega_t)-\ot{Q}_6(\omega_t)\right)^+ ,\label{DCP_6}\\
&\hspace{2cm}\text{s.t. }  \ot{Q}_4(\omega_t) + \frac{\mu_{A}}{\mu_{B}} \ot{Q}_5(\omega_t) = \ot{W}_4(\omega_t), \nonumber\\
& \hspace{2cm}\qquad \ot{Q}_k(\omega_t)\geq 0,\quad\text{for all }k\in\{4,5\}. \nonumber
\end{align}
We can solve the optimization problem \eqref{DCP_6} by the following lemma whose proof is presented in Appendix \ref{opt_solution_2_proof}.

\begin{lemma}\label{opt_solution_2}
Consider the optimization problem
\begin{subequations}\label{opt_problem}
\begin{align}
\min\quad &h_4 q_4 + h_5 q_5 + h_7\left(q_4-q_3\right)^+ + h_8\left(q_3-q_4\right)^+ + h_9\left(q_6-q_5\right)^+ + h_{10}\left(q_5-q_6\right)^+ ,\label{opt_problem_1}\\
&\hspace{2cm}\text{s.t.}\quad  q_4 + \frac{\mu_{A}}{\mu_{B}} q_5 = w_4, \label{opt_problem_2}\\
&\hspace{2cm}\qquad\;\; q_4\geq 0,\;q_5\geq 0, \label{opt_problem_3}
\end{align}
\end{subequations}
where $q_4$ and $q_5$ are the decision variables and all of the parameters are nonnegative. Then, there exists an optimal solution among the four solutions given in Table \ref{sln_table} below.
\begin{table}[htbp]\scriptsize
\centering \caption{An optimal solution set for the optimization problem \eqref{opt_problem}.}
\begin{tabular}{c c c}
\hline
\textbf{Solution \#}  & $\bm{q_4}$ & $\bm{q_5}$ \\
\hline
1 & 0 & $(\mu_{B}/\mu_{A}) w_4$   \\

2 & $q_3\wedge w_4$ &  $(\mu_{B}/\mu_{A}) (w_4-q_3)^+ $ \\

3 & $\left(w_4-(\mu_{A}/\mu_{B}) q_6\right)^+$ & $q_6\wedge (\mu_{B}/\mu_{A}) w_4 $ \\

4 & $w_4$ & 0   \\
\hline
\end{tabular}\label{sln_table}
\end{table}
\end{lemma}

\begin{remark}\label{sln_remark}
The optimal solution of Lemma \ref{opt_solution_2} is strictly dependent on the parameters $(q_3,q_6,w_4)$, i.e., the optimal solution can change as $q_3$, or $q_6$, or $w_4$ changes. For example, consider an example in which $\mu_A=\mu_B=1$, $h_4=h_5=1$, $h_8=h_{10}=2$, $h_7=h_9=3$, and $(q_3,q_6,w_4)=(0.01, 0,1)$. Then, the objective function values corresponding to the four solutions in Table \ref{sln_table} are $(3.02,2.98,3.97,3.97)$, respectively. This implies that an optimal solution is the second solution in Table \ref{sln_table}, which is $(q_4,q_5)=(0.01,0.99)$. Next, we consider the same example but this time $(q_3,q_6,w_4)=(1,0.01,1)$, i.e., $q_3$ and $q_6$ change. In this case, the objective function values corresponding to the four solutions in Table \ref{sln_table} are $(4.98,1.03,1.02,1.03)$, respectively. This implies that an optimal solution is the third solution in Table \ref{sln_table}, which is $(q_4,q_5)=(0.99,0.01)$. \end{remark}

By Lemma \ref{opt_solution_2}, for each $t\in\R_+$ and $\omega$ in $\Omega$ except a null set, an optimal solution of DCP \eqref{DCP_6} is among the four solutions presented below, which then motivates a control policy.
\begin{enumerate}
\item If 
\begin{equation}\label{task_sln_1}
\left(\ot{Q}_4(\omega_t),\ot{Q}_5(\omega_t)\right)=\left(0,\; \frac{\mu_{B}}{\mu_{A}} \ot{W}_4(\omega_t) \right),
\end{equation}
then server 4 should give static priority to type $a$ jobs. 

\item If
\begin{equation}\label{task_sln_2}
\left(\ot{Q}_4(\omega_t),\ot{Q}_5(\omega_t)\right)= \left(\ot{Q}_3(\omega_t)\wedge\ot{W}_4(\omega_t),\;\frac{\mu_{B}}{\mu_{A}}\left(\ot{W}_4(\omega_t)-\ot{Q}_3(\omega_t)\right)^+\right),
\end{equation}
then server 4 should use the proposed policy. Note that the solution in \eqref{task_sln_2} is the same as the one in Proposition \ref{lb_DCP}.

\item If
\begin{equation}\label{task_sln_3}
\left(\ot{Q}_4(\omega_t),\ot{Q}_5(\omega_t)\right)= \left(\left(\ot{W}_4(\omega_t)-\frac{\mu_{A}}{\mu_{B}}\ot{Q}_6(\omega_t)\right)^+,\; \ot{Q}_6(\omega_t)\wedge \frac{\mu_{B}}{\mu_{A}} \ot{W}_4(\omega_t) \right),
\end{equation}
then server 4 should use the proposed policy but this time it should prioritize type $b$ jobs over type $a$ jobs. In other words, if $\mu_B\leq \mu_5$, server 4 should give static priority to type $b$ jobs; otherwise, it should use the SDP policy to pace the departure process of type $b$ jobs from buffer 5 with the ones from buffer 6.

\item If
\begin{equation}\label{task_sln_4}
\left(\ot{Q}_4(\omega_t),\ot{Q}_5(\omega_t)\right)= \left(\ot{W}_4(\omega_t),\;0\right),
\end{equation}
then server 4 should give static priority to type $b$ jobs. 
\end{enumerate}

Since the optimal solution can change as the process $\big(\ot{Q}_3(\omega_t),\ot{Q}_6(\omega_t),\ot{W}_4(\omega_t)\big)$ changes with time (cf. Remark \ref{sln_remark}), we interpret the following dynamic control policy: Whenever server 4 makes a service completion, the system controller computes the objective function values of DCP \eqref{DCP_6} corresponding to the four solutions in \eqref{task_sln_1} - \eqref{task_sln_4} and implements the control policy corresponding to the solution with the minimum objective function value until the next service completion epoch in server 4.

\subsection{Networks with Arbitrary Number of Forks}\label{ex_fork}

\begin{figure}[htb]
\begin{center}
\includegraphics[width=0.8\textwidth]{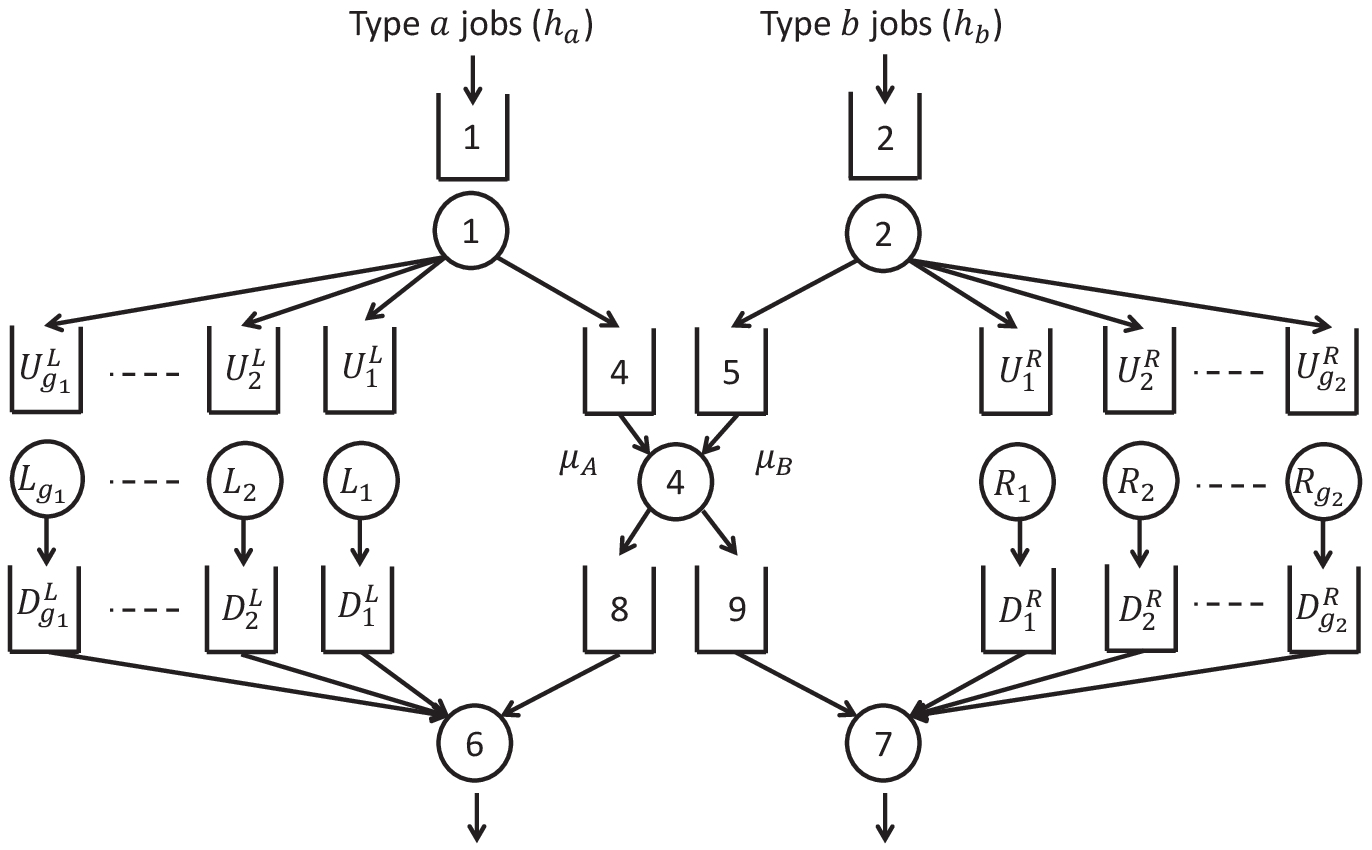}
\caption{A fork-join processing network with two job types and arbitrary number of forks.}\label{fj_network_2}
\end{center}
\end{figure}

In this section, we consider the network presented in Figure \ref{fj_network_2}, in which type $a$ and $b$ jobs fork into $g_1+1$ and $g_2+1$ number of jobs, respectively, where $g_1,g_2\in\N_+$; and server 4 is the only server which processes both job types.

Parallel with objective \eqref{obj_5}, we consider the objective of minimizing
\begin{equation}\label{obj_8}
\pr\left(h_a (g_1+1)\left(Q_{U_1^L}^r(t)+Q_{D_1^L}^r(t)\right)+h_b (g_2+1)\left(Q_{U_1^R}^r(t)+Q_{D_1^R}^r(t)\right)>x \right),
\end{equation}
for all $t\in\R_+$ and $x>0$. We assume that servers 6 and 7 are in light traffic, hence we can assume that the service processes in these servers are instantaneous as in Section \ref{DCP}. Then, similar to \eqref{queue_inst}, we have
\begin{equation}\label{ins_queue}
Q_{D_1^L}^r =\left(Q_4^r\vee \left(\max_{i\in\{2,\ldots,g_1\}} Q_{U_i^L}^r\right)-Q_{U_1^L}^r\right)^+, \qquad Q_{D_1^R}^r =\left(Q_5^r\vee \left(\max_{i\in\{2,\ldots,g_2\}} Q_{U_i^R}^r\right)-Q_{U_1^R}^r\right)^+.
\end{equation}
From \eqref{ins_queue}, defining $\ot{h}_a := h_a (g_1+1)$ and $\ot{h}_b := h_b (g_2+1)$, and the fact that $x+(y-x)^+=x\vee y$ for all $x,y\in\R$, the objective \eqref{obj_8} is equivalent to minimizing
\begin{equation}\label{obj_9}
\pr\left(\ot{h}_a \left(Q_4^r(t)\vee \left(\max_{i\in\{1,\ldots,g_1\}} Q_{U_i^L}^r(t)\right)\right) +\ot{h}_b\left(Q_5^r(t)\vee \left(\max_{i\in\{1,\ldots,g_2\}} Q_{U_i^R}^r(t)\right)\right)>x \right),
\end{equation}
for all $t\in\R_+$ and $x>0$.

Parallel with Proposition \ref{general_conv}, we can prove the following weak convergence result under any work-conserving control policy:
\begin{align}
&\left(\oh{Q}^r_1,\oh{Q}^r_2,\;\oh{Q}^r_{U_i^L},i\in\{1,2,\ldots,g_1\},\;\oh{Q}^r_{U_j^R},j\in\{1,2,\ldots,g_2\},\;\oh{W}^r_4\right) \nonumber\\
&\hspace{2cm} \Longrightarrow  \left(\ot{Q}_1,\ot{Q}_2,\;\ot{Q}_{U_i^L},i\in\{1,2,\ldots,g_1\},\;\ot{Q}_{U_j^R},j\in\{1,2,\ldots,g_2\},\;\ot{W}_4\right),\label{general_conv_3}
\end{align} 
where the limiting process is the zero process for the buffers whose corresponding dedicated server is in light traffic and an SRBM for the buffers whose corresponding dedicated server is in heavy traffic and the workload process in server 4. Then, by \eqref{obj_9}, \eqref{general_conv_3}, and using the technique that we use to derive the DCP \eqref{DCP_4}, we construct the following DCP for this network. For each $x>0$ and $t\geq 0$,
\begin{align}
\min\quad &\pr\left(\ot{h}_a \left(\ot{Q}_4(t)\vee \left(\max_{i\in\{1,\ldots,g_1\}} \ot{Q}_{U_i^L}(t)\right)\right) +\ot{h}_b\left(\ot{Q}_5(t)\vee \left(\max_{i\in\{1,\ldots,g_2\}} \ot{Q}_{U_i^R}(t)\right)\right)>x \right),\nonumber\\
& \hspace{2cm} \text{s.t.}\quad\;  \ot{Q}_4(t) +\frac{\mu_{A}}{\mu_{B}} \ot{Q}_5(t) =  \ot{W}_4(t), \label{DCP_13}\\
& \hspace{3cm} \ot{Q}_k(t)\geq 0,\quad\text{for all }k\in\{4,5\}.\nonumber
\end{align}
When we consider the DCP \eqref{DCP_13} path-wise, we have the following optimization problem for all $t\in\R_+$ and $\omega$ in $\Omega$ except a null set:
\begin{subequations}\label{DCP_14}
\begin{align}
\min\quad &\ot{h}_a \left(\ot{Q}_4(\omega_t)\vee \left(\max_{i\in\{1,\ldots,g_1\}} \ot{Q}_{U_i^L}(\omega_t)\right)\right) +\ot{h}_b\left(\ot{Q}_5(\omega_t)\vee \left(\max_{i\in\{1,\ldots,g_2\}} \ot{Q}_{U_i^R}(\omega_t)\right)\right),\label{DCP_14_1}\\
\text{s.t.}\quad & \ot{Q}_4(\omega_t) + \frac{\mu_{A}}{\mu_{B}} \ot{Q}_5(\omega_t) = \ot{W}_4(\omega_t), \label{DCP_14_2}\\
& \ot{Q}_k(\omega_t)\geq 0,\quad\text{for all }k\in\{4,5\}. \label{DCP_14_3}
\end{align}
\end{subequations}
where  $\ot{Q}_k(\omega_t)$, $k\in\{4,5\}$ are the decision variables. Note that the objective function \eqref{DCP_14_1} is equal to
\begin{align}
& \ot{h}_a \left(\max_{i\in\{1,\ldots,g_1\}} \ot{Q}_{U_i^L}(\omega_t)\right)+  \ot{h}_a\left(\ot{Q}_4(\omega_t)- \left(\max_{i\in\{1,\ldots,g_1\}} \ot{Q}_{U_i^L}(\omega_t)\right)\right)^+ \nonumber\\
&\hspace{2cm} + \ot{h}_b  \left(\max_{i\in\{1,\ldots,g_2\}} \ot{Q}_{U_i^R}(\omega_t)\right)+\ot{h}_b\left(\ot{Q}_5(\omega_t)- \left(\max_{i\in\{1,\ldots,g_2\}} \ot{Q}_{U_i^R}(\omega_t)\right)\right)^+.\label{DCP_15}
\end{align}
Since all of the buffers $U_i^L$, $i\in\{1,\ldots,g_1\}$ and $U_i^R$, $i\in\{1,\ldots,g_2\}$ are independent of the control, minimizing the objective in \eqref{DCP_15} is equivalent to minimizing the objective \eqref{DCP_16_1} below, thus the optimization problem \eqref{DCP_14} is equivalent to the following one:
\begin{subequations}\label{DCP_16}
\begin{align}
\min\quad &\ot{h}_a\left(\ot{Q}_4(\omega_t)- \left(\max_{i\in\{1,\ldots,g_1\}} \ot{Q}_{U_i^L}(\omega_t)\right)\right)^+ + \ot{h}_b\left(\ot{Q}_5(\omega_t)- \left(\max_{i\in\{1,\ldots,g_2\}} \ot{Q}_{U_i^R}(\omega_t)\right)\right)^+,\label{DCP_16_1}\\
\text{s.t.}\quad & \ot{Q}_4(\omega_t) + \frac{\mu_{A}}{\mu_{B}} \ot{Q}_5(\omega_t) = \ot{W}_4(\omega_t), \label{DCP_16_2}\\
& \ot{Q}_k(\omega_t)\geq 0,\quad\text{for all }k\in\{4,5\}. \label{DCP_16_3}
\end{align}
\end{subequations}
Without loss of generality, let us assume that $\ot{h}_a \mu_{A}\geq \ot{h}_b \mu_{B}$. Then by Lemma \ref{opt_solution}, an optimal solution of the optimization problem \eqref{DCP_16} is
\begin{equation}\label{opt_sln_fork}
\left(\ot{Q}_4,\ot{Q}_5\right) = \left(\ot{W}_4\wedge \left(\max_{i\in\{1,\ldots,g_1\}} \ot{Q}_{U_i^L}\right), \; \frac{\mu_{B}}{\mu_{A}} \left(\ot{W}_4 - \left(\max_{i\in\{1,\ldots,g_1\}} \ot{Q}_{U_i^L}\right)\right)^+ \right).
\end{equation}
We can interpret the solution \eqref{opt_sln_fork} in the following way: If the processing capacity of each of the servers $L_i$, $i\in\{1,2,\ldots,g_1\}$ is greater than or equal to $\mu_A$, then server 4 should give static priority to type $a$ jobs all the time. Otherwise, server 4 should give priority to type $a$ jobs whenever the number of jobs in buffer 4 is strictly greater than the maximum of the number of jobs in buffers $U_i^L$, $i\in\{1,2,\ldots,g_1\}$. In other words, server 4 should pace the departure process of type $a$ jobs from buffer 4 with the minimum of the ones from the buffers $U_i^L$, $i\in\{1,2,\ldots,g_1\}$. Hence, we see a slightly different version of the SDP policy.

\subsection{Networks with More Than Two Job Types}\label{ex_type}

In this section we consider fork-join networks with more than two job types. We first consider a network with an arbitrary number of job types and a single shared server in Section \ref{ex_type_1}, and then consider a network with three job types and two shared servers in Section \ref{ex_type_2}.

\subsubsection{Networks with Arbitrary Job Types and a Single Shared Sever}\label{ex_type_1}

\begin{figure}[htb]
\begin{center}
\includegraphics[width=0.7\textwidth]{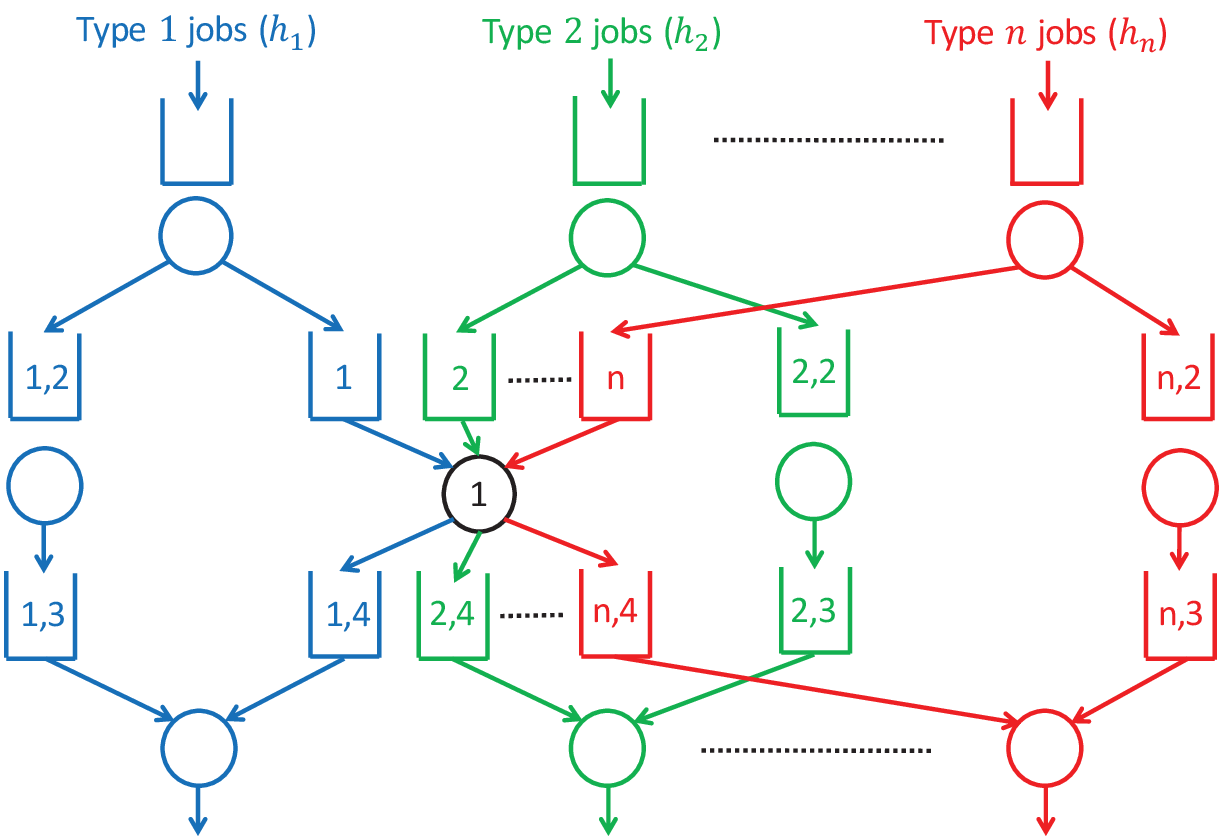}
\caption{(Color online) A fork-join processing network with arbitrary job types and a single shared server.}\label{fj_network_5}
\end{center}
\end{figure}

In this section, we consider the network presented in Figure \ref{fj_network_5}, where there are $n$ job types such that $n$ is arbitrary and $n\geq 2$. Server 1 is a shared server which processes all job types, whereas all other servers process only a single job type. Upon arriving to the system, each job is first processed in a server, then it is forked into two jobs which are processed in the shared server and in a dedicated server, respectively, lastly the two forked jobs are joined in the corresponding downstream server and leave the system. 

A type $j$ job can wait in buffers $j$, $(j,2)$, $(j,3)$, and $(j,4)$, where buffer $j$ feeds the shared server, buffer $(j,2)$ feeds the dedicated server, buffer $(j,3)$ is fed by the dedicated server, and buffer $(j,4)$ is fed by the shared server (cf. Figure \ref{fj_network_5}). Server 1, the shared server, processes type $j$ jobs with rate $\mu_j$ for all $j\in\{1,2,\ldots,n\}$. Let $h_j$ denote the cost per unit time to hold a type $j$ job in the system. Without loss of generality, we assume that $h_1\mu_1\geq h_2\mu_2\geq \ldots \geq h_n\mu_n$. We assume that server 1 is in heavy traffic and all of the downstream servers are in light traffic. Parallel with \eqref{workload_process}, let
\begin{equation*}
 \oh{W}^r := \sum_{j=1}^n \frac{\oh{Q}_j^r}{\mu_j^r}
\end{equation*}
denote the diffusion scaled workload process that server 1 sees. Then, parallel with Proposition \ref{general_conv}, we can prove the following weak convergence result under any work-conserving control policy:
\begin{equation}\label{general_conv_4}
\left(\oh{Q}^r_{j,2},j\in\{1,2,\ldots,n\},\;\oh{W}^r\right) \Longrightarrow  \left(\ot{Q}_{j,2},j\in\{1,2,\ldots,n\},\;\ot{W}\right),
\end{equation}
where the limiting process is the zero process for the buffers whose corresponding dedicated server is in light traffic and an SRBM for the buffers whose corresponding dedicated server is in heavy traffic and the workload process in server 1. Since we assume that the downstream servers are in light traffic, we can assume that the service processes in these servers are instantaneous. Then, similar to \eqref{queue_inst}, we have
\begin{equation*}
Q_{j,3}^r=\left(Q_j^r-Q_{j,2}^r\right)^+,\quad Q_{j,4}^r=\left(Q_{j,2}^r-Q_j^r\right)^+,\qquad\forall j\in\{1,2,\ldots,n\}.
\end{equation*}
By using the technique that we use to derive DCP \eqref{DCP_4}, we construct the following DCP for this network. For each $x>0$ and $t\geq 0$,
\begin{align}
& \min\quad \pr\left(\sum_{j=1}^n h_j \left(\ot{Q}_{j,2}(t)+\left(\ot{Q}_j(t) - \ot{Q}_{j,2}(t)\right)^+\right) > x \right),\nonumber\\
& \hspace{2cm} \text{s.t.}\quad\; \sum_{j=1}^n \ot{Q}_j(t)/\mu_j =  \ot{W}(t), \label{DCP_7}\\
& \hspace{3cm} \ot{Q}_j(t)\geq 0,\quad\text{for all }j\in\{1,2,\ldots,n\}.\nonumber
\end{align}

When we consider the DCP \eqref{DCP_7} path-wise, we have the following optimization problem:
\begin{align}
& \min\quad \sum_{j=1}^n h_j \left(\ot{Q}_j (\omega_t) - \ot{Q}_{j,2} (\omega_t)\right)^+ ,\nonumber\\
&\hspace{0.5cm} \text{s.t.}\quad\; \sum_{j=1}^n \ot{Q}_j(\omega_t)/\mu_j =  \ot{W}(\omega_t), \label{DCP_8}\\
&\hspace{1.5cm}  \ot{Q}_j(\omega_t)\geq 0,\quad\text{for all }j\in\{1,2,\ldots,n\}.\nonumber
\end{align}
where  $\ot{Q}_j(\omega_t)$, $j\in\{1,2,\ldots,n\}$ are the decision variables. Note that, since $\ot{Q}_{j,2}$ is an exogenous process for all $j\in\{1,2,\ldots,n\}$ (cf. \eqref{general_conv_4}), we neglect the term $ \sum_{j=1}^n h_j \ot{Q}_{j,2}(\omega_t)$ in the objective function of the optimization problem \eqref{DCP_8}. We can solve the optimization problem \eqref{DCP_8} by the following lemma whose proof is presented in Appendix \ref{opt_solution_3_proof}.

\begin{lemma}\label{opt_solution_3}
Consider the optimization problem
\begin{subequations}\label{opt_problem_1_0}
\begin{align}
& \min\quad \sum_{j=1}^n h_j \left(q_j -q_{j,2} \right)^+ ,\label{opt_problem_1_1}\\
&\hspace{0.5cm} \text{s.t.}\quad\; \sum_{j=1}^n q_j /\mu_j =  w, \label{opt_problem_1_2}\\
&\hspace{1.5cm}  q_j\geq 0,\quad\text{for all }j\in\{1,2,\ldots,n\}.\label{opt_problem_1_3}
\end{align}
\end{subequations}
where $q_j$,  $j\in\{1,2,\ldots,n\}$ are the decision variables, all of the parameters are nonnegative, and $h_1\mu_1\geq h_2\mu_2\geq \ldots \geq h_n\mu_n$. Then, an optimal solution of the optimization problem \eqref{opt_problem_1_0} is
\begin{subequations}\label{opt_problem_1_4}
\begin{align}
& q_1 = q_{1,2} \wedge (\mu_1 w),\label{opt_problem_1_41}\\
& q_j = q_{j,2} \wedge \left[\mu_j \left(\cdots\left(\left( w - \frac{q_{1,2}}{\mu_1}\right)^+ -  \frac{q_{2,2}}{\mu_2}\right)^+ - \cdots -  \frac{q_{j-1,2}}{\mu_{j-1}}\right)^+\right],\quad\forall j\in\{2,3,\ldots,n-1\},\label{opt_problem_1_42}\\
& q_n = \mu_n \left(\cdots\left(\left( w - \frac{q_{1,2}}{\mu_1}\right)^+ -  \frac{q_{2,2}}{\mu_2}\right)^+ - \cdots -  \frac{q_{n-1,2}}{\mu_{n-1}}\right)^+.\label{opt_problem_1_43}
\end{align}
\end{subequations}
\end{lemma}

\begin{remark}
When $n=2$, the solution of Lemma \ref{opt_solution_3} is the same as the solution of Lemma \ref{opt_solution} by the fact that the notation $(\mu_1,\mu_2,q_{1,2},q_{2,2},w)$ in Lemma \ref{opt_solution_3} corresponds to the notation $(\mu_A,\mu_B,q_3,q_6,w_4/\mu_A)$ in Lemma \ref{opt_solution}.
\end{remark}

\begin{remark}\label{comment_sln}
In the optimal solution \eqref{opt_problem_1_4}, we see that $h_j \left(q_j -q_{j,2} \right)^+=0$ for all $j\in\{1,2,\ldots,n-1\}$ but $h_n \left(q_n -q_{n,2} \right)^+$ can be strictly positive depending on the parameters. This implies that the objective function \eqref{opt_problem_1_1} can become positive only because of type $n$ jobs. In other words, cost can occur only because of the cheapest job type.
\end{remark}

Lemma \ref{opt_solution_3} implies that an optimal solution of the optimization problem \eqref{DCP_8} is
\begin{subequations}\label{opt_problem_1_5}
\begin{align}
& \ot{Q}_1 = \ot{Q}_{1,2} \wedge (\mu_1 \ot{W}),\label{opt_problem_1_51}\\
& \ot{Q}_j = \ot{Q}_{j,2} \wedge \left[\mu_j \left(\cdots\left(\left( \ot{W} - \frac{\ot{Q}_{1,2}}{\mu_1}\right)^+ -  \frac{\ot{Q}_{2,2}}{\mu_2}\right)^+ - \cdots -  \frac{\ot{Q}_{j-1,2}}{\mu_{j-1}}\right)^+\right],\forall j\in\{2,3,\ldots,n-1\},\label{opt_problem_1_52}\\
& \ot{Q}_n = \mu_n \left(\cdots\left(\left( \ot{W} - \frac{\ot{Q}_{1,2}}{\mu_1}\right)^+ -  \frac{\ot{Q}_{2,2}}{\mu_2}\right)^+ - \cdots -  \frac{\ot{Q}_{n-1,2}}{\mu_{n-1}}\right)^+.\label{opt_problem_1_53}
\end{align}
\end{subequations}
Then, we can interpret the following control policy from the solution \eqref{opt_problem_1_5}: 

\begin{itemize}

\item The priority ranking of the job types is $1\geq 2\geq \ldots \geq n$.

\item \textit{\textbf{Type 1 jobs}}: If $\mu_1$ is less than or equal to the processing rate of the dedicated server for type 1 jobs, server 1 should give static priority to type 1 jobs all the time. Otherwise, server 1 should give priority to buffer 1 only when the number of jobs in buffer 1,2 becomes less than the one in buffer 1. During the remaining time, server 1 should process the job types $j\in\{2,\ldots,n\}$.

\item \textit{\textbf{Type $\bm{j}$ jobs, $\bm{j\in\{2,\ldots,n-1\}$}}}: Server 1 gives priority to type $j$ jobs only when the higher priority job types ($\{1,\ldots,j-1\}$) do not require any processing. Consider the time intervals in which server 1 gives priority to type $j$ jobs. If $\mu_j$ is less than or equal to the processing rate of the dedicated server for type $j$ jobs, server 1 should give static priority to type $j$ jobs all the time. Otherwise, server 1 should give priority to buffer $j$ only when the number of jobs in buffer $j,2$ becomes less than the one in buffer $j$. During the remaining time, server 1 should process the job types $i\in\{j+1,\ldots,n\}$.

\item \textit{\textbf{Type $\bm{n}$ jobs}} are processed during the remaining time.

\end{itemize}

Therefore, we see a \textit{chained implementation} of the proposed policy. The throughput rate of type $j$ jobs is maximized by keeping buffer $j$ less than or equal to buffer $(j,2)$ with minimum effort for all $j\in\{1,\ldots,n-1\}$ (cf. \eqref{opt_problem_1_51}, \eqref{opt_problem_1_52}, and Remark \ref{comment_sln}). By this way, server 4 gives as much as priority to type $n$ jobs.

\begin{remark}\label{most_general_network}
It is possible to extend the network in Figure \ref{fj_network_5} to the case where each job type forks into arbitrary number of jobs (cf. Figure \ref{fj_network_0}). Suppose that type $j$ jobs fork into $g_j+1$ jobs and $\ot{Q}_{j,2,i}$, $i\in\{1,2,\ldots,g_j\}$ denote the buffers in front of the dedicated servers corresponding to the type $j$ jobs. Then by the DCP \eqref{DCP_13} and the optimization problem \eqref{DCP_16} in Section \ref{ex_fork}, we need to replace $\ot{Q}_{j,2}$ with $\max_{i\in\{1,\ldots,g_j\}}\ot{Q}_{j,2,i}$ in both the DCP \eqref{DCP_7} and the optimization problem \eqref{DCP_8}. We can still solve the modified version of the optimization problem \eqref{DCP_8} with Lemma \ref{opt_solution_3}. The control policy that we can interpret is a modification of the control policy that we interpret from the solution \eqref{opt_problem_1_5} with the one from the solution \eqref{opt_sln_fork}.
\end{remark}

\subsubsection{Networks with Three Job Types and Two Shared Servers}\label{ex_type_2}

\begin{figure}[htb]
\begin{center}
\includegraphics[width=0.7\textwidth]{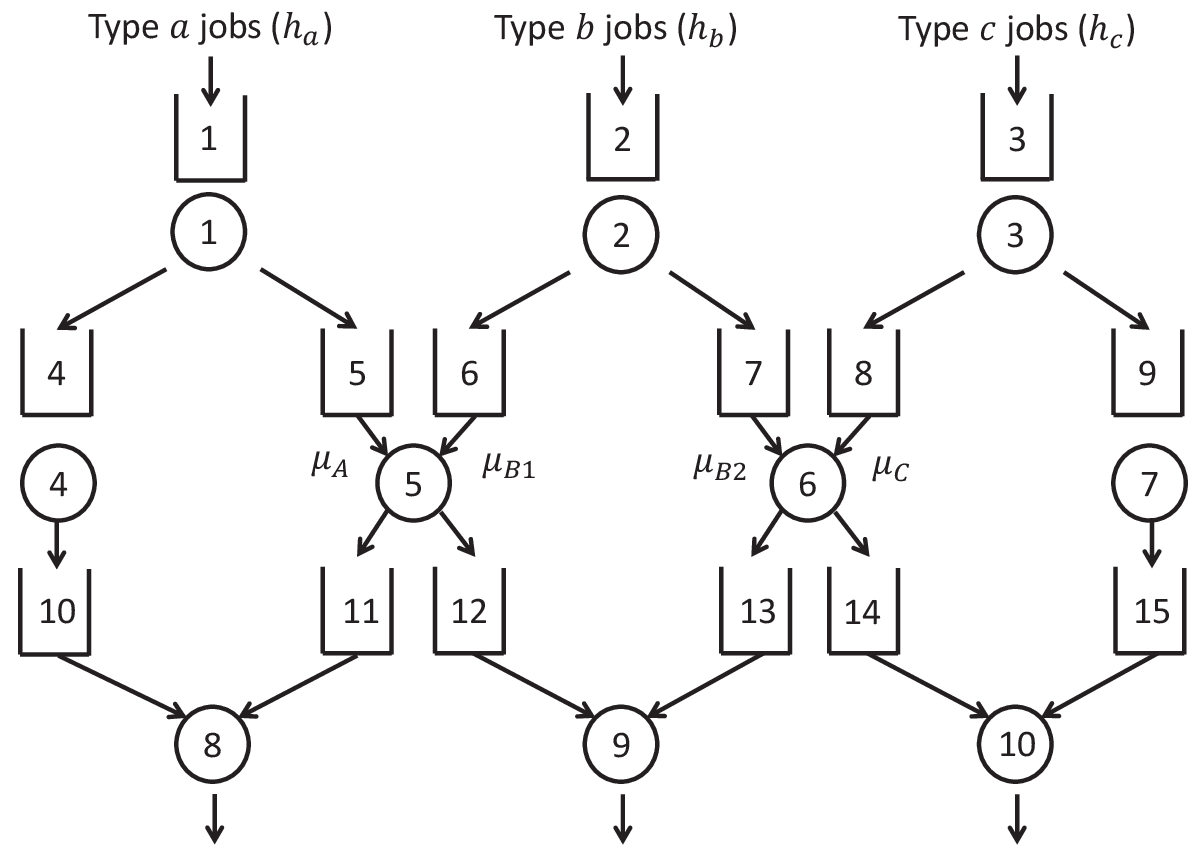}
\caption{A fork-join processing network with three job types and two shared servers.}\label{fj_network_3}
\end{center}
\end{figure}

In this section, we consider the network presented in Figure \ref{fj_network_3}, where there are three job types, namely type $a$, $b$, and $c$ jobs. There are two shared servers, which are servers 5 and 6. We formulate and solve the approximating DCP in order to derive heuristic control policies. In contrast with the networks in Figures \ref{fj_network_2} and \ref{fj_network_5}, the DCP solution is not a straightforward extension of the DCP solution presented in Section \ref{DCP}. 

Server 5 processes both type $a$ and $b$ jobs with rates $\mu_A$ and $\mu_{B1}$, respectively; and server 6 processes both type $b$ and $c$ jobs with rates $\mu_{B2}$ and $\mu_C$, respectively. Let $h_a$, $h_b$, and $h_c$ denote the holding cost rate per job per unit time for type $a$, $b$, and $c$ jobs, respectively. Suppose that both servers 5 and 6 are in heavy traffic, and servers 8, 9, and 10 are in light traffic. Let 
\begin{equation*}
 \oh{W}_5^r := \oh{Q}_5^r +\frac{\mu_{A}^r}{\mu_{B1}^r} \oh{Q}_6^r,\hspace{1cm} \oh{W}_6^r := \oh{Q}_7^r +\frac{\mu_{B2}^r}{\mu_{C}^r} \oh{Q}_8^r
\end{equation*}
denote the workload process (up to a constant scale factor) in servers 5 and 6, respectively. Then, parallel with Proposition \ref{general_conv}, we can prove the following weak convergence result under any work-conserving control policy:
\begin{equation}\label{general_conv_2}
\left(\oh{Q}^r_1,\oh{Q}^r_2,\oh{Q}^r_3,\oh{Q}^r_4,\oh{Q}^r_9,\oh{W}^r_5,\oh{W}^r_6\right) \Longrightarrow  \left(\ot{Q}_1,\ot{Q}_2,\ot{Q}_3,\ot{Q}_4,\ot{Q}_9,\ot{W}_5,\ot{W}_6 \right),
\end{equation}
where the limiting process is the zero process for the buffers whose corresponding dedicated server is in light traffic and an SRBM for the buffers whose corresponding dedicated server is in heavy traffic and the workload processes in servers 5 and 6. Since servers 8, 9, and 10 are in light traffic, we can assume that the service processes in these servers are instantaneous; and similar to \eqref{queue_inst}, we have
\begin{equation*}
Q_{10}^r=\left(Q_5^r-Q_4^r\right)^+, \; Q_{15}^r=\left(Q_8^r-Q_9^r\right)^+, \; Q_{6}^r + Q_{12}^r=Q_{6}^r  + \left(Q_7^r-Q_6^r\right)^+ = Q_6^r \vee Q_7^r.
\end{equation*}
Then, by using the technique that we use to derive DCP \eqref{DCP_4}, we construct the following DCP for this network. For any $x>0$ and $t\geq 0$,
\begin{align}
& \min\quad \pr\Bigg(h_a \left(\ot{Q}_4(t)+\left(\ot{Q}_5(t)-\ot{Q}_4(t)\right)^+\right)+h_b \left(\ot{Q}_6(t)\vee\ot{Q}_7(t)\right) \nonumber\\
&\hspace{7cm} + h_c \left(\ot{Q}_9(t)+\left(\ot{Q}_8(t)-\ot{Q}_9(t)\right)^+\right)>x \Bigg),\nonumber\\
& \hspace{2cm} \text{s.t.}\quad\;  \ot{Q}_5(t) +\frac{\mu_{A}}{\mu_{B1}} \ot{Q}_6(t) =  \ot{W}_5(t), \label{DCP_11}\\
& \hspace{3cm} \ot{Q}_7(t) +\frac{\mu_{B2}}{\mu_{C}} \ot{Q}_8(t) =  \ot{W}_6(t), \nonumber\\
& \hspace{3cm} \ot{Q}_k(t)\geq 0,\quad\text{for all }k\in\{5,6,7,8\}.\nonumber
\end{align}
When we consider the DCP \eqref{DCP_11} path-wise, we have the following optimization problem:
\begin{align}
& \min\quad h_a \left(\ot{Q}_5(\omega_t)-\ot{Q}_4(\omega_t)\right)^+ +h_b \left(\ot{Q}_6(\omega_t)\vee\ot{Q}_7(\omega_t)\right) + h_c \left(\ot{Q}_8(\omega_t)-\ot{Q}_9(\omega_t)\right)^+,\nonumber\\
& \hspace{2cm} \text{s.t.}\quad\;  \ot{Q}_5(\omega_t) +\frac{\mu_{A}}{\mu_{B1}} \ot{Q}_6(\omega_t) =  \ot{W}_5(\omega_t), \label{DCP_12}\\
& \hspace{3cm} \ot{Q}_7(\omega_t) +\frac{\mu_{B2}}{\mu_{C}} \ot{Q}_8(\omega_t) = \ot{W}_6(\omega_t), \nonumber\\
& \hspace{3cm} \ot{Q}_k(\omega_t)\geq 0,\quad\text{for all }k\in\{5,6,7,8\},\nonumber
\end{align}
where  $\ot{Q}_k(\omega_t)$, $k\in\{5,6,7,8\}$ are the decision variables. Note that, since servers 4 and 7 are independent of control, we neglect the term $h_a\ot{Q}_4(\omega_t)+h_c\ot{Q}_9(\omega_t)$ in the objective function of the optimization problem \eqref{DCP_12}. Next, we will solve the optimization problem \eqref{DCP_12} case by case.

First, let us consider the case $h_b \mu_{B1} \geq h_a \mu_A $ and $h_b \mu_{B2} \geq h_c \mu_C$. We present the optimal solution in the following result.

\begin{lemma}\label{three_type_lemma}
When $h_b \mu_{B1} \geq h_a \mu_A $ and $h_b \mu_{B2} \geq h_c \mu_C$, an optimal solution of the optimization problem \eqref{DCP_12} is the following:
\begin{enumerate}
\item If $h_b \geq h_a \mu_{A}/ \mu_{B1} + h_c \mu_{C} / \mu_{B2}$, an optimal solution is 
\begin{equation}\label{DCP_sln_4}
\left(\ot{Q}_5,\ot{Q}_6,\ot{Q}_7,\ot{Q}_8\right) = \left(\ot{W}_5,\bm{0},\bm{0},  \frac{\mu_{C}}{\mu_{B2}}\ot{W}_6 \right).
\end{equation}

\item If $h_b < h_a \mu_{A}/ \mu_{B1} + h_c \mu_{C} / \mu_{B2}$, an optimal solution is 
\begin{subequations}\label{DCP_sln_5}
\begin{align}
& \ot{Q}_5= \max \left\{ \ot{Q}_4\wedge \ot{W}_5,\;  \ot{W}_5 - \frac{\mu_{A}}{\mu_{B1}}\left(\ot{W}_6 - \frac{\mu_{B2}}{\mu_{C}} \ot{Q}_9\right)^+ \right\},\label{DCP_sln_5_1}\\
& \ot{Q}_6 = \ot{Q}_7 =  \min \left\{ \frac{\mu_{B1}}{\mu_{A}} \left( \ot{W}_5 -  \ot{Q}_4\right)^+,\; \left( \ot{W}_6 -  \frac{\mu_{B2}}{\mu_{C}} \ot{Q}_9\right)^+\right\},\label{DCP_sln_5_2}\\
& \ot{Q}_8= \max \left\{  \ot{Q}_9\wedge \frac{\mu_{C}}{\mu_{B2}} \ot{W}_6 ,\; \frac{\mu_{C}}{\mu_{B2}} \left( \ot{W}_6 -  \frac{\mu_{B1}}{\mu_{A}}  \left(\ot{W}_5 -  \ot{Q}_4\right)^+\right)\right\}.\label{DCP_sln_5_3}
\end{align}
\end{subequations}
\end{enumerate}
\end{lemma}
The proof of Lemma \ref{three_type_lemma} is provided in Appendix \ref{three_type_lemma_proof}. Note that the optimal solution has different structures depending on the cost parameters. On the one hand, $h_b \geq h_a \mu_{A}/ \mu_{B1} + h_c \mu_{C} / \mu_{B2}$ intuitively implies that the holding cost of a type $b$ job is greater than the sum of the holding costs of a type $a$ and type $c$ job. Hence, $(\ot{Q}_6,\ot{Q}_7) =(\bm{0},\bm{0})$ in \eqref{DCP_sln_4} and this implies that we should give static preemptive priority to type $b$ jobs in servers 5 and 6 all the time.

On the other hand, $\max\{h_a \mu_{A}/ \mu_{B1}, h_c \mu_{C} / \mu_{B2}\} \leq h_b < h_a \mu_{A}/ \mu_{B1} + h_c \mu_{C} / \mu_{B2}$ intuitively implies that the holding cost of a type $b$ job is greater than the one of a type $a$ or type $c$ job but less than the sum of the holding costs of a type $a$ and type $c$ job. In this case we interpret the solution \eqref{DCP_sln_5} in the following way. If we decrease both $\ot{Q}_6$ and $\ot{Q}_7$ in \eqref{DCP_sln_5_2}, which is equivalent to giving more priority to type $b$ jobs, then jobs will accumulate in buffers 10 and 15 waiting for the corresponding jobs to arrive buffers 11 and 14, respectively. Hence, the throughput rate of type $b$ jobs will increase but the throughput rate of both type $a$ and $c$ jobs will decrease, which is not desired because $h_b < h_a \mu_{A}/ \mu_{B1} + h_c \mu_{C} / \mu_{B2}$. If we decrease only $\ot{Q}_6$ in \eqref{DCP_sln_5_2}, then the throughput rate of type $a$, $b$, and $c$ jobs decrease, stays the same, and stays the same, respectively; and this is not desired. The case in which we decrease only $\ot{Q}_7$ in \eqref{DCP_sln_5_2} follows similarly. Next, suppose that we increase $\ot{Q}_6(t)$ or $\ot{Q}_7(t)$ in \eqref{DCP_sln_5_2}, which is equivalent to giving less priority type to $b$ jobs. Suppose that $\ot{Q}_6(t) = \ot{Q}_7(t) = (\mu_{B1}/\mu_{A}) \big( \ot{W}_5(t) -  \ot{Q}_4(t)\big)^+$ in \eqref{DCP_sln_5_2} for some $t\in\R_+$. Since $\ot{Q}_{10}(t)=\big(\ot{Q}_5(t)-\ot{Q}_4(t)\big)^+$, by substitution $\ot{Q}_{10}(t)=0$ in this case, i.e., throughput rate of type $a$ jobs is in the maximum level. Then, when we increase $\ot{Q}_6(t)$ or $\ot{Q}_7(t)$, the throughput rate of type $a$, $b$, and $c$ jobs stays the same, decreases, and (in the best case) increases, respectively. Since $h_c \mu_{C} / \mu_{B2} \leq h_b$, this result is not desired. The case $\ot{Q}_6(t) = \ot{Q}_7(t) = \big( \ot{W}_6(t) - (\mu_{B2}/\mu_{C}) \ot{Q}_9(t)\big)^+$ in \eqref{DCP_sln_5_2} follows similarly (for more intuition see the proof of Lemma \ref{three_type_lemma} and Remark \ref{type_3_intuition} in Appendix \ref{three_type_lemma_proof}). 

Therefore, \eqref{DCP_sln_5} is the optimal solution and we interpret the following policy from it: Whenever \eqref{DCP_sln_5_2} does not hold, if $\ot{Q}_6$ ($\ot{Q}_7$) is strictly greater than the right hand side of \eqref{DCP_sln_5_2}, then server 5 (6) should give preemptive priority to type $b$ jobs; if $\ot{Q}_6$ ($\ot{Q}_7$) is less than or equal to the right hand side of \eqref{DCP_sln_5_2}, then server 5 (6) should give preemptive priority to type $a$ ($c$) jobs.

Second, let us consider the case $h_a \mu_A \geq h_b \mu_{B1}$ and $h_b \mu_{B2} \geq h_c \mu_C$. We can see that 
\begin{align}
&\left(\ot{Q}_5,\ot{Q}_6,\ot{Q}_7,\ot{Q}_8\right) = \Bigg(\ot{Q}_4\wedge\ot{W}_5,\; \frac{\mu_{B1}}{\mu_{A}}\left(\ot{W}_5-\ot{Q}_4\right)^+, \nonumber\\
&\hspace{4cm} \frac{\mu_{B1}}{\mu_{A}}\left(\ot{W}_5-\ot{Q}_4\right)^+ \wedge \ot{W}_6,\; \frac{\mu_{C}}{\mu_{B2}}\left(\ot{W}_6- \frac{\mu_{B1}}{\mu_{A}}\left(\ot{W}_5-\ot{Q}_4\right)^+ \right)^+ \Bigg) \label{DCP_sln_2}
\end{align}
is an optimal solution of the optimization problem \eqref{DCP_12} for all $t\geq0$ and $\omega$ in $\Omega$ except a null set. Since the derivation of \eqref{DCP_sln_2} is very similar to the proof of Lemma \ref{three_type_lemma}, we skip it.

We can interpret a control policy from \eqref{DCP_sln_2} in the following way. If server 4 processes type $a$ jobs with a faster rate than server 5 does, then server 5 should give static priority to type $a$ jobs. Otherwise, server 5 should use the SDP policy to pace the departure process of the type $a$ jobs from buffer 5 with the one from buffer 4. Similarly, server 6 should use the SDP policy to pace the departure process of type $b$ jobs from buffer 7 with the one from buffer 6.

Third, let us consider the case $h_a \mu_A \geq h_b \mu_{B1}$ and $h_c \mu_C \geq h_b \mu_{B2}$. We can see that 
\begin{equation}\label{DCP_sln_3}
\left(\ot{Q}_5,\ot{Q}_6,\ot{Q}_7,\ot{Q}_8\right) = \left(\ot{Q}_4\wedge\ot{W}_5, \frac{\mu_{B1}}{\mu_{A}}\left(\ot{W}_5-\ot{Q}_4\right)^+, \left(\ot{W}_6-\frac{\mu_{B2}}{\mu_{C}}\ot{Q}_9\right)^+, \frac{\mu_{C}}{\mu_{B2}}\ot{W}_6\wedge\ot{Q}_9 \right)
\end{equation}
is an optimal solution of the optimization problem \eqref{DCP_12}. Since the derivation of \eqref{DCP_sln_3} is very similar to the proof of Lemma \ref{three_type_lemma}, we skip it.

We can interpret a control policy from \eqref{DCP_sln_3} in the following way. If server 4 (7) processes type $a$ ($c$) jobs with a faster rate than server 5 (6) does, then server 5 (6) should give static priority to type $a$ ($c$) jobs. Otherwise, server 5 (6) should use the SDP policy to pace the departure process of type $a$ ($c$) jobs from buffer 5 (8) with the one from buffer 4 (9).

Lastly, the case $h_b \mu_{B1} \geq h_a \mu_A$ and $h_c \mu_C \geq h_b \mu_{B2}$ is equivalent to the second case that we consider, hence we skip it.

It is possible to construct DCPs for networks with more than three job types and more than two shared servers by the same methodology that we use to construct DCP \eqref{DCP_11}. However, finding a closed-form optimal solution and interpreting a heuristic control policy is not trivial. As the number of shared servers increases in the network, the dimension of the corresponding DCP, which is the number of workload constraints, increases and finding a closed-form optimal solution becomes quite challenging. For example, there is only a single shared server in the networks considered in Figure \ref{fj_network} and in Sections \ref{task_dep}, \ref{ex_fork}, and \ref{ex_type_1}. As a result, all of the corresponding DCPs (cf. \eqref{DCP_4}, \eqref{DCP_5}, \eqref{DCP_13}, and \eqref{DCP_7}) are single dimensional, i.e., each DCP has a single workload constraint. However, when there are two shared servers as in the network in Figure \ref{fj_network_3}, the corresponding DCP \eqref{DCP_11} has two workload constraints and is more difficult to solve than the single dimensional ones.

In the queueing control literature, there are many studies (cf. \citet{har98, har99, bel01,man04, ata05, dai08}) in which the corresponding DCP is first converted to a single dimensional equivalent DCP and control policies are obtained by solving the latter DCP. However, in general, how to control stochastic networks whose equivalent lower dimensional DCP has more than one dimension is an open problem.

\begin{figure}[htb]
\begin{center}
\includegraphics[width=0.9\textwidth]{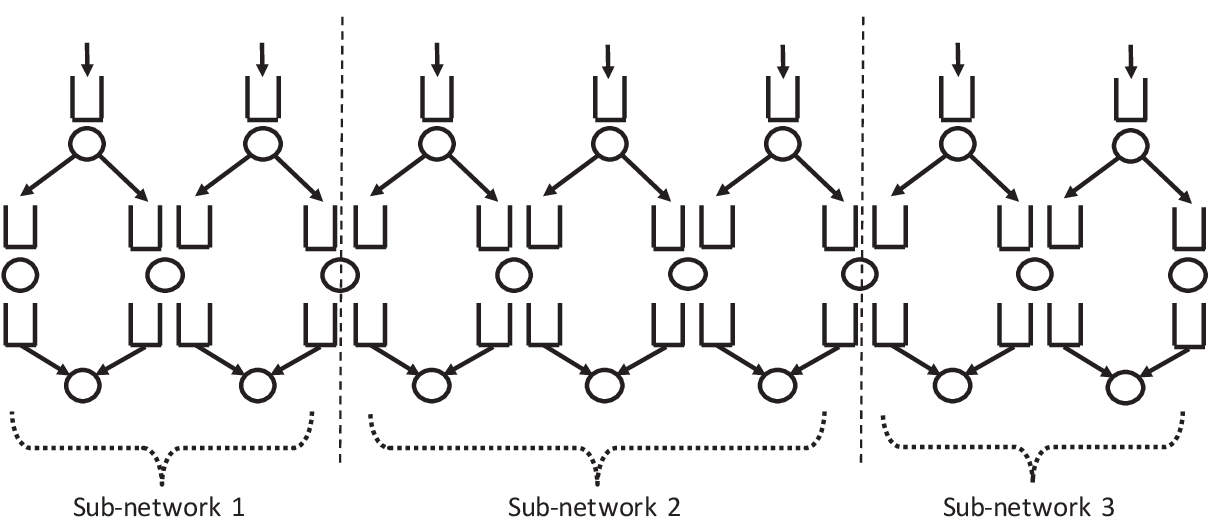}
\caption{A fork-join processing network with seven job types and six shared servers.}\label{fj_network_4}
\end{center}
\end{figure}

In order to overcome the curse of dimensionality in the networks with more than three job types and more than two shared servers, one possible solution is to decompose the network into the ones with at most three job types and use the control policies mentioned in this study to control each sub-network. For example, in Figure \ref{fj_network_4}, there are seven different job types arriving to the system, and the network is divided into three sub-networks. Then the important questions are: 1) How can a large network be decomposed into smaller ones? and 2) Which control policies should be used in the servers which belong to two different sub-networks? If there exists a server which processes more than two different job types and is in light traffic, the network can clearly be decomposed from this server and any work-conserving control policy can be used in this server, because its processing capacity is high. Otherwise, the answers are not trivial and requires further research.

\section{Concluding Remarks}\label{conclusion}

The synchronization constraints in fork-join networks complicate their analysis. In this paper, we have formulated and solved approximating diffusion control problems (DCPs) for a variety of fork-join networks (depicted in Figures \ref{fj_network_0}, \ref{fj_network}, \ref{fj_network_2}, \ref{fj_network_5}, and \ref{fj_network_3}). The DCP solutions have motivated either a slow departure pacing control or a static priority control, depending on the network parameters. We have rigorously proved the asymptotic optimality of our proposed control for the network in Figure \ref{fj_network}, and we conjecture that the proofs will extend in a straightforward manner to the networks in Figures \ref{fj_network_0}, \ref{fj_network_2}, \ref{fj_network_5}, and \ref{fj_network_3} and the network in Figure \ref{fj_network} with task dependent holding costs.   

Our proposed controls will not change when light traffic queues\footnote{What we mean by a queue is a server and its corresponding buffer, e.g., server 3 and buffer 3 is a queue in Figure \ref{fj_network}.} are added to the fork-join networks that we have studied. Furthermore, minor modifications to our proposed policies will accommodate the addition of heavy traffic queues before the shared or dedicated servers.  For example, if we add a heavy traffic queue between server 1 and buffer 3 in Figure \ref{fj_network}, even though our proofs do not cover this case, we can see that the proposed policy (or a slightly modified version of it) is still asymptotically optimal. In this modified system, if server 3 is in light traffic, then the modified network is equivalent to the one where server 3 is replaced by the newly added heavy traffic queue and all the other servers have the same processing capacities with the ones in the original network. If server 3 is in heavy traffic in this modified system, then the proposed policy behaves in the following way: When the total number of jobs in buffer 3 and the buffer corresponding to the newly added queue is less than the one in buffer 4, then server 4 gives priority to buffer 4; otherwise server 4 gives priority to buffer 5. 

The complicated case is when there are heavy-traffic queues after the shared servers. Then, it is not clear either what the proposed policy should be or how to prove an asymptotic optimality result. An excellent topic for future research is to develop control policies for the broader class of fork-join networks with heterogeneous customer populations described in \citet{ngu94}. More specifically, that paper assumes FCFS scheduling, but we believe other control policies can lead to better performance.

\section*{Acknowledgements}

We would like to thank Itai Gurvich for proposing the fork-join network structure in Figure \ref{fj_network} as an interesting prototype network to study. We would like to thank the SAMSI 2013-2014 working group on patient flow for valuable discussions on the ideas in this paper. We would like to thank the associate editor and the three anonymous referees for their constructive comments, which improved the paper significantly.


\appendix
\section*{APPENDIX}

The appendix is organized into three sections: \ref{light_traffic_conv_rate_proof}, \ref{standard_results}, and \ref{detailed_numerical_results}. Section \ref{light_traffic_conv_rate_proof} derives the convergence rate of the queue length process in a light traffic $GI/GI/1$ queue, which is used in the weak convergence proof in the main body (the proof of Proposition \ref{light_traffic_conv_rate_lemma}). Section \ref{standard_results} presents the proofs of all results in the paper with standard methodology. Section \ref{detailed_numerical_results} provides the detailed results of all simulation experiments mentioned in Section \ref{numerical_analyses}.

\section{Proof of Proposition \ref{light_traffic_conv_rate_lemma}}\label{light_traffic_conv_rate_proof}

Without loss of generality, we will prove Proposition \ref{light_traffic_conv_rate_lemma} only for $r\in\N_+$. Let $U(0):=0$, $V(0):=0$,
\begin{align*}
U(k):=\sum_{i=1}^k u_i,&\quad\forall k\in\N_+,\quad\text{and}\quad A(t):=\max\{k\in\N:U(k)\leq t\},\quad\forall t\geq 0,\\
V(k):=\sum_{i=1}^k v_i,&\quad\forall k\in\N_+,\quad\text{and}\quad S(t):=\max\{k\in\N:V(k)\leq t\},\quad\forall t\geq 0.
\end{align*}
Let $T(t)$ ($I(t)$) denote the cumulative amount of time the server is busy (idle) up to time $t$, $t\geq 0$. Then, $T(t)+I(t)=t$ for all $t\geq 0$ and
\begin{equation*}
Q(t)= Q(0)+A(t)-S(T(t)),\quad\forall t\geq 0.
\end{equation*}
Let us define the following shifted processes. For all $t\geq 0$ and $i\in\N$,
\begin{align*}
I^{r,i}(t)&:= I(ir+t)-I(ir), &T^{r,i}(t)&:= T(ir+t)-T(ir),\\
A^{r,i}(t)&:= A(ir+t)-A(ir)-\lambda t, &S^{r,i}(t)&:= S(T(ir+t))-S(T(ir))-\mu T^{r,i}(t),\\
Q^{r,i}(t)&:= Q(ir+t), & X^{r,i}(t)&:= Q(ir)+A^{r,i}(t)-S^{r,i}(t)+(\lambda-\mu)t,
\end{align*}
Then, after some algebra, for all $t\geq 0$ and $i\in\N$, we have
\begin{equation*}
Q^{r,i}(t)=X^{r,i}(t) + \mu I^{r,i}(t).
\end{equation*}
For all $i\in\N$ and $r\in\N_+$, since $I^{r,i}(0)=0$, $I^{r,i}(\cdot)$ is nondecreasing, and $\int_0^\infty Q^{r,i}(t) \dr I^{r,i}(t)=0$, we have
\begin{equation*}
\left(Q^{r,i}, \mu I^{r,i}\right) =\left(\Phi,\Psi\right)(X^{r,i}).
\end{equation*}
 
Next, for all $r\in\N_+$, let us define the following sets:
\begin{align}
\mG^r &= \left\{ Q(0) \leq \frac{((\mu-\lambda)\wedge\epsilon)r}{2} \right\} \cap \bigcap_{i=0}^{r^{n-1}-1} \left\{ \left\Vert A^{r,i} \right\Vert_r \vee \left\Vert S^{r,i} \right\Vert_r  \leq \frac{((\mu-\lambda)\wedge\epsilon)r}{8} \right\},\label{good_g}\\
\mH^r &=  \bigcap_{i=0}^{r^{n-1}-1} \left\{ Q(ir) \leq \frac{((\mu-\lambda)\wedge\epsilon)r}{2}\right\}. \label{good_h}
\end{align}
We have the following result related to \eqref{good_g} whose proof is presented in Appendix \ref{proof_lemma_exp_conv}.

\begin{lemma}\label{exp_conv_r}
For all $\epsilon>0$, there exist $r_2\geq 0$ and strictly positive constants $C_4$, $C_5$, $C_6$, and $C_7$ which are independent of $r$ such that if $r\geq r_2$,
\begin{align}
&\pr \left(  \left\Vert A^{r,i} \right\Vert_r >\epsilon r \right) \leq C_4r^n\e^{-C_5r}, \label{exp_conv_r_1}\\
&\pr \left(  \left\Vert S^{r,i} \right\Vert_r >\epsilon r \right) \leq C_6r^n\e^{-C_7r}, \label{exp_conv_r_2}
\end{align}
for all $i\in\{0,1,\ldots,r^{n-1}-1\}$.
\end{lemma}

Next, we will show that $\mG^r \subset \mH^r$ by induction. First note in the set $\mG^r $, $\left\{ Q(0) \leq ((\mu-\lambda)\wedge\epsilon)r/2 \right\}$. Next, suppose that, in the set $\mG^r $, $\left\{ Q(jr) \leq ((\mu-\lambda)\wedge\epsilon)r/2 \right\}$ for all $j\in\{1,2,\ldots,i\}$ for some $i\in\{0,1,\ldots,r^{n-1}-2\}$. Then, in the set $\mG^r $,
\begin{align*}
&Q((i+1)r)=Q^{r,i}(r) = X^{r,i}(r) + \mu I^{r,i}(r),\\
&\quad = Q(ir)+A^{r,i}(r)-S^{r,i}(r)+(\lambda-\mu)r+ \sup_{0\leq t \leq r} \left( -Q(ir)-A^{r,i}(t)+S^{r,i}(t)+(\mu-\lambda)t\right)^+,\\
&\quad \leq Q(ir)+\frac{((\mu-\lambda)\wedge\epsilon)r}{4}+(\lambda-\mu)r+ \sup_{0\leq t \leq r} \left( -A^{r,i}(t)+S^{r,i}(t)\right)^+ + \sup_{0\leq t \leq r} \left( -Q(ir)+(\mu-\lambda)t\right)^+,\\
&\quad  \leq \frac{((\mu-\lambda)\wedge\epsilon)r}{2}+ Q(ir)+(\lambda-\mu)r + \sup_{0\leq t \leq r} \left( -Q(ir)+(\mu-\lambda)t\right)^+,\\
&\quad  = \frac{((\mu-\lambda)\wedge\epsilon)r}{2}+ Q(ir)+(\lambda-\mu)r + \left(-Q(ir)+(\mu-\lambda)r\right)\vee 0, \\
&\quad  = \frac{((\mu-\lambda)\wedge\epsilon)r}{2}+ Q(ir)+(\lambda-\mu)r -Q(ir)+(\mu-\lambda)r, \\
&\quad  = \frac{((\mu-\lambda)\wedge\epsilon)r}{2},
\end{align*}
where the first and second inequalities are by \eqref{good_g} and second to last equality is by the induction hypothesis. Therefore, $\mG^r \subset \mH^r$.

Next, note that the probability in left hand side of \eqref{l_conv_rate} is equal to 
\begin{align}
 \pr \left(  \bigcup_{i=0}^{r^{n-1}-1} \left\{ \sup_{0\leq t \leq r} Q^{r,i}(t)>r  \epsilon \right\} \right) &\leq \pr \left(  \bigcup_{i=0}^{r^{n-1}-1} \left\{\sup_{0\leq t \leq r}  Q^{r,i}(t)>r  \epsilon \right\} ,\;\mG^r\right)+ \pr \left( \left(\mG^r\right)^c\right), \nonumber\\
&= \pr \left(  \bigcup_{i=0}^{r^{n-1}-1} \left\{\sup_{0\leq t \leq r}  Q^{r,i}(t)>r  \epsilon \right\} ,\;\mG^r\cap\mH^r\right)+ \pr \left( \left(\mG^r\right)^c\right).\label{l_conv_eq_1}
\end{align}
First, let us consider the first probability in \eqref{l_conv_eq_1}. Let us fix an arbitrary $r\in\N_+$. In the set $\mG^r \cap \mH^r$, for all $i\in\{0,1,\ldots,r^{n-1}-1\}$, we have
\begin{align}
&\sup_{0\leq t \leq r}  Q^{r,i}(t) = \sup_{0\leq t \leq r} \left\{X^{r,i}(t) + \mu I^{r,i}(t) \right\}, \nonumber\\
&\quad =  \sup_{0\leq t \leq r} \left\{ Q(ir)+A^{r,i}(t)-S^{r,i}(t)+(\lambda-\mu)t+ \sup_{0\leq s \leq t} \left( -Q(ir)-A^{r,i}(s)+S^{r,i}(s)+(\mu-\lambda)s\right)^+ \right\}, \nonumber\\
&\quad  \leq \frac{((\mu-\lambda)\wedge\epsilon)r}{2}+  \sup_{0\leq t \leq r} \left\{  Q(ir)+(\lambda-\mu)t + \sup_{0\leq s \leq t} \left( -Q(ir)+(\mu-\lambda)s\right)^+ \right\},\nonumber\\
&\quad  =\frac{((\mu-\lambda)\wedge\epsilon)r}{2}+ \sup_{0\leq t \leq r}\Phi\left( Q(ir)+(\lambda-\mu)e \right)(t), \nonumber\\
&\quad  =\frac{((\mu-\lambda)\wedge\epsilon)r}{2}+Q(ir), \nonumber\\
&\quad  \leq ((\mu-\lambda)\wedge\epsilon)r,\label{l_conv_eq_2}
\end{align}
where the first and last inequalities are by \eqref{good_g} and \eqref{good_h}, respectively. Therefore,
\begin{equation}\label{l_conv_eq_3}
 \bigcup_{i=0}^{r^{n-1}-1} \left\{\sup_{0\leq t \leq r}  Q^{r,i}(t)>r  \epsilon \right\} \cap \mG^r\cap\mH^r = \emptyset,
\end{equation}
by \eqref{l_conv_eq_2}, which implies that the first probability in \eqref{l_conv_eq_1} is equal to 0.

Second, let us consider the second probability in \eqref{l_conv_eq_1}, which is less than or equal to 
\begin{align}
&\pr\left(Q(0) > \frac{((\mu-\lambda)\wedge\epsilon)r}{2} \right)\nonumber\\
&\hspace{1cm}+ \sum_{i=0}^{r^{n-1}-1}  \pr\left( \left\Vert A^{r,i} \right\Vert_r> \frac{((\mu-\lambda)\wedge\epsilon)r}{8} \right)+\sum_{i=0}^{r^{n-1}-1} \pr\left(\left\Vert S^{r,i} \right\Vert_r  > \frac{((\mu-\lambda)\wedge\epsilon)r}{8} \right),\label{l_conv_eq_4}
\end{align}
by \eqref{good_g}. Let $r_1:= r_0\vee r_2$. Then, by \eqref{l_in_q} and Lemma \ref{exp_conv_r}, if $r\geq r_1$, \eqref{l_conv_eq_4} is less than or equal to
\begin{equation}
C_0r^{2n-1}\e^{-C_1r} + C_4r^{2n-1}\e^{-C_5r} + C_6r^{2n-1}\e^{-C_7r} =  C_2r^{2n-1}\e^{-C_3r},\label{l_conv_eq_5}
\end{equation}
where $C_2:=C_0+C_4+C_6$ and $C_3:=\min\{C_1,C_5,C_7\}$. Therefore, \eqref{l_conv_eq_1}, \eqref{l_conv_eq_3}, and \eqref{l_conv_eq_5} prove Proposition \ref{light_traffic_conv_rate_lemma}.

\subsection{Proof of Lemma \ref{exp_conv_r}}\label{proof_lemma_exp_conv}

We first present some preliminary results in Appendix \ref{exp_conv_r_prem}, second prove \eqref{exp_conv_r_2} in Appendix \ref{exp_conv_r_2_proof}, and lastly prove \eqref{exp_conv_r_1} in Appendix \ref{exp_conv_r_1_proof}. Let us fix an arbitrary $i\in\{0,1,\ldots,r^{n-1}-1\}$. Without loss of generality, we choose $\epsilon>0$ such that $\epsilon<\mu$.

\subsubsection{Preliminary Results}\label{exp_conv_r_prem}

Note that 
\begin{align}
\left\Vert S^{r,i} \right\Vert_r = \sup_{0\leq t\leq r} \left| S^{r,i}(t) \right| &= \sup_{0\leq t\leq r} \left| S(T(ir+t))-S(T(ir))-\mu T(ir+t) + \mu T(ir) \right|,\nonumber\\
&\leq  \sup_{0\leq t\leq r} \left| S(T(ir)+t)-S(T(ir))-\mu (T(ir)+ t) + \mu T(ir) \right|,\nonumber\\
&=  \sup_{0\leq t\leq r} \left| S(T(ir)+t)-S(T(ir))-\mu t \right|,\label{exp_conv_eq_1}
\end{align}
where the first inequality is by the fact that for all $t\geq 0$, there exists an $s\in[0,t]$ such that $T(ir+t)=T(ir)+s$, because $T(\cdot)$ is nondecreasing and $T(t)\leq t$ for all $t\geq 0$. For all $m\in\N$, let
\begin{equation*}
\mF_m :=\sigma\left\{ u_k, k\in\N_+;\; v_k, k\in\{1,2,\ldots,m+1\}\right\}.
\end{equation*}
Clearly, $\{\mF_m,m\in\N\}$ is a filtration. Then, similar to Lemma 8.3 of \citet{wil98a} and Lemmas 7.5 and 7.6 of \citet{bel01}, we have the following result.

\begin{lemma}\label{st_lemma}
For all $t\geq 0$, $S(T(t))$ is a stopping time with respect to the filtration $\{\mF_m,m\in\N\}$. The $\sigma$-algebra associated with the stopping time $S(T(t))$ is 
\begin{equation*}
\mF_{S(T(t))} := \left\{A\in\mF: A\cap \{S(T(t))=m\}\in\mF_m, \forall m\in\N \right\}.
\end{equation*}
Then, $v_{S(T(t))+1}\in\mF_{S(T(t))}$. Let $\{A_k,k\in\N_+\}$ be a sequence of sets such that $A_k\in\mB(\R)$ for all $k\in\N_+$, and $B$ be a set such that $B\in\mF_{S(T(t))}$. Then,
\begin{align*}
&\pr \left( v_{S(T(t))+k}\in A_k, k\in\{2,3,\ldots\},\;B,\; S(T(t)) <\infty \right) \\
&\hspace{3cm}= \pr \left( v_{k-1}\in A_k, k\in\{2,3,\ldots\}\right)\times\pr\left( B,\; S(T(t)) <\infty \right)
\end{align*}
\end{lemma}
Lemma \ref{st_lemma} can be proven similarly to Lemma 8.3 of \citet{wil98a} and Lemma 7.6 of \citet{bel01} are proven, hence we skip the proof. 

Next, we present some preliminary results related to large deviations theory. For all $\alpha>0$, let 
\begin{equation}\label{exp_conv_eq_2}
\ell(\alpha):=\ln \E\left[ \exp\left\{\alpha\left(\frac{1}{\mu}- v_1\right)\right\} \right].
\end{equation}
Then, $e^{\ell(\alpha)}<\infty$ for all $\alpha\in(-\ol{\alpha},\ol{\alpha})$ by the exponential moment assumption on $v_1$. For $x\geq 0$, let
\begin{equation}\label{exp_conv_eq_3}
\Lambda_1(x):=\sup_{\alpha>0}\left\{ \alpha x-\ell(\alpha)\right\},\quad\text{and}\quad \Lambda_2(x):=\sup_{\alpha>0}\left\{ \alpha x-\ell(-\alpha)\right\}.
\end{equation}
Note that $\Lambda_1$ and $\Lambda_2$ are not exactly but very similar to Fenchel-Legendre transform of $\ell(\alpha)$ and $\ell(-\alpha)$, respectively (cf. Definition 2.2.2 of \citet{dem98}). Then, we have the following result.

\begin{lemma}\label{fl_transform}
Both $\Lambda_1$ and $\Lambda_2$ are convex and nondecreasing in $[0,\infty)$, $\Lambda_1(0)=\Lambda_2(0)=0$, and $\Lambda_1(x)>0$ and $\Lambda_2(x)>0$ for all $x> 0$. 
\end{lemma}

\begin{proof}
First, let us consider $\Lambda_1$. Let 
\begin{equation}\label{le_1}
\Lambda_1^*(x):=\sup_{\alpha\in\R}\left\{ \alpha x-\ell(\alpha)\right\},\quad\forall x\in\R,
\end{equation}
denote the Fenchel-Legendre transform of $\ell(\alpha)$ (cf. Definition 2.2.2 of \citet{dem98}). Then, by Parts (a) and (b) of Lemma 2.2.5 of \citet{dem98}, $\Lambda_1^*$ is convex, $\Lambda_1^*(0)=0$, $\Lambda_1^*$ is nondecreasing in $[0,\infty)$, and for all $x\geq 0$,
\begin{equation}\label{le_2}
\Lambda_1^*(x)=\sup_{\alpha\geq 0}\left\{ \alpha x-\ell(\alpha)\right\} = \Lambda_1(x) \vee 0,
\end{equation}
where the last equality in \eqref{le_2} is by \eqref{exp_conv_eq_3}. Moreover, by Parts (a) and (c) of Lemma 2.2.5 of \citet{dem98}, $\ell$ is convex in $\R$, $\ell$ is differentiable in $(-\ol{\alpha},\ol{\alpha})$, and $\ell(0)=\ell'(0)=0$, where $\ell'$ is the derivative of $\ell$. Then, $\ell$ achieves the global minimum at point $0$; and since it is convex, $\ell$ is nondecreasing in the interval $[0,\infty)$. Therefore, for any fixed $x>0$, there exists an $\alpha^*\in(0,\ol{\alpha})$ such that $\Lambda_1^*(x)\geq \alpha^*x-\ell(\alpha^*)>0$. As a result, $\Lambda_1(x) = \Lambda_1^*(x) >0$ for all $x>0$ by \eqref{le_2}. Lastly, since $\ell$ is convex and achieves its global minimum at point $0$ in its domain $\R$ and continuous in the interval $(-\ol{\alpha},\ol{\alpha})$, 
\begin{equation*}
\Lambda_1(0)= -\inf_{\alpha>0} \ell(\alpha) = -\ell(0) = 0 = \Lambda_1^*(0),
\end{equation*}
where the first equality is by \eqref{exp_conv_eq_3}. Hence, $\Lambda_1(x) = \Lambda_1^*(x)$ for all $x\geq 0$, which implies that $\Lambda_1$ is convex and nondecreasing in $[0,\infty)$, $\Lambda_1(0)=0$, and $\Lambda_1(x)>0$ for all $x>0$.

The proofs for $\Lambda_2$ follows with exactly the same way. The only difference is that we consider the random variable $v_1-1/\mu$ instead of $1/\mu-v_1$.
\end{proof}

Next, we will derive three more preliminary results.

\begin{lemma}\label{ld_result_1}
For all $\epsilon>0$, there exists an $r_3\geq 0$ such that if $r\geq r_3$, we have
\begin{equation*}
\pr \left( S(r) > (\mu+\epsilon) r \right) \leq C_8\e^{-C_9 r},
\end{equation*}
where $C_8$ and $C_9$ are strictly positive constants independent of $r$.
\end{lemma}

\begin{proof}
\begin{align}
&\pr \left( S(r) > (\mu+\epsilon) r \right) \leq \pr \left( S(r) > \rd{(\mu+\epsilon) r} \right)=  \pr \left( V(\rd{(\mu+\epsilon) r})<r \right),\nonumber\\
&\quad =  \pr \left( \sum_{j=1}^{\rd{(\mu+\epsilon)r}} \left(\frac{1}{\mu}-v_j \right) > \frac{\rd{(\mu+\epsilon)r}}{\mu} -r\right) \leq \pr \left( \alpha\sum_{j=1}^{\rd{(\mu+\epsilon)r}} \left(\frac{1}{\mu}-v_j \right) > \alpha \frac{\epsilon r-1}{\mu} \right),\label{ld_result_eq_1}\\
&\quad = \pr \left( \exp\left\{\alpha \sum_{j=1}^{\rd{(\mu+\epsilon)r}} \left(\frac{1}{\mu}-v_j \right) \right\}> \exp\left\{\alpha \frac{\epsilon r-1}{\mu} \right\}\right),\nonumber\\
&\quad \leq \E\left[  \exp\left\{\alpha \sum_{j=1}^{\rd{(\mu+\epsilon)r}} \left(\frac{1}{\mu}-v_j \right) \right\}\right] \exp\left\{-\alpha \frac{\epsilon r-1}{\mu} \right\},\label{ld_result_eq_2}\\
& \quad= \prod_{j=1}^{\rd{(\mu+\epsilon)r}}  \E\left[  \exp\left\{\alpha \left(\frac{1}{\mu}-v_j \right) \right\}\right] \exp\left\{-\alpha \frac{\epsilon r-1}{\mu} \right\} = \exp\left\{\rd{(\mu+\epsilon)r}\ell(\alpha)-\alpha \frac{\epsilon r-1}{\mu}\right\},\label{ld_result_eq_3}\\
&\quad= \exp\left\{-\rd{(\mu+\epsilon)r}\left(\alpha \frac{\epsilon r-1}{\mu\rd{(\mu+\epsilon)r}}-\ell(\alpha)\right)\right\},\label{ld_result_eq_4}
\end{align}
where $\alpha>0$ is an arbitrary constant in \eqref{ld_result_eq_1}, \eqref{ld_result_eq_2} is by Markov's inequality, and the equalities in \eqref{ld_result_eq_3} are by the i.i.d. property of the sequence $\{v_j,j\in\N_+\}$ and \eqref{exp_conv_eq_2}. Since \eqref{ld_result_eq_4} is valid for all $\alpha>0$,
\begin{align}
\pr \left( S(r) > (\mu+\epsilon) r \right) &\leq  \exp\left\{-\rd{(\mu+\epsilon)r}\sup_{\alpha>0}\left(\alpha \frac{\epsilon r-1}{\mu\rd{(\mu+\epsilon)r}}-\ell(\alpha)\right)\right\},\nonumber\\
&=   \exp\left\{-\rd{(\mu+\epsilon)r}\Lambda_1\left( \frac{\epsilon r-1}{\mu\rd{(\mu+\epsilon)r}}\right)\right\},\label{ld_result_eq_5}
\end{align}
by \eqref{exp_conv_eq_3}. Note that there exists an $r_3>0$ and $\epsilon_1>0$ such that if $r\geq r_3$, we have
\begin{equation*}
 \frac{\epsilon r-1}{\mu\rd{(\mu+\epsilon)r}}>\epsilon_1.
 \end{equation*}
By Lemma \ref{fl_transform}, $\Lambda_1(\epsilon_1)>0$ and $\Lambda_1(x)$ is nondecreasing for all $x\geq 0$, thus the term in \eqref{ld_result_eq_5} converges to 0 with exponential rate as $r\rightarrow\infty$. To complete the proof, let $C_8:=\exp(\Lambda_1(\epsilon_1))$ and $C_9:=(\mu+\epsilon)\Lambda_1(\epsilon_1)$. 
\end{proof}

\begin{lemma}\label{ld_result_2}
For all $t\in\R_+$, $\pr \left( S(t) =\infty \right) =0$.
\end{lemma}

\begin{proof}
Fix arbitrary $t\in\R_+$ and $0<\epsilon<1/\mu$. Then,
\begin{equation*}
\pr \left( S(t) =\infty \right) = \pr \left( V(k) <t,\;\forall k\in\N_+\right)= \pr \left( \frac{V(k)}{k} <\frac{t}{k},\;\forall k\in\N_+\right)=0,
\end{equation*}
because for each $\omega\in\Omega$ except a null set, there exists a $k_0\in\N_+$ such that if $k\geq k_0$, $V(k,\omega)/k > 1/\mu-\epsilon$ by SLLN. When $k$ is sufficiently large, $V(k,\omega)/k > 1/\mu-\epsilon>t/k$ for almost all $\omega\in\Omega$, which completes the proof. 
\end{proof}

\begin{lemma}\label{ld_result_3}
For all $\epsilon>0$ and $r> 0$, we have
\begin{equation*}
\pr \left( v_1-\frac{1}{\mu} > \epsilon r \right) \leq C_{10}\e^{-C_{11} r},
\end{equation*}
where $C_{10}$ and $C_{11}$ are strictly positive constants independent of $r$.
\end{lemma}

\begin{proof}
\begin{align}
 \pr \left( v_1-\frac{1}{\mu} > \epsilon r \right) &= \pr \left( \frac{\bar{\alpha}}{2}\left(v_1-\frac{1}{\mu}\right) >  \frac{\bar{\alpha}}{2}\epsilon r \right) = \pr \left( \exp\left\{ \frac{\bar{\alpha}}{2}\left(v_1-\frac{1}{\mu}\right)\right\} > \exp\left\{ \frac{\bar{\alpha}\epsilon r}{2}\right\} \right), \nonumber\\
&\leq \E\left[ \exp\left\{ \frac{\bar{\alpha}}{2}\left(v_1-\frac{1}{\mu}\right)\right\}\right]\exp\left\{- \frac{\bar{\alpha}\epsilon r}{2}\right\}=\exp\left\{\ell\left(- \frac{\bar{\alpha}}{2}\right)-\frac{\bar{\alpha}\epsilon r}{2}\right\},\label{ld_result_eq_7}
\end{align}
where the inequality in \eqref{ld_result_eq_7} is by Markov's inequality, and the equality in \eqref{ld_result_eq_7} is by \eqref{exp_conv_eq_2}. Since $\ell\left(- \bar{\alpha}/2\right)<\infty$ because of the exponential moment assumption, \eqref{ld_result_eq_7} gives us the desired result. To complete the proof, let $C_{10}:=\exp\{\ell(- \bar{\alpha}/2)\}$ and $C_{11}:=\bar{\alpha}\epsilon/2$. 
\end{proof}

Next, we will prove Lemma \ref{exp_conv_r} by Lemmas \ref{st_lemma}, \ref{fl_transform}, \ref{ld_result_1}, \ref{ld_result_2}, and \ref{ld_result_3} in the following section.

\subsubsection{Proof of \eqref{exp_conv_r_2}}\label{exp_conv_r_2_proof}

Let 
\begin{equation*}
\eta^r:=\inf \left\{ t\in [0,r]: \left| S(T(ir)+t)-S(T(ir))-\mu t \right| > \epsilon r\right\},
\end{equation*}
where $\inf \{\emptyset\}:= \infty$ for completeness. Then,
\begin{equation}\label{exp_conv_eq_4}
\pr \left( \left\Vert S^{r,i}\right\Vert_r > \epsilon r \right)\leq \pr \left( \sup_{0\leq t\leq r }\left| S(T(ir)+t)-S(T(ir))-\mu t \right| > \epsilon r \right) = \pr \left(\eta^r \leq r \right),
\end{equation}
by \eqref{exp_conv_eq_1}. Let 
\begin{equation}\label{exp_conv_eq_4_0}
\ot{V}^{i,r} (k) := \sum_{j=S(T(ir))+2}^{S(T(ir))+k} v_j,\quad \forall k\in\{2,3,\ldots\},\qquad \ot{V}^{i,r}_1 (k) := \sum_{j=S(T(ir))+1}^{S(T(ir))+k} v_j,\quad \forall k\in\N_+,
\end{equation}
$\ot{V}^{i,r} (1):=0$ and $\ot{V}^{i,r} (k)=\ot{V}^{i,r}_1 (k):=0$ for all $k\in\{\ldots,-2,-1,0\}$. Then, for all $t\geq 0$,
\begin{align*}
& \left\{ \left| S(T(ir)+t)-S(T(ir))-\mu t \right| > \epsilon r\right\} \\
&\hspace{1cm}= \left\{  S(T(ir)+t)-S(T(ir))> \mu t + \epsilon r\right\} \cup  \left\{ S(T(ir)+t)-S(T(ir))< \mu t - \epsilon r\right\} \\
&\hspace{1cm}\subseteq \left\{  S(T(ir)+t)-S(T(ir))> \rd{\mu t + \epsilon r}\right\} \cup  \left\{ S(T(ir)+t)-S(T(ir))< \ru{\mu t - \epsilon r} \right\} \\
&\hspace{1cm} \subseteq   \left\{ \ot{V}^{i,r} (\rd{\mu t + \epsilon r}) < t \right\} \cup  \left\{\ot{V}^{i,r}_1 ( \ru{\mu t - \epsilon r}) > t \right\}.
\end{align*}
Let us define
\begin{align*}
&\eta_1^r:=\inf \left\{ t\in [0,r]: \ot{V}^{i,r}\left(\rd{ \mu t+ \epsilon r } \right) < t \right\},\\
&\eta_2^r:=\inf \left\{ t\in [0,r]: \ot{V}^{i,r}_1\left(\ru{\mu t- \epsilon r} \right) > t \right\}.
\end{align*}
Then, $\eta^r \geq \eta_1^r\wedge \eta_2^r$, so
\begin{equation}\label{exp_conv_eq_5}
\pr \left(\eta^r \leq r\right)\leq \pr \left(\eta_1^r\wedge \eta_2^r \leq r \right) \leq  \pr \left(\eta_1^r \leq r \right) +  \pr \left( \eta_2^r\leq r\right).
\end{equation}

First, 
\begin{align}
\pr \left(\eta_1^r< r \right) &=  \pr \left( \inf_{0\leq t\leq r} \left(\ot{V}^{i,r}\left(\rd{ \mu t+ \epsilon r } \right) - t\right) <0 \right)\nonumber\\
&\leq  \pr \left( \min_{j=\rd{\epsilon r},\ldots, \rd{ (\mu+ \epsilon) r}} \left(\ot{V}^{i,r}\left(j \right) - \frac{j+1-\epsilon r}{\mu}\right) <0 \right)\nonumber\\
&=  \pr \left( \min_{j=\rd{\epsilon r},\ldots, \rd{ (\mu+ \epsilon) r}} \left(\ot{V}^{i,r}\left(j \right) - \frac{j-1}{\mu}\right) < \frac{2-\epsilon r}{\mu} \right)\nonumber\\
&\leq  \pr \left( \min_{j=\rd{\epsilon r},\ldots, \rd{ (\mu+ \epsilon) r}} \left(\ot{V}^{i,r}\left(j \right) - \frac{j-1}{\mu}\right) < \frac{2-\epsilon r}{\mu},\; S(T(ir))<\infty  \right)+\pr \left( S(T(ir))=\infty\right)\nonumber\\
&=  \pr \left( \min_{j=\rd{\epsilon r},\ldots, \rd{ (\mu+ \epsilon) r}} \left(\ot{V}^{i,r}\left(j \right) - \frac{j-1}{\mu}\right) < \frac{2-\epsilon r}{\mu},\; S(T(ir))<\infty  \right)\label{exp_conv_eq_6}\\
&=  \pr \left( \min_{j=\rd{\epsilon r},\ldots, \rd{ (\mu+ \epsilon) r}} \left(V\left(j-1 \right) - \frac{j-1}{\mu}\right) < \frac{2-\epsilon r}{\mu},\; S(T(ir))<\infty  \right)\label{exp_conv_eq_7}\\
&=  \pr \left( \min_{j=\rd{\epsilon r}-1,\ldots, \rd{ (\mu+ \epsilon) r}-1} \left(V\left(j \right) - \frac{j}{\mu}\right) < \frac{2-\epsilon r}{\mu}  \right)\nonumber\\
&=  \pr \left( -\left(\min_{j=\rd{\epsilon r}-1,\ldots, \rd{ (\mu+ \epsilon) r}-1} \left(V\left(j \right) - \frac{j}{\mu}\right)\right) > \frac{\epsilon r-2}{\mu} \right)\nonumber\\
&=  \pr \left( \max_{j=\rd{\epsilon r}-1,\ldots, \rd{ (\mu+ \epsilon) r}-1} \left(\frac{j}{\mu}-V\left(j \right)\right) > \frac{\epsilon r-2}{\mu} \right)\nonumber\\
&=  \pr \left( \max_{j=\rd{\epsilon r}-1,\ldots, \rd{ (\mu+ \epsilon) r}-1} \alpha\left(\frac{j}{\mu}-V\left(j \right)\right) > \alpha\frac{\epsilon r-2}{\mu} \right)\label{exp_conv_eq_8}\\
&=  \pr \left( \exp\left\{\max_{j=\rd{\epsilon r}-1,\ldots, \rd{ (\mu+ \epsilon) r}-1} \alpha\left(\frac{j}{\mu}-V\left(j \right)\right)\right\} > \exp\left\{\alpha\frac{\epsilon r-2}{\mu}\right\} \right)\nonumber\\
&=  \pr \left( \max_{j=\rd{\epsilon r}-1,\ldots, \rd{ (\mu+ \epsilon) r}-1}\exp\left\{ \alpha\left(\frac{j}{\mu}-V\left(j \right)\right)\right\} > \exp\left\{\alpha\frac{\epsilon r-2}{\mu}\right\} \right)\nonumber\\
&\leq  \E \left[ \exp\left\{ \alpha\left(\frac{\rd{ (\mu+ \epsilon) r}-1}{\mu}-V\left(\rd{ (\mu+ \epsilon) r}-1 \right)\right)\right\}\right] \exp\left\{-\alpha\frac{\epsilon r-2}{\mu}\right\}\label{exp_conv_eq_9}\\
&=  \E \left[ \exp\left\{\sum_{j=1}^{\rd{ (\mu+ \epsilon) r}-1} \alpha\left(\frac{1}{\mu}-v_j\right)\right\}\right] \exp\left\{-\alpha\frac{\epsilon r-2}{\mu}\right\}\nonumber\\
&=  \E \left[ \prod_{j=1}^{\rd{ (\mu+ \epsilon) r}-1}\exp\left\{ \alpha\left(\frac{1}{\mu}-v_j\right)\right\}\right] \exp\left\{-\alpha\frac{\epsilon r-2}{\mu}\right\}\nonumber\\
&= \prod_{j=1}^{\rd{ (\mu+ \epsilon) r}-1} \E \left[ \exp\left\{ \alpha\left(\frac{1}{\mu}-v_j\right)\right\}\right] \exp\left\{-\alpha\frac{\epsilon r-2}{\mu}\right\}\label{exp_conv_eq_10}\\
&= \E \left[ \exp\left\{ \alpha\left(\frac{1}{\mu}-v_1\right)\right\}\right]^{\rd{ (\mu+ \epsilon) r}-1} \exp\left\{-\alpha\frac{\epsilon r-2}{\mu}\right\}\label{exp_conv_eq_11}\\
&= \left(\e^{\ell(\alpha)}\right)^{\rd{ (\mu+ \epsilon) r}-1} \exp\left\{-\alpha\frac{\epsilon r-2}{\mu}\right\}\label{exp_conv_eq_12}\\
&= \exp\left\{-(\rd{ (\mu+ \epsilon) r}-1)\left(\alpha\frac{\epsilon r-2}{\mu(\rd{ (\mu+ \epsilon) r}-1)}-\ell(\alpha)\right)\right\},\label{exp_conv_eq_13}
\end{align}
where \eqref{exp_conv_eq_6} is by the fact that $T(t)\leq t$ for all $t\geq 0$ and Lemma \ref{ld_result_2}; \eqref{exp_conv_eq_7} is  by Lemma \ref{st_lemma} and \eqref{exp_conv_eq_4_0}, $\alpha>0$ is an arbitrary constant in \eqref{exp_conv_eq_8}, we use Doob's inequality (cf. Theorem 5.4.2 of \citet{dur10}) in order to obtain the inequality in \eqref{exp_conv_eq_9}, \eqref{exp_conv_eq_10} and \eqref{exp_conv_eq_11} are by the i.i.d. property of the sequence $\{v_j,j\in\N_+\}$, \eqref{exp_conv_eq_12} is by \eqref{exp_conv_eq_2}. Since \eqref{exp_conv_eq_13} is valid for all $\alpha>0$, by \eqref{exp_conv_eq_3} we have 
\begin{equation}\label{exp_conv_eq_14}
 \pr \left(\eta_1^r< r \right) \leq  \exp\left\{-(\rd{ (\mu+ \epsilon) r}-1)\Lambda_1\left(\frac{\epsilon r-2}{\mu(\rd{ (\mu+ \epsilon) r}-1)}\right)\right\}.
\end{equation}

Second,
\begin{align}
&\pr \left(\eta_2^r< r \right) \nonumber\\
&=  \pr \left(  \sup_{0\leq t\leq r} \left( \ot{V}^{i,r}_1\left(\ru{\mu t- \epsilon r} \right) - t \right) >0 \right)\nonumber\\
&\leq  \pr \left( \max_{j= 1,\ldots, \ru{ (\mu- \epsilon) r}} \left(\ot{V}^{i,r}_1(j) - \frac{j-1+\epsilon r}{\mu}\right) >0 \right)\label{exp_conv_eq_15}\\
&=  \pr \left(\max_{j= 1,\ldots, \ru{ (\mu- \epsilon) r}} \left(\ot{V}^{i,r}_1(j) - \frac{j}{\mu}\right) > \frac{\epsilon r-1}{\mu} \right)\nonumber\\
&=  \pr \left(\left( v_{S(T(ir))+1}-\frac{1}{\mu}\right)\vee \left[ \left( v_{S(T(ir))+1}-\frac{1}{\mu}\right)+ \max_{j= 2,\ldots, \ru{ (\mu- \epsilon) r}} \left(\ot{V}^{i,r}(j) - \frac{j-1}{\mu}\right)\right] > \frac{\epsilon r-1}{\mu} \right)\label{exp_conv_eq_16}\\
&\leq  2\pr \left(v_{S(T(ir))+1}-\frac{1}{\mu} > \frac{\epsilon r-1}{2\mu}\right) + \pr\left(\max_{j= 2,\ldots, \ru{ (\mu- \epsilon) r}} \left(\ot{V}^{i,r}(j) - \frac{j-1}{\mu}\right) > \frac{\epsilon r-1}{2\mu} \right),\label{exp_conv_eq_17}
\end{align}
where \eqref{exp_conv_eq_15} is by the fact that $\epsilon<\mu$, \eqref{exp_conv_eq_16} is by \eqref{exp_conv_eq_4_0}. By the same way we derive \eqref{exp_conv_eq_14}, we can derive that the second probability in \eqref{exp_conv_eq_17} is less than or equal to 
\begin{equation}\label{exp_conv_eq_18}
 \exp\left\{-(\ru{(\mu- \epsilon) r}-1)\Lambda_2\left(\frac{\epsilon r-1}{2\mu(\ru{ (\mu- \epsilon) r}-1)}\right)\right\}.
\end{equation}

Next, let us consider the first probability in \eqref{exp_conv_eq_17}, which is less than or equal to
\begin{align}
& 2\pr \left(v_{S(T(ir))+1}-\frac{1}{\mu} > \frac{\epsilon r-1}{2\mu},\;S(r^n-r)\leq 2\mu (r^n-r)-1\right) + 2\pr \left(S(r^n-r)> 2\mu (r^n-r)-1\right) \nonumber\\
& \hspace{2cm}\leq 2\pr \left(\max_{j\in\{1,2,\ldots,2\mu (r^n-r)\}} \left(v_j-\frac{1}{\mu} \right)> \frac{\epsilon r-1}{2\mu}\right) +  2\pr \left(S(r^n-r)> 2\mu (r^n-r)-1\right) \nonumber \\
& \hspace{2cm}\leq 4\mu (r^n-r) \pr \left(\left(v_1-\frac{1}{\mu} \right)> \frac{\epsilon r-1}{2\mu}\right) +  2\pr \left(S(r^n-r)> 2\mu (r^n-r)-1\right) \label{exp_conv_eq_19}.
\end{align}
By Lemmas \ref{ld_result_1} and \ref{ld_result_3}, there exists an $r_4>0$ and strictly positive constants $C_{12}$ and $C_{13}$ independent of $r$ and $i$ such that if $r\geq r_4$, the sum in \eqref{exp_conv_eq_19} is less than or equal to $C_{12} (r^n-r)\e^{-C_{13}r}$. Next, let us consider \eqref{exp_conv_eq_14} and \eqref{exp_conv_eq_18}. There exists an $r_5>0$ and $\epsilon_2>0$ such that if $r\geq r_5$, we have
\begin{equation*}
\frac{\epsilon r-2}{\mu(\rd{ (\mu+ \epsilon) r}-1)} >\epsilon_2,\quad\text{and}\quad \frac{\epsilon r-1}{2\mu(\ru{ (\mu- \epsilon) r}-1)}>\epsilon_2.
\end{equation*}
By Lemma \ref{fl_transform}, $\Lambda_i(\epsilon_2)>0$ and $\Lambda_i(x)$ is nondecreasing for all $x\geq 0$ and $i\in\{1,2\}$. Hence, if $r>r_5$, the sum of the terms in \eqref{exp_conv_eq_14} and \eqref{exp_conv_eq_18} is less than or equal to $C_{14}\e^{-C_{15}r}$, where $C_{14}:=2\exp\{\Lambda_1(\epsilon_2)\vee\Lambda_2(\epsilon_2)\}$ and $C_{15}:=(\mu- \epsilon) (\Lambda_1(\epsilon_2)\wedge\Lambda_2(\epsilon_2))$. Lastly, let $r_2:=r_4\vee r_5$, $C_6:= C_{12}\vee C_{14}$, and $C_7:= C_{13}\wedge C_{15}$, then we obtain \eqref{exp_conv_r_2}.

\subsubsection{Proof of \eqref{exp_conv_r_1}}\label{exp_conv_r_1_proof}

Note that proving \eqref{exp_conv_r_1} is equivalent to proving that there exist $r_2\geq 0$ and strictly positive constants $C_4$ and $C_5$ which are independent of $r$ such that if $r\geq r_2$,
\begin{equation}\label{exp_conv_eq_20}
\pr \left(  \sup_{0\leq t \leq r} \left| S(ir+t)-S(ir) -\mu t \right | >\epsilon r \right) \leq C_4r^n\e^{-C_5r}.
\end{equation}
Hence we will explain how to prove \eqref{exp_conv_eq_20}. Let for all $m\in\N$
\begin{equation*}
\ot{\mF}_m :=\sigma\left\{v_k, k\in\{1,2,\ldots,m+1\}\right\}.
\end{equation*}
Then, it is easy to see that $S(ir)$ is a stopping time with respect to the filtration $\{\ot{\mF}_m,m\in\N\}$. The rest of the proof is very similar to the one of \eqref{exp_conv_r_2} (cf. Appendix \ref{exp_conv_r_2_proof}), the major differences are that 
\begin{enumerate}
\item the definitions in \eqref{exp_conv_eq_4_0} should be updated as
\begin{equation*}
\ot{V}^{i,r} (k) := \sum_{j=S(ir)+2}^{S(ir)+k} v_j,\quad \forall k\in\{2,3,\ldots\},\qquad \ot{V}^{i,r}_1 (k) := \sum_{j=S(ir)+1}^{S(ir)+k} v_j,\quad \forall k\in\N_+,
\end{equation*}

\item  we need to use Theorem 4.1.3 of \citet{dur10} instead of Lemma \ref{st_lemma} in order to obtain \eqref{exp_conv_eq_7}.

\end{enumerate}

\section{Proofs of the Results Stated in the Main Body}\label{standard_results}

In this section, we present the proofs of some of the results stated in the main body of the paper.

\subsection{Proof of Proposition \ref{fluid_result_prop}}\label{fluid_results}

We will present the proof for the case in which servers 1 and 2 are in light traffic and servers 3 and 5 are in heavy traffic. The other cases follow similarly. We will derive the fluid limit results first for server 1, second for servers 3 and 4, lastly for server 6. Results associated with the other servers follow similarly. By FSLLN, we have as $r\rightarrow \infty$
\begin{equation}\label{FSLLN_results}
\left(\ol{S}^r_j, j\in\mJ\cup\mA\right)\xrightarrow{a.s.}\left(\ol{S}_j, j\in\mJ\cup\mA\right)\qquad\text{u.o.c.},
\end{equation}
where for all $t\geq 0$, $\ol{S}_a(t)=\lambda_a t$, $\ol{S}_b(t)=\lambda_b t$, and $\ol{S}_j(t)=\mu_j t$ for all $j\in\mA$. After some algebra, we have for all $t\geq 0$,
\begin{align*}
\ol{Q}^r_1(t)&= \ol{X}^r_1(t)+ \mu_1^r\ol{I}^r_1(t)\\
\ol{X}^r_1(t)&:= \left(\ol{S}_a^r(t)-\lambda_a^r t\right)-\left(\ol{S}_1^r(\ol{T}_1^r(t))-\mu_1^r \ol{T}_1^r(t)\right)+(\lambda_a^r-\mu_1^r)t.
\end{align*} 
Since server 1 works in a work-conserving fashion (cf. \eqref{w_conserving_1}) and by the uniqueness of the solution of the Skorokhod problem (cf. Proposition B.1 of \citet{bel01}) with respect to $\ol{X}^r_1$, we have $\mu_1^r\ol{I}^r_1=\Psi(\ol{X}_1^r)$ and $\ol{Q}^r_1=\Phi(\ol{X}_1^r)$ for all $r$. By \eqref{FSLLN_results}, $\sup_{0\leq t\leq T}\left|\ol{S}_a^r(t)-\lambda_a^r t\right|\xrightarrow{a.s.} 0$ for all $T\geq 0$. Since $\ol{T}_1^r(t)\leq t$ by definition, we have for all $T\geq 0$
\begin{equation}\label{s1_fluid_1}
\sup_{0\leq t\leq T}\left|\ol{S}_1^r(\ol{T}_1^r(t))-\mu_1^r \ol{T}_1^r(t) \right|\leq \sup_{0\leq t\leq T}\left|\ol{S}_1^r(t)-\mu_1^r t \right| \xrightarrow{a.s.} 0,\quad\text{ as $r\rightarrow \infty$}.
\end{equation}
Therefore, $\ol{X}^r_1 \xrightarrow{a.s.} \ol{X}_1$ u.o.c., where $\ol{X}_1(t)=(\lambda_a-\mu_1)t$. Since the mappings $\Psi$ and $\Phi$ are continuous, we have $\mu_1^r\ol{I}^r_1\xrightarrow{a.s.}\Psi(\ol{X}_1)=:\mu_1\ol{I}_1$ u.o.c. and $\ol{Q}^r_1\xrightarrow{a.s.}\ol{Q}_1=\Phi(\ol{X}_1)$ u.o.c. Since $\lambda_a<\mu_1$, then $\ol{I}_1 (t)=(1-\lambda_a/\mu_1)t$ for all $t\geq 0$ and $\ol{Q}_1=\textbf{0}$. Then by \eqref{idle_def_1}, $\ol{T}_1^r\xrightarrow{a.s.}\ol{T}_1$ u.o.c. where $\ol{T}_1 (t)=(\lambda_a/\mu_1)t$.

Next, let us consider server 3. After some algebra, we have for each $t\geq 0$
\begin{align*}
\ol{Q}^r_3(t)&= \ol{X}^r_3(t)+ \mu_3^r\ol{I}^r_3(t)\\
\ol{X}^r_3(t)&:= \left(\ol{S}_1^r(\ol{T}_1^r(t))-\mu_1^r \ol{T}_1^r(t)\right)-\left(\ol{S}_3^r(\ol{T}_3^r(t))-\mu_3^r \ol{T}_3^r(t)\right)+\mu_1^r \ol{T}_1^r(t)-\mu_3^rt.
\end{align*} 
By the uniqueness of the solution of the Skorokhod problem with respect to $\ol{X}^r_3$, we have $\mu_3^r\ol{I}^r_3=\Psi(\ol{X}_3^r)$ and $\ol{Q}^r_3=\Phi(\ol{X}_3^r)$ for all $r$. By \eqref{s1_fluid_1}, we have $\left(\ol{S}_1^r(\ol{T}_1^r)-\mu_1^r \ol{T}_1^r\right)\xrightarrow{a.s.} \textbf{0}$ u.o.c. Since $\ol{T}_3^r(t)\leq t$ for all $t\geq 0$, we have $\left(\ol{S}_3^r(\ol{T}_3^r)-\mu_3^r \ol{T}_3^r\right)\xrightarrow{a.s.} \textbf{0}$ u.o.c. by \eqref{FSLLN_results}. Since $\lambda_a=\mu_3$ and $\ol{T}_1 (t)=(\lambda_a/\mu_1)t$ for all $t\geq 0$, we have $\ol{X}^r_3 \xrightarrow{a.s.} \ol{X}_3=\textbf{0}$ u.o.c. Then, we have $\mu_3^r\ol{I}^r_3\xrightarrow{a.s.} \Psi(\ol{X}_3) =\textbf{0}$ u.o.c. and $\ol{Q}^r_3\xrightarrow{a.s.}\ol{Q}_3=\Phi(\ol{X}_3)=\textbf{0}$ u.o.c. Then by \eqref{idle_def_1}, $\ol{T}_3^r\xrightarrow{a.s.}\ol{T}_3$ u.o.c. where $\ol{T}_3 (t)=t=(\lambda_a/\mu_3) t$ for all $t\geq 0$.

Let us consider the workload process in server 4. For all $t\geq 0$,
\begin{align*}
\ol{W}^r_4(t)&= \ol{Q}^r_4(t)+\frac{\mu_{A}^r}{\mu_{B}^r}\ol{Q}^r_5(t)= \ol{U}^r_4(t) + \mu_{A}^r\ol{I}^r_4(t),\\
\ol{U}^r_4(t)&:= \ol{S}_1^r(\ol{T}_1^r(t))-\left(\ol{S}_{A}^r(\ol{T}_{A}^r(t))-\mu_{A}^r \ol{T}_{A}^r(t)\right) + \frac{\mu_{A}^r}{\mu_{B}^r}\Big(\ol{S}_2^r(\ol{T}_2^r(t))-\left(\ol{S}_{B}^r(\ol{T}_{B}^r(t))-\mu_{B}^r \ol{T}_{B}^r(t)\right)\Big)-\mu_{A}^r t.
\end{align*}
Since we consider work-conserving policies, server 4 can be idle only if there is no workload in bufers 4 and 5. Therefore, we can use the uniqueness of the solution of the Skorokhod problem with respect to $\ol{U}^r_4$ and we have $\mu_{A}^r\ol{I}^r_4=\Psi(\ol{U}_4^r)$ and $\ol{W}^r_4=\Phi(\ol{U}_4^r)$ for all $r$. By \eqref{FSLLN_results} and the logic of \eqref{s1_fluid_1}, we have $\left(\ol{S}_{A}^r(\ol{T}_{A}^r)-\mu_{A}^r \ol{T}_{A}^r\right)\xrightarrow{a.s.} \textbf{0}$ u.o.c. and $\left(\ol{S}_{B}^r(\ol{T}_{B}^r)-\mu_{B}^r \ol{T}_{B}^r\right)\xrightarrow{a.s.} \textbf{0}$ u.o.c. By the fluid limit results of servers 1 and 2 and Theorem 5.3 of \citet{che01} (Random Time-Change theorem), we have $\ol{S}_1^r(\ol{T}_1^r)\xrightarrow{a.s.}\ol{S}_1(\ol{T}_1)$ u.o.c. and $\ol{S}_2^r(\ol{T}_2^r)\xrightarrow{a.s.} \ol{S}_2(\ol{T}_2)$ u.o.c. Therefore, $\ol{U}_4^r\xrightarrow{a.s.} \ol{U}_4$ u.o.c. where $\ol{U}_4=\textbf{0}$ by Assumption \ref{assumption_rate} part \ref{ht_assumption_s4_1}. Then, by the continuity of the mappings $\Psi$ and $\Phi$, we have $\mu_{A}^r\ol{I}^r_4\xrightarrow{a.s.} \textbf{0}$ u.o.c. and $\ol{W}^r_4\xrightarrow{a.s.}\textbf{0}$ u.o.c. Since $0\leq \ol{Q}^r_4\leq \ol{W}^r_4$ and $0\leq \ol{Q}^r_5\leq (\mu_{B}^r/\mu_{A}^r)\ol{W}^r_4$, we have $\ol{Q}^r_4\xrightarrow{a.s.}\textbf{0}$ u.o.c. and $\ol{Q}^r_5\xrightarrow{a.s.}\textbf{0}$ u.o.c. Note that, 
\begin{align}
\ol{Q}^r_4(t)&= \ol{S}_1^r(\ol{T}_1^r(t))-\left(\ol{S}_{A}^r(\ol{T}_{A}^r(t))-\mu_{A}^r \ol{T}_{A}^r(t)\right)-\mu_{A}^r \ol{T}_{A}^r(t),\label{s4_fluid_1}\\
\ol{Q}^r_5(t)&= \ol{S}_2^r(\ol{T}_2^r(t))-\left(\ol{S}_{B}^r(\ol{T}_{B}^r(t))-\mu_{B}^r \ol{T}_{B}^r(t)\right)-\mu_{B}^r \ol{T}_{B}^r(t).\label{s4_fluid_2}
\end{align}
By substituting the above results in \eqref{s4_fluid_1} and \eqref{s4_fluid_2}, we have $\ol{T}_{A}^r\xrightarrow{a.s.}\ol{T}_{A}$ u.o.c. where $\ol{T}_{A} (t)=(\lambda_a/\mu_{A}) t$ and $\ol{T}_{B}^r\xrightarrow{a.s.}\ol{T}_{B}$ u.o.c. where $\ol{T}_{B} (t)=(\lambda_b/\mu_{B}) t$ for all $t\geq 0$.

Next, we consider server 6. Let $\ol{Z}^r_6:=\ol{Q}^r_7\wedge \ol{Q}^r_8$. After some algebra, for each $T\geq 0$
\begin{align*}
\ol{Z}^r_6(t)&= \ol{X}^r_6(t)+\mu_6^r\ol{I}^r_6(t),\\
\ol{X}^r_6(t)&:= \ol{S}_3^r(\ol{T}_3^r(t))\wedge \ol{S}_{A}^r(\ol{T}_{A}^r(t))- \left(\ol{S}_6^r(\ol{T}_6^r(t))-\mu_6^r \ol{T}_6^r(t)\right)-\mu_6^r t.
\end{align*}
Since server 6 can be idle only if $\ol{Z}^r_6(t)=0$ (cf. \eqref{w_conserving_2}), we can use the solution of the Skorokhod problem with respect to $\ol{X}^r_6$ and we have $\mu_6^r\ol{I}^r_6=\Psi(\ol{X}^r_6)$ and $\ol{Z}^r_6=\Phi(\ol{X}^r_6)$. By \eqref{FSLLN_results} and the logic of \eqref{s1_fluid_1}, we have $\left(\ol{S}_6^r(\ol{T}_6^r)-\mu_6^r \ol{T}_6^r\right)\xrightarrow{a.s.} \textbf{0}$ u.o.c. By the fluid limit results of servers 3 and 4 and random time-change theorem, we have $\ol{S}_3^r(\ol{T}_3^r)\xrightarrow{a.s.}\ol{S}_3(\ol{T}_3)$ u.o.c. and $\ol{S}_{A}^r(\ol{T}_{A}^r)\xrightarrow{a.s.} \ol{S}_{A}(\ol{T}_{A})$ u.o.c. Therefore, $\ol{X}_6^r\xrightarrow{a.s.} \ol{X}_6$ u.o.c. where $\ol{X}_6(t)=(\lambda_a -\mu_6)t$ for all $t\geq 0$.  Then, by the fact that $\lambda_a<\mu_6$ and the continuity of the mappings $\Psi$ and $\Phi$, we have $\ol{I}^r_6\xrightarrow{a.s.} \ol{I}_6$ u.o.c. where $\ol{I}_6(t)=(1-\lambda_a/\mu_6)t$ and $\ol{Z}^r_6\xrightarrow{a.s.}\textbf{0}$ u.o.c. By \eqref{idle_def_1}, $\ol{T}_6^r\xrightarrow{a.s.}\ol{T}_6$ u.o.c. where $\ol{T}_6 (t)=(\lambda_a/\mu_6)t$. Note that,
\begin{equation}\label{s6_fluid_1}
\ol{Q}^r_7(t)= \ol{S}_3^r(\ol{T}_3^r(t))-\ol{S}_6^r(\ol{T}_6^r(t)),\qquad \ol{Q}^r_8(t)= \ol{S}_{A}^r(\ol{T}_{A}^r(t))-\ol{S}_6^r(\ol{T}_6^r(t)).
\end{equation}
Therefore, we have $\ol{Q}^r_7\xrightarrow{a.s.}\textbf{0}$ u.o.c. and $\ol{Q}^r_8\xrightarrow{a.s.}\textbf{0}$ u.o.c. by \eqref{s6_fluid_1} and random time-change theorem.  

The convergence results related to $\ol{Q}^r_k$, $k\in\{2,6,9,10\}$ and $\ol{T}^r_j$, $j\in\{2,5,7\}$ follow similarly, hence we skip them.

\subsection{Proof of Proposition \ref{general_conv}}\label{proof_general_conv}

We provide the proof for the case $\mH=\{2,3,5\}$, i.e servers 2, 3, and 5 are in heavy traffic together with server 4 but server 1 is in light traffic. The proofs of the other cases follow similarly. We derive the weak convergence result first for server 1, second for servers 2, 3, 4, and 5. 

We use the Skorokhod's representation theorem to obtain the equivalent distributional representations of the processes in \eqref{Skorokhod_1_0} (for which we use the same symbols and call ``Skorokhod represented versions'') such that all Skorokhod represented versions of the processes are defined in the same probability space and the weak convergence in \eqref{Skorokhod_1_0} is replaced by almost sure convergence. Then we have \eqref{Skorokhod_1_1} and let us consider the Skorokhod represented versions of the processes in \eqref{Skorokhod_1_1}. 

We first consider server 1. By \eqref{Skorokhod_1_1}, random time-change theorem, and Theorem 4.1 of \citet{whi80} (continuity of addition), we have $\oh{S}_a^r-\oh{S}^r_1(\ol{T}_1^r) \xrightarrow{a.s.}\ot{S}_a-\ot{S}_1(\ol{T}_1)$ u.o.c. Since server 1 works in a work-conserving fashion and is in light traffic, we have $\oh{Q}_1^r=\Phi\big(\oh{X}_1^r\big)$ and $\mu_1^r \oh{I}_1^r=\Psi\big(\oh{X}_1^r\big)$ by \eqref{netput_eq_1} and \eqref{queue_matrix}. Then by Lemma 6.4 (ii) of \citet{che01}, we have
\begin{equation}\label{s1_diffusion_1}
\left(\mu_1^r\oh{I}_1^r+r(\lambda_a^r-\mu_1^r)e,\; \oh{Q}_1^r\right)\xrightarrow{a.s.} \left(-\ot{S}_a+\ot{S}_1(\ol{T}_1),\;\bm{0}\right),\qquad\text{u.o.c.}
\end{equation}

Next let us consider servers 2, 3, 4, and 5. Let $\oh{\bm{Q}}^r_{\mH}$, $\oh{\bm{X}}^r_{\mH}$, and $\oh{\bm{I}}^r_{\mH}$ be $|\mH |$-dimensional vectors derived from the vectors $\oh{\bm{Q}}^r$, $\oh{\bm{X}}^r$, and $\oh{\bm{I}}^r$  (cf. \eqref{vector_def_1}) by deleting the rows corresponding to each $i$, $i\notin\mH$, respectively. Let $\bm{R}_{\mH}^r$ denote the $|\mH |\times |\mH |$-dimensional matrix derived from $\bm{R}^r$ (cf. \eqref{vector_def_2}) by deleting the rows and columns corresponding to each $i$, $i\notin\mH$. Then by \eqref{vector_def_1}, \eqref{vector_def_2}, and \eqref{queue_matrix}, we have
\begin{equation*}
\bm{\oh{Q}^r}_{\mH} = \bm{\oh{X}^r}_{\mH} + \bm{R^r}_{\mH} \bm{\oh{I}^r}_{\mH}.
\end{equation*}
By the fact that all servers work in a work-conserving fashion (cf. \eqref{w_conserving}) and Theorem 7.2 of \citet{che01}, $\bm{\oh{Q}^r}_{\mH} =\phi\big(\bm{\oh{X}^r}_{\mH}\big)$ and $\bm{\oh{I}^r}_{\mH} =\psi\big(\bm{\oh{X}^r}_{\mH}\big)$, where $(\phi,\psi)$ is the multidimensional reflection mapping which is Lipschitz continuous under the uniform norm. Hence let us first focus on $\bm{\oh{X}^r}_{\mH}$. By Assumption \ref{assumption_rate} Parts \ref{ht_assumption_s4_1}, \ref{ht_assumption_s4_2}, and \ref{ht_assumption_s1235},  \eqref{workload_process}, \eqref{brew_process}, \eqref{netput_eq}, \eqref{Skorokhod_1_1}, random time-change theorem, continuity of addition, and the fact that $\oh{Q}_1^r\xrightarrow{a.s.} \bm{0}$ u.o.c., 
\begin{align}
&\oh{\bm{X}}_{\mH}^r=\begin{bmatrix}
 \oh{X}_2^r  \\
 \oh{X}_3^r  \\
 \oh{X}_6^r  \\
 \oh{X}_4^r +\frac{\mu_{A}^r}{\mu_{B}^r} \oh{X}_5^r 
\end{bmatrix}\xrightarrow{a.s.}
\begin{bmatrix}
 \ot{S}_b - \ot{S}_2\big(\ol{T}_2\big)+\theta_2 e  \\
\ot{S}_a- \ot{S}_3\big(\ol{T}_3\big)+ \theta_3 e \\
\ot{S}_2\big(\ol{T}_2\big) - \ot{S}_5\big(\ol{T}_5\big)+ (\theta_5-\theta_2) e  \\
\ot{S}_a- \ot{S}_A\big(\ol{T}_A\big) + \frac{\mu_{A}}{\mu_{B}}\left( \ot{S}_2\big(\ol{T}_2\big) - \ot{S}_B\big(\ol{T}_B\big)\right) + \left(\theta_4- \frac{\mu_{A}}{\mu_{B}}\theta_2\right) e
\end{bmatrix}\nonumber\\
&\hspace{12cm}=: \ot{\bm{X}}_{\mH}\qquad \text{u.o.c.}\label{netput_conv}
\end{align}

After some algebra, it is possible to see that  $\ot{\bm{X}}_{\mH}$ is a Brownian motion starting from $0_{|\mH|}$ with drift vector $\bm{\theta}_{\mH}$ and covariance matrix $\Sigma_{\mH}$. By the continuity of the multidimensional reflection mapping, we have  $\bm{\ot{Q}}_{\mH} =\phi\big(\bm{\ot{X}}_{\mH}\big)$ and $\bm{\ot{I}}_{\mH} =\psi\big(\bm{\ot{X}}_{\mH}\big)$, where
\begin{equation*}
\bm{\ot{Q}}_{\mH}=\left(\ot{Q}_2,\ot{Q}_3,\ot{Q}_6,\ot{W}_4\right)',\qquad \bm{\ot{I}}_{\mH}=\left(\ot{I}_2,\ot{I}_3,\ot{I}_5,\ot{I}_4\right)'.
\end{equation*}
By Definition 3.1 of \citet{wil98b}, $\bm{\ot{Q}}_{\mH}$ is an SRBM associated with the data $\big(P_{|\mH |},\bm{\theta}_{\mH},\Sigma_{\mH},\bm{R}_{\mH},0_{|\mH |}\big)$.

 Since the Skorokhod represented version of the processes have the same joint distribution with the original ones, when the Skorokhod represented versions of the processes converge almost surely u.o.c., then the original processes weakly converge. In other words, corresponding to the original processes, we have
\begin{equation}\label{diffusion_result_2}
\left(\oh{Q}^r_1,\oh{Q}^r_2,\oh{Q}^r_3,\oh{Q}^r_6,\oh{W}^r_4\right) \Longrightarrow \left(\textbf{0},\ot{Q}_2,\ot{Q}_3,\ot{Q}_6,\ot{W}_4\right).
\end{equation}


\subsection{Proof of Lemma \ref{up_lemma}}\label{up_lemma_proof}

\begin{align}
& S_{A}^r(T_{A}^r(\tau_{2n-1}^r)+t)=\sup\left\{k\in \N:\sum_{l=1}^k v_{A}^r(l) \leq T_{A}^r(\tau_{2n-1}^r)+t\right\}\nonumber\\
&\hspace{2cm}= \sup\left\{k+A_n^r:k\in \N,\;\sum_{l=1}^{A_n^r} v_{A}^r(l) + \sum_{l=1}^{k} v_{A}^r(l+A_n^r) \leq T_{A}^r(\tau_{2n-1}^r)+t\right\}\nonumber,\\
&\hspace{2cm}\geq\sup\left\{k+A_n^r:k\in \N,\;T_{A}^r(\tau_{2n-1}^r)+ \sum_{l=1}^{k} v_{A}^r(l+A_n^r) \leq T_{A}^r(\tau_{2n-1}^r)+t\right\}\label{lemma_up_eq_3}\\
&\hspace{2cm}= A_n^r+ \sup\left\{k\in \N: \sum_{l=1}^{k} v_{A}^r(l+A_n^r) \leq t\right\} = A_n^r + E_{A}^{r,n}(t),\label{lemma_up_eq_4}
\end{align}
where \eqref{lemma_up_eq_3} is by the fact that $\sum_{l=1}^{A_n^r} v_{A}^r(l) \leq  T_{A}^r(\tau_{2n-1}^r)$ by definition of $A_n^r$. Similar to the derivation of \eqref{lemma_up_eq_4}, we can get the following result:
\begin{equation}\label{lemma_up_eq_5}
S_{3}^r(T_{3}^r(\tau_{2n-1}^r)+t)= B_n^r+ \sup\left\{k\in \N: \sum_{l=1}^{k} v_{3}^r(l+B_n^r) \leq t\right\} = B_n^r + E_{3}^{r,n}(t).
\end{equation}
Note that, we have equality sign in \eqref{lemma_up_eq_5} unlike the greater than or equal to sign in \eqref{lemma_up_eq_3}. This is because there is a service completion in server 3 exactly at $\tau_{2n-1}^r$ by construction (cf. \eqref{tau_odd}). Therefore, we have $\sum_{l=1}^{B_n^r} v_{3}^r(l) =  T_{3}^r(\tau_{2n-1}^r)$, and this gives us the equality sign in \eqref{lemma_up_eq_5}. Then,  for all $t\geq 0$, 
 \begin{equation*}
 S_{3}^r(T_{3}^r(\tau_{2n-1}^r)+t) - S_{A}^r(T_{A}^r(\tau_{2n-1}^r)+t) \leq E_{3}^{r,n}(t) - E_{A}^{r,n}(t) +1,
\end{equation*}
by \eqref{lemma_up_eq_4}, \eqref{lemma_up_eq_5}, and the fact that $B_n^r-A_n^r=1$.

\subsection{Proof of Lemma \ref{ata_modification}}\label{proof_of_ata_modification}

First, we present the following results from the literature which will be used later in the proof.

\begin{lemma}\label{lemma_ata_kumar}
(\textbf{Lemma 9 of \citet{ata05}}) Given $\epsilon>0$ and $T>2/\epsilon$, we have for each $j\in\mJ\cup\mA$ and $\alpha>0$
\begin{equation*}
\pr \left(\sup_{0\leq t \leq T} |S_j(t)-x_jt|\geq \epsilon T \right)\leq \frac{C_j(\epsilon,\alpha)}{T^{1+\alpha}},
\end{equation*}
where
\begin{align}
 &C_j(\epsilon,\alpha) := \left(\frac{2+2\alpha}{1+2\alpha}\right)^{2+2\alpha}\left(\frac{18(2+2\alpha)^{3/2}}{(1+2\alpha)^{1/2}}\right)^{2+2\alpha} \nonumber\\
 &\hspace{5cm}  \times\E\left[|v_j(1)-1/x_j|^{2+2\alpha} \right] \left[ \left(\frac{4x_j^2(x_j+\epsilon)}{\epsilon^2}\right)^{1+\alpha}+\left(\frac{4x_j^3}{\epsilon^2}\right)^{1+\alpha}\right], \label{ata_constant}
\end{align}
if $j\in\mJ$, then $x_j=\lambda_j$; and if $j\in\mA$, then $x_j=\mu_j$.
\end{lemma}
The proof of Lemma \ref{lemma_ata_kumar} can be seen in \citet{ata05} (Note that we have fixed a small typo in \eqref{ata_constant} by replacing $\big((2+2\alpha)/(1+2\alpha)\big)$ with $\big((2+2\alpha)/(1+2\alpha)\big)^{2+2\alpha}$. This typo does not affect the results of \citet{ata05}). Since we assume exponential moment condition for the service times (cf. Assumption \ref{assumption_moment}), Lemma \ref{lemma_ata_kumar} holds for all $\alpha>0$. However, the proof presented in \citet{ata05} requires a weaker moment assumption.

\begin{lemma}\label{hall_and_heyde}
(\textbf{Equation 3.67 of \citet{hal80}})For any martingale with differences $Z_i,1\leq i\leq n$ and any $p\geq 2$, we have
\begin{equation*}
\E \left[ \left| \sum_{i=1}^n Z_i \right|^p \right]\leq (18pq^{1/2})^{p} n^{p/2} \max_{1\leq i\leq n} \E \left[ \left|Z_i \right|^p \right],
\end{equation*}
where $q=(1-p^{-1})^{-1}$.
\end{lemma}

The proof of Lemma \ref{ata_modification} is a modification of the proof of Lemma \ref{lemma_ata_kumar}. We present the proof for the case $j=3$, the other case follows similarly. Let us define
\begin{equation}\label{lemma_up_eq_8}
\mE_{n,3}^r := \left\{B_n^r\leq \ru{(\mu_3^r+\epsilon_9)r^2T} \right\},
\end{equation}
where $\epsilon_9>0$ is an arbitrary constant. Then, 
\begin{align}
& \sum_{n=1} ^{N^r}  \pr \left( \sup_{0\leq t<r^{\gamma}T} \left|E_{3}^{r,n}(t) - \mu_{3}^r t\right| > \epsilon_2r^{\gamma}T,\;\tau_{2n-1}^r \leq r^2T \right)\nonumber\\
& \hspace{1cm} \leq \sum_{n=1} ^{N^r} \pr \left( \tau_{2n-1}^r \leq r^2T,\; \left(\mE_{n,3}^r\right)^c \right) +\sum_{n=1} ^{N^r}  \pr \left( \sup_{0\leq t<r^{\gamma}T} \left|E_{3}^{r,n}(t) - \mu_{3}^r t\right| > \epsilon_2r^{\gamma}T,\;\mE_{n,3}^r \right),\label{lemma_up_eq_9}
\end{align}
Let us consider the first sum in \eqref{lemma_up_eq_9}.
\begin{align}
& \sum_{n=1} ^{N^r} \pr \left( \tau_{2n-1}^r \leq r^2T,\; \left(\mE_{n,3}^r\right)^c \right)\leq \sum_{n=1} ^{N^r} \pr \left( \tau_{2n-1}^r \leq r^2T,\; S_{3}^r(T_{3}^r(\tau_{2n-1}^r)) > (\mu_3^r+\epsilon_9)r^2T\right)\nonumber\\
& \hspace{1cm} \leq  \sum_{n=1} ^{N^r} \pr \left( S_{3}^r(r^2T)> (\mu_3^r+\epsilon_9)r^2T\right)\leq  \sum_{n=1} ^{N^r} \pr \left( \sup_{0\leq t \leq r^2T} \left|S_{3}^r(t)-\mu_3^r t\right| > \epsilon_9 r^2T\right)\nonumber\\
&\hspace{1cm} = N^r \pr \left( \sup_{0\leq t \leq r^2T} \left|S_{3}^r(t)-\mu_3^r t\right| > \epsilon_9 r^2T\right) \leq N^r\frac{C_3^r(\epsilon_9,\alpha_2)}{(r^2T)^{1+\alpha_2}}\rightarrow 0\qquad\text{as $r\rightarrow\infty$},\label{lemma_up_eq_10}
\end{align}
where $\alpha_2>0$ is an arbitrary constant, and \eqref{lemma_up_eq_10} is by Lemma \ref{lemma_ata_kumar} and the fact that $N^r=O(r^2)$. Moreover, for each $r$, $C_3^r(\epsilon_9,\alpha_2)$ is the constant defined in Lemma \ref{lemma_ata_kumar} associated with the renewal process $S_3^r$. It is straightforward to see that $C_3^r(\epsilon_9,\alpha_2)$ is bounded above by a constant uniformly in $r$, hence we get the convergence result in \eqref{lemma_up_eq_10}.

Next, let us consider the second sum in \eqref{lemma_up_eq_9}. Let us fix $n$ and $r$. Let 
\begin{equation*}
\eta:=\inf \left\{ t\in [0,r^{\gamma}T): |E_{3}^{r,n}(t) - \mu_{3}^r t| > \epsilon_2r^{\gamma}T\right\},
\end{equation*}
where $\inf \{\emptyset\}:= \infty$ for completeness. Then
\begin{equation*}
\pr \left( \sup_{0\leq t<r^{\gamma}T} \left|E_{3}^{r,n}(t) - \mu_{3}^r t\right| > \epsilon_2r^{\gamma}T,\;\mE_{n,3}^r \right)= \pr \left(\eta < r^{\gamma}T,\;\mE_{n,3}^r \right).
\end{equation*}
Let $\ot{V}_3^r(-k):= 0$ for all $k\in\N$ and $\ot{V}_3^r(k):= \sum_{i=1}^k v_3^r(B_n^r+i)$. Then,
\begin{align*}
\{\left|E_{3}^{r,n}(t) - \mu_{3}^r t\right| > \epsilon_2r^{\gamma}T\}&=\{E_{3}^{r,n}(t) > \mu_{3}^r t + \epsilon_2r^{\gamma}T \}\cup \{E_{3}^{r,n}(t) < \mu_{3}^r t - \epsilon_2r^{\gamma}T \}\\
& \subseteq \{E_{3}^{r,n}(t) > \rd{\mu_{3}^r t + \epsilon_2r^{\gamma}T} \}\cup \{E_{3}^{r,n}(t) < \ru{\mu_{3}^r t - \epsilon_2r^{\gamma}T} \}\\
& = \{\ot{V}_3^r(\rd{\mu_{3}^r t + \epsilon_2r^{\gamma}T}) < t \}\cup \{\ot{V}_3^r(\ru{\mu_{3}^r t - \epsilon_2r^{\gamma}T}) > t \}.
\end{align*}
Let us define
\begin{align*}
&\eta_1:=\inf \left\{ t\in [0,r^{\gamma}T): \ot{V}_3^r(\rd{\mu_{3}^r t + \epsilon_2r^{\gamma}T}) < t \right\},\\
&\eta_2:=\inf \left\{ t\in [0,r^{\gamma}T): \ot{V}_3^r(\ru{\mu_{3}^r t - \epsilon_2r^{\gamma}T}) > t \right\}.
\end{align*}
Then $\eta \geq \eta_1\wedge \eta_2$. Thus
\begin{equation}\label{lemma_up_eq_13_1}
\pr \left(\eta < r^{\gamma}T,\;\mE_{n,3}^r \right)\leq \pr \left(\eta_1\wedge \eta_2< r^{\gamma}T,\;\mE_{n,3}^r \right) \leq  \pr \left(\eta_1< r^{\gamma}T,\;\mE_{n,3}^r \right) +  \pr \left( \eta_2< r^{\gamma}T,\;\mE_{n,3}^r \right).
\end{equation}
We will consider the two probabilities after the second inequality in \eqref{lemma_up_eq_13_1} separately. First,
\begin{equation*}
\eta_1=\inf \left\{ t\in [0,r^{\gamma}T): \ot{V}_3^r(\rd{\mu_{3}^r t + \epsilon_2r^{\gamma}T})-\frac{\rd{\mu_{3}^r t + \epsilon_2r^{\gamma}T}}{\mu_{3}^r} < t - \frac{\rd{\mu_{3}^r t + \epsilon_2r^{\gamma}T}}{\mu_{3}^r}  \right\}.
\end{equation*}
Next, let us define
\begin{equation}\label{eta_tilda}
\ot{\eta}_1:=\inf \left\{ t\in [0,r^{\gamma}T): \ot{V}_3^r(\rd{\mu_{3}^r t + \epsilon_2r^{\gamma}T})-\frac{\rd{\mu_{3}^r t + \epsilon_2r^{\gamma}T}}{\mu_{3}^r} < -\frac{\epsilon_2r^{\gamma}T}{2\mu_{3}^r}  \right\}.
\end{equation}
When $r$ is sufficiently large, $\epsilon_2r^{\gamma}T >2$. This implies
\begin{equation*}
-\frac{\epsilon_2r^{\gamma}T}{2\mu_{3}^r} > t - \frac{\rd{\mu_{3}^r t + \epsilon_2r^{\gamma}T}}{\mu_{3}^r} 
\end{equation*} 
and $\pr \left( \eta_1< r^{\gamma}T,\;\mE_{n,3}^r \right)\leq \pr \left( \ot{\eta}_1< r^{\gamma}T,\;\mE_{n,3}^r \right)$ when $r$ is sufficiently large. Then, 
\begin{align}
&\pr \left( \ot{\eta}_1<r^{\gamma}T,\;\mE_{n,3}^r \right) \leq \pr \left( \sup_{i=1,2,\ldots,\rd{(\mu_{3}^r + \epsilon_2)r^{\gamma}T}} \left|\ot{V}_3^r(i)-\frac{i}{\mu_{3}^r}\right| > \frac{\epsilon_2r^{\gamma}T}{2\mu_{3}^r} ,\;\mE_{n,3}^r \right)\nonumber\\
&\leq \E \left[ \left(\sup_{i=1,2,\ldots,\rd{(\mu_{3}^r + \epsilon_2)r^{\gamma}T}} \left|\ot{V}_3^r(i)-\frac{i}{\mu_{3}^r}\right|\right)^{2+2\alpha}\I\left(\mE_{n,3}^r\right) \right] \left(\frac{\epsilon_2r^{\gamma}T}{2\mu_{3}^r}\right)^{-(2+2\alpha)}\nonumber \\
& = \E \left[ \sum_{j=1}^{\ru{(\mu_3^r+\epsilon_9)r^2T}}\left(\sup_{i=1,2,\ldots,\rd{(\mu_{3}^r + \epsilon_2)r^{\gamma}T}} \left|\ot{V}_3^r(i)-\frac{i}{\mu_{3}^r}\right|\I\left( B_n^r=j \right)\right)^{2+2\alpha} \right] \left(\frac{\epsilon_2r^{\gamma}T}{2\mu_{3}^r}\right)^{-(2+2\alpha)} \nonumber \\
& = \sum_{j=1}^{\ru{(\mu_3^r+\epsilon_9)r^2T}}  \E \left[ \left(\sup_{i=1,2,\ldots,\rd{(\mu_{3}^r + \epsilon_2)r^{\gamma}T}} \left | V_3^r(i+j)-V_3^r(j)-\frac{i}{\mu_{3}^r}\right|\I\left( B_n^r=j \right)\right)^{2+2\alpha} \right] \left(\frac{\epsilon_2r^{\gamma}T}{2\mu_{3}^r}\right)^{-(2+2\alpha)}\nonumber \\
& \leq \ru{(\mu_3^r+\epsilon_9)r^2T}\;\E \left[ \left(\sup_{i=1,2,\ldots,\rd{(\mu_{3}^r + \epsilon_2)r^{\gamma}T}} \left | V_3^r(i)-\frac{i}{\mu_{3}^r}\right|\right)^{2+2\alpha} \right] \left(\frac{\epsilon_2r^{\gamma}T}{2\mu_{3}^r}\right)^{-(2+2\alpha)}, \label{lemma_up_eq_14}
\end{align}
where $\alpha>0$ is an arbitrary constant such that $4/(1+\alpha)<\gamma$; the first inequality is by the definition of $\ot{\eta}_1$ (cf. \eqref{eta_tilda}); the second inequality is by Markov's inequality; the first equality is by the definition of $\mE_{n,3}^r$ (cf. \eqref{lemma_up_eq_8}); the second equality is by the fact that $(V_3^r(j+i)-V_3^r(j))\I\left( B_n^r=j \right)=\ot{V}_3^r(i)\I\left( B_n^r=j \right)$ for all $i,j\in\N_+$; and the last inequality is due to the fact that $\I\left( B_n^r=j \right)\leq 1$ and $V_3^r(j+i)-V_3^r(j)\overset{d}{=}V_3^r(i)$ by the i.i.d. property of $\{v_3^r(i),i\in\N_+\}$. Note that
\begin{equation*}
\left\{\left(V_3^r(i)-\frac{i}{\mu_{3}^r}\right),i\in\N_+\right\}
\end{equation*}
is a martingale. Then, by the $L^p$ maximum inequality (cf. Theorem 5.4.3 of \citet{dur10}), we see that the term in \eqref{lemma_up_eq_14} is less than or equal to
\begin{equation}\label{lemma_up_eq_15}
\ru{(\mu_3^r+\epsilon_9)r^2T} \left(\frac{2+2\alpha}{1+2\alpha}\right)^{(2+2\alpha)}\E \left[   \left | V_3^r(\rd{(\mu_{3}^r + \epsilon_2)r^{\gamma}T})-\frac{\rd{(\mu_{3}^r + \epsilon_2)r^{\gamma}T}}{\mu_{3}^r}\right |^{2+2\alpha} \right] \left(\frac{\epsilon_2r^{\gamma}T}{2\mu_{3}^r }\right)^{-(2+2\alpha)}.
\end{equation}
Then, by using Lemma \ref{hall_and_heyde} on \eqref{lemma_up_eq_15}, the term in \eqref{lemma_up_eq_15} is less than or equal to
\begin{align}
&\ru{(\mu_3^r+\epsilon_9)r^2T}  \left(\frac{2+2\alpha}{1+2\alpha}\right)^{(2+2\alpha)} \left(\frac{18(2+2\alpha)^{3/2}}{(1+2\alpha)^{1/2}}\right)^{(2+2\alpha)} \left(\rd{(\mu_{3}^r + \epsilon_2)r^{\gamma}T}\right)^{1+\alpha} \nonumber\\
& \hspace{7cm} \times\E \left[ \left |v_3^r(1)-\frac{1}{\mu_{3}^r}\right |^{2+2\alpha} \right] \left(\frac{\epsilon_2r^{\gamma}T}{2\mu_{3}^r}\right)^{-(2+2\alpha)}, \label{lemma_up_eq_16}
\end{align}
so does $\pr \left( \ot{\eta}_1<r^{\gamma}T,\;\mE_{n,3}^r \right) $.

Now, we will consider the second probability in \eqref{lemma_up_eq_13_1}. First,
\begin{equation*}
\eta_2=\inf \left\{ t\in [0,r^{\gamma}T): \ot{V}_3^r(\ru{\mu_{3}^r t - \epsilon_2r^{\gamma}T})-\frac{\ru{\mu_{3}^r t - \epsilon_2r^{\gamma}T}}{\mu_{3}^r} > t - \frac{\ru{\mu_{3}^r t - \epsilon_2r^{\gamma}T}}{\mu_{3}^r}  \right\}.
\end{equation*}
Next, let us define
\begin{equation*}
\ot{\eta}_2:=\inf \left\{ t\in [0,r^{\gamma}T): \ot{V}_3^r(\ru{\mu_{3}^r t - \epsilon_2r^{\gamma}T})-\frac{\ru{\mu_{3}^r t - \epsilon_2r^{\gamma}T}}{\mu_{3}^r} > \frac{\epsilon_2r^{\gamma}T}{2\mu_{3}^r}  \right\}.
\end{equation*}
When $r$ is sufficiently large, $\epsilon_2r^{\gamma}T >2$. This implies
\begin{equation*}
\frac{\epsilon_2r^{\gamma}T}{2\mu_{3}^r} < t - \frac{\ru{\mu_{3}^r t - \epsilon_2r^{\gamma}T}}{\mu_{3}^r} 
\end{equation*} 
and $\pr \left( \eta_2< r^{\gamma}T,\;\mE_{n,3}^r \right)\leq \pr \left( \ot{\eta}_2< r^{\gamma}T,\;\mE_{n,3}^r \right)$ when $r$ is sufficiently large. Moreover, since $\mu_3^r>\epsilon_2$ when $r$ is sufficiently large, $\ru{(\mu_{3}^r - \epsilon_2)r^{\gamma}T}>0$ when $r$ is sufficiently large. Then, 
\begin{align}
\pr \left( \ot{\eta}_2<r^{\gamma}T,\;\mE_{n,3}^r \right) &\leq \pr \left( \sup_{i=1,2,\ldots,\ru{(\mu_{3}^r - \epsilon_2)r^{\gamma}T}} \left|\ot{V}_3^r(i)-\frac{i}{\mu_{3}^r}\right| > \frac{\epsilon_2r^{\gamma}T}{2\mu_{3}^r} ,\;\mE_{n,3}^r \right).\nonumber\\
& \leq \pr \left( \sup_{i=1,2,\ldots,\rd{(\mu_{3}^r + \epsilon_2)r^{\gamma}T}} \left|\ot{V}_3^r(i)-\frac{i}{\mu_{3}^r}\right| > \frac{\epsilon_2r^{\gamma}T}{2\mu_{3}^r} ,\;\mE_{n,3}^r \right)\label{lemma_up_eq_19_0}\\
&\leq \ru{(\mu_3^r+\epsilon_9)r^2T}  \left(\frac{2+2\alpha}{1+2\alpha}\right)^{(2+2\alpha)}\left(\frac{18(2+2\alpha)^{3/2}}{(1+2\alpha)^{1/2}}\right)^{(2+2\alpha)} \nonumber\\
& \hspace{1.5cm} \times \left(\ru{(\mu_{3}^r + \epsilon_2)r^{\gamma}T}\right)^{1+\alpha} \E \left[ \left |v_3^r(1)-\frac{1}{\mu_{3}^r}\right |^{2+2\alpha} \right] \left(\frac{\epsilon_2r^{\gamma}T}{2\mu_{3}^r}\right)^{-(2+2\alpha)}, \label{lemma_up_eq_19}
\end{align}
where the inequality in \eqref{lemma_up_eq_19_0} is by the fact that $\epsilon_2r^{\gamma}T\geq 1$ when $r$ is sufficiently large (this implies $\ru{(\mu_{3}^r - \epsilon_2)r^{\gamma}T}\leq \rd{(\mu_{3}^r + \epsilon_2)r^{\gamma}T}$ when $r$ is sufficiently large) and \eqref{lemma_up_eq_19} is by \eqref{lemma_up_eq_14} and \eqref{lemma_up_eq_16}.

By the exponential moment assumption (cf. Assumption \ref{assumption_moment}), both of the right hand sides in \eqref{lemma_up_eq_16} and \eqref{lemma_up_eq_19} are $O(r^{2-\gamma(1+\alpha)})$. By \eqref{lemma_up_eq_13_1} and the fact that $4/(1+\alpha)<\gamma$, we have
\begin{align*}
&\sum_{n=1} ^{N^r}  \pr \left( \sup_{0\leq t<r^{\gamma}T} \left|E_{3}^{r,n}(t) - \mu_{3}^r t\right| > \epsilon_2r^{\gamma}T,\;\mE_{n,3}^r \right)\\
&\hspace{3cm} \leq \sum_{n=1} ^{N^r} \left(\pr \left(\eta_1< r^{\gamma}T,\;\mE_{n,3}^r \right) +  \pr \left( \eta_2< r^{\gamma}T,\;\mE_{n,3}^r \right)\right) \\
&\hspace{3cm} = N^r O(r^{2-\gamma(1+\alpha)}) = O(r^{4-\gamma(1+\alpha)}) \rightarrow 0,\qquad\text{as $r\rightarrow\infty$}.
\end{align*}

By the fact that $A_n^r= B_n^r-1$ and using the same technique, we can also prove that
\begin{equation*}
\sum_{n=1} ^{N^r} \pr \left( \sup_{0\leq t<r^{\gamma}T} \left|E_{A}^{r,n}(t) - \mu_{A}^r t\right| > \epsilon_2r^{\gamma}T,\;\tau_{2n-1}^r\leq r^2T\right)\rightarrow 0,\qquad\text{as $r\rightarrow\infty$}.
\end{equation*}

\subsection{Proof of Lemma \ref{down_lemma}}\label{down_lemma_proof}

Consider the sum in \eqref{lemma_down_eq_9}, which is less than or equal to
\begin{align}
& \sum_{n=1} ^{N_2^r} \pr \left(\tau_{2n}^{r,*} \leq r^2(T+\epsilon_4),\; S_{B}^r(T_{B}^r(\tau_{2n}^{r,*})) > N_3^r\right) \label{lemma_down_eq_9_1}\\
& \hspace{2cm} + \sum_{n=1} ^{N_2^r} \pr \left(\tau_{2n}^{r,*} \leq r^2(T+\epsilon_4),\; S_{b}^r(\tau_{2n}^{r,*}) > N_3^r\right) \label{lemma_down_eq_9_2}\\
& \hspace{3.5cm}+ \sum_{n=1} ^{N_2^r} \pr \left(\tau_{2n}^{r,*} \leq r^2(T+\epsilon_4),\; Q_2^r(\tau_{2n}^{r,*})> \left(\left( \mu_2\wedge\mu_B-\kappa-\lambda_b \right)\wedge \epsilon_5 \right)r/2 \right) . \label{lemma_down_eq_9_3}
\end{align}
First, consider the sum in \eqref{lemma_down_eq_9_1}, which is less than or equal to
\begin{align}
& \sum_{n=1} ^{N_2^r} \pr \left(\tau_{2n}^{r,*} \leq r^2(T+\epsilon_4),\; S_{B}^r(T_{B}^r(\tau_{2n}^{r,*})) > (\lambda_b^r+\mu_{B}^r)r^2(T+\epsilon_4)\right)\nonumber\\
& \hspace{3cm}\leq  \sum_{n=1} ^{N_2^r} \pr \left( S_{B}^r(r^2(T+\epsilon_4))> (\lambda_b^r+\mu_{B}^r)r^2(T+\epsilon_4)\right)\nonumber\\
& \hspace{3cm}\leq  N_2^r \pr \left( \sup_{0\leq t \leq r^2(T+\epsilon_4)} \left|S_{B}^r(t)-\mu_{B}^r t\right| > \lambda_b^r r^2(T+\epsilon_4)\right)\nonumber\\
& \hspace{3cm}\leq N_2^r\frac{C_{B}^r(\lambda_b^r,\alpha_1)}{\left[r^2(T+\epsilon_4)\right]^{1+\alpha_1}}\rightarrow 0\qquad\text{as $r\rightarrow\infty$},\label{lemma_down_eq_10}
\end{align}
where $\alpha_1>0$ is any arbitrary constant and \eqref{lemma_down_eq_10} is by Lemma \ref{lemma_ata_kumar} (cf. Appendix \ref{proof_of_ata_modification}). For each $r$, $C_B^r(\lambda_b^r,\alpha_1)$ is the constant defined in Lemma \ref{lemma_ata_kumar} associated with the renewal process $S_B^r$. It is straightforward to see that $C_B^r(\lambda_b^r,\alpha_1)$ is bounded above by a constant uniformly in $r$, hence we get the convergence result in \eqref{lemma_down_eq_10}. We can show that the sum in \eqref{lemma_down_eq_9_2} converges to $0$ by the same way we derive \eqref{lemma_down_eq_10}. 

Next let us look at the sum in \eqref{lemma_down_eq_9_3}, which is less than or equal to
\begin{equation}\label{lemma_down_eq_10_01}
\sum_{n=1} ^{N_2^r} \pr \left(\sup_{0\leq t \leq r^2(T+\epsilon_4)} Q_2^r(t)> \frac{\left(\left( \mu_2\wedge\mu_B-\kappa-\lambda_b \right)\wedge \epsilon_5 \right)r}{2} \right).
\end{equation} 
Let us construct a hypothetical $GI/GI/1$ queue, where the buffer length process is denoted by $Q^{(1)}_2$, the interarrival and service time sequences are $\left\{ v_b^r(i)(\lambda_b^r/(\lambda_b+\kappa/2)),\; i\in\N_+ \right\}$ and $\left\{ v^r_2(i)(\mu_2^r/(\mu_2-\kappa/2)),\; i\in\N_+ \right\}$, respectively\footnote{Since $v_b^r(i)(\lambda_b^r/(\lambda_b+\kappa/2))= \bar{v}_b(i)/(\lambda_b+\kappa/2)$ and $v_2^r(i)(\mu_2^r/(\mu_2-\kappa/2))= \bar{v}_2(i)/(\mu_2-\kappa/2)$ for all $i\in\N_+$ (cf. Section \ref{seq_of_systems}), the corresponding sequences are independent of $r$}, and $Q^{(1)}_2(0):=Q_2^r(0)=0$. Hence, the arrival and service rates in the hypothetical $GI/GI/1$ queue are $\lambda_b+\kappa/2$ and $\mu_2-\kappa/2$, respectively. By Assumption \ref{assumption_rate} Parts \ref{lambda_rate_assumption} and \ref{mu_rate_assumption}, there exists an $r_6\geq 1$ such that if $r\geq r_6$, the term in \eqref{lemma_down_eq_10_01} is less than or equal to 
\begin{equation}\label{lemma_down_eq_10_1}
\sum_{n=1} ^{N_2^r} \pr \left(\sup_{0\leq t \leq r^2(T+\epsilon_4)} Q^{(1)}_2(t)> \frac{\left(\left( \mu_2\wedge\mu_B-\kappa-\lambda_b \right)\wedge \epsilon_5 \right)r}{2} \right).
\end{equation} 
The term in \eqref{lemma_down_eq_10_1} converges to 0 as $r\rightarrow\infty$ by the fact the corresponding server is in light traffic and Proposition \ref{light_traffic_conv_rate_lemma}. 
\subsection{Proof of Lemma \ref{opt_solution_2}}\label{opt_solution_2_proof}

In the optimization problem \eqref{opt_problem}, replacing $q_4$ with $w_4-(\mu_{A}/\mu_{B}) q_5$ gives us the following equivalent optimization problem which has only one decision variable.
\begin{subequations}\label{opt_problem_4}
\begin{align}
\min\quad &h_4 \left(w_4 - \frac{\mu_{A}}{\mu_{B}} q_5\right) + h_5 q_5 + h_7\left(w_4 - q_3 - \frac{\mu_{A}}{\mu_{B}} q_5\right)^+ + h_8\left(q_3-w_4 + \frac{\mu_{A}}{\mu_{B}} q_5\right)^+  \nonumber\\
&\hspace{8cm}  + h_9\left(q_6-q_5\right)^+ + h_{10}\left(q_5-q_6\right)^+ ,\label{opt_problem_4_1}\\
&\hspace{2cm}\text{s.t.}\quad 0\leq q_5 \leq \frac{\mu_{B}}{\mu_{A}} w_4. \label{opt_problem_4_2}
\end{align}
\end{subequations}
The objective function \eqref{opt_problem_4_1} is the sum of six different functions each of which is convex, continuous, and piecewise linear with respect to the decision variable $q_5$. Since sum of finitely many convex (continuous, piecewise linear) functions is convex (continuous, piecewise linear), the objective function \eqref{opt_problem_4_1} is also convex, continuous, and piecewise linear with respect to the decision variable $q_5$. Then, an optimal solution should be either in the boundaries of the feasible region of $q_5$, which are $0$ and $(\mu_{B}/\mu_{A}) w_4$ (cf. \eqref{opt_problem_4_2}), or in one of the break points which are in the interval $\big[0,(\mu_{B}/\mu_{A}) w_4\big]$ of the convex, continuous, and piecewise linear objective function \eqref{opt_problem_4_1}, which are $(\mu_{B}/\mu_{A}) (w_4-q_3)^+$ and $q_6\wedge (\mu_{B}/\mu_{A}) w_4$. Therefore, the optimal solution, denoted by $q_5^*$, is such that 
\begin{equation*}
q_5^*\in\left\{0,\; \frac{\mu_{B}}{\mu_{A}} w_4,\;  \frac{\mu_{B}}{\mu_{A}} (w_4-q_3)^+,\; q_6\wedge  \frac{\mu_{B}}{\mu_{A}} w_4 \right\}.
\end{equation*}
Then the optimal solution set in Table \ref{sln_table} follows by plugging $q_5^*$ in the equality constraint \eqref{opt_problem_2}.

\subsection{Proof of Lemma \ref{opt_solution_3}}\label{opt_solution_3_proof}

For each $j\in\{1,2,\ldots,n\}$, $h_j \left(q_j -q_{j,2} \right)^+$ is a convex and continuous function of the decision variables $q_i$, $i\in\{1,2,\ldots,n\}$. Since sum of finitely many convex and continuous functions is also convex and continuous, then the objective function \eqref{opt_problem_1_1} is also convex and continuous function of the decision variables $q_j$, $j\in\{1,2,\ldots,n\}$. Then, a local optimum solution of the optimization problem \eqref{opt_problem_1_0} is also a global optimum and we will prove that the solution \eqref{opt_problem_1_4} is a local optimum. First, it is easy to see that the solution \eqref{opt_problem_1_4} satisfies the constraints \eqref{opt_problem_1_2} and \eqref{opt_problem_1_3}, so it is feasible. We will prove that the solution \eqref{opt_problem_1_4} is a local optimum by showing that any deviation from this solution does not improve the objective function value. 

Let us fix an arbitrary $j\in\{1,2,\ldots,n-1\}$ and consider $q_j$. Note that the cost incurred due to $q_j$ is 0 (cf. Remark \ref{comment_sln}). Thus, decreasing $q_j$ cannot decrease $h_j \left(q_j -q_{j,2} \right)^+$ but may increase the objective function value because at least one $q_i$, $i\in\{1,2,\ldots,n\}\backslash\{j\}$ will increase by \eqref{opt_problem_1_2}. Therefore, decreasing $q_j$, $j\in\{1,2,\ldots,n-1\}$ does not improve the objective function value.

Next, let us increase $q_j$ by $\epsilon$ where $\epsilon>0$ but sufficiently small. Increasing $q_j$ by $\epsilon$ increases the objective function value by at most $h_j\epsilon$. Since $h_i \left(q_i -q_{i,2} \right)^+=0$ for all $i\in\{1,2,\ldots,n-1\}\backslash\{j\}$, decreasing $q_i$, $i\in\{1,2,\ldots,n-1\}\backslash\{j\}$ will not decrease the objective function value. Therefore, increasing $q_j$ but decreasing $q_i$, $i\in\{1,2,\ldots,n-1\}\backslash\{j\}$ does not improve the objective function value. Next, let us decrease $q_n$. There are two cases to consider. First, suppose that $q_j\geq q_{j,2}$. Then, increasing $q_j$ by $\epsilon$ increases the objective function value exactly by $h_j\epsilon$. In this case, decreasing $q_n$ by $\mu_n \epsilon/\mu_j$ (cf. \eqref{opt_problem_1_2}) can decrease the objective function value by at most $h_n\mu_n \epsilon/\mu_j$ which is less than or equal to $h_j\epsilon$ because $h_j\mu_j\geq h_n\mu_n$, thus the net change in the objective function value is nonnegative. Second, suppose that $q_j< q_{j,2}$. Then, increasing $q_j$ by $\epsilon$ will not increase the objective function value for sufficiently small $\epsilon$. Note that, when $q_j< q_{j,2}$, by \eqref{opt_problem_1_41} and \eqref{opt_problem_1_42}
\begin{equation}\label{eq_opt_problem_1_sln}
\frac{q_{j,2}}{\mu_j} > \frac{q_{j}}{\mu_j} =
\begin{cases} 
w,&\mbox{if $j=1$},\\
\left(\cdots\left(\left( w - \frac{q_{1,2}}{\mu_1}\right)^+ -  \frac{q_{2,2}}{\mu_2}\right)^+ - \cdots -  \frac{q_{j-1,2}}{\mu_{j-1}}\right)^+ ,&\mbox{if $j\in\{2,3,\ldots,n-1\}$}.
\end{cases}
\end{equation}
This implies that if $j=1$, $q_n=0$ by \eqref{opt_problem_1_2} and \eqref{opt_problem_1_3}, and if $j\in\{2,3,\ldots,n-1\}$, then by \eqref{opt_problem_1_43} and \eqref{eq_opt_problem_1_sln},
\begin{align*}
q_n &= \mu_n \left(\cdots\left(\left(\cdots\left(\left( w - \frac{q_{1,2}}{\mu_1}\right)^+ -  \frac{q_{2,2}}{\mu_2}\right)^+ - \cdots - \frac{q_{j-1,2}}{\mu_{j-1}}\right)^+- \frac{q_{j,2}}{\mu_j}\right)^+ -\cdots  \frac{q_{n-1,2}}{\mu_{n-1}}\right)^+, \\
& = \mu_n \left(\cdots\left(\left(\frac{q_{j}}{\mu_j} -  \frac{q_{j,2}}{\mu_j}\right)^+  -  \frac{q_{j+1,2}}{\mu_{j+1}}\right)^+- \cdots -  \frac{q_{n-1,2}}{\mu_{n-1}}\right)^+, \\
& = \mu_n \left(\cdots\left( 0 -  \frac{q_{j+1,2}}{\mu_{j+1}}\right)^+- \cdots -  \frac{q_{n-1,2}}{\mu_{n-1}}\right)^+, \\
& = 0.
\end{align*}
Hence, when $q_j< q_{j,2}$, $q_n=0$ so it cannot be decreased by \eqref{opt_problem_1_3}. As a result, increasing $q_j$ by $\epsilon$ does not improve the objective function value. Therefore, increasing or decreasing $q_j$ does not improve the objective function value for all $j\in\{1,2,\ldots,n-1\}$.

Lastly, let us consider $q_n$. If we increase (decrease) $q_n$, then some of the $q_j$, $j\in\{1,2,\ldots,n-1\}$ must increase (decrease) by \eqref{opt_problem_1_2} and \eqref{opt_problem_1_3}. Since the latter change does not improve the objective function value, changing the value of $q_n$ does not improve the objective function value. Therefore, the solution \eqref{opt_problem_1_4} is a local optimum and also a global optimum.

\subsection{Proof of Lemma \ref{three_type_lemma}}\label{three_type_lemma_proof}

By decreasing the number of decision variables from 4 to 2 by using the equality constraints and then simplifying the notation in the optimization problem \eqref{DCP_12}, we get the following equivalent optimization problem:
\begin{subequations}\label{opt_type_3}
\begin{align}
& \min\quad h_a \left(w_5 - \frac{\mu_{A}}{\mu_{B1}} q_6  - q_4\right)^+ +h_b \left(q_6\vee q_7\right) + h_c \left(\frac{\mu_{C}}{\mu_{B2}}(w_6-q_7) - q_9 \right)^+,\label{opt_type_3_1}\\
& \hspace{2cm} \text{s.t.}\quad\;  0\leq  q_6 \leq   \frac{\mu_{B1}}{\mu_{A}} w_5,\qquad 0\leq q_7\leq w_6, \label{opt_type_3_2}
\end{align}
\end{subequations}
where the decision variables are $q_6$ and $q_7$. Note that increasing $q_6$ or $q_7$ decrease the first and third terms but increase the second term in the objective function \eqref{opt_type_3_1}. Hence, the key point is to compare the total decrease in the first and third terms with the increase in the second term in \eqref{opt_type_3_1}, when we increase $q_6$ and $q_7$ from 0 to their corresponding upper bounds. We solve the optimization problem \eqref{opt_type_3} case by case.  

 First consider the case $q_4\geq w_5$. Then, the first term in the objective function \eqref{opt_type_3_1} is 0 for all values of $q_6$ satisfying the constraints in \eqref{opt_type_3_2}. Since increasing $q_6$ increases the second term in \eqref{opt_type_3_1}, $q_6=0$ in the optimal solution. Since $h_b \geq h_c \mu_C/ \mu_{B2}$ by assumption, then increasing $q_7$ does not decrease the objective function value, thus $q_7=0$ in the optimal solution. In the case $q_9\geq (\mu_{C}/\mu_{B2})w_6$, $q_6=q_7=0$ in the optimal solution by the same logic and the fact that $h_b \geq h_a \mu_A/  \mu_{B1}$. Therefore, the last case to consider is $w_5> q_4$ and $(\mu_{C}/\mu_{B2})w_6  > q_9$ and we will consider it in two cases.

\paragraph{Case 1: $\bm{h_b \geq h_a \mu_{A}/ \mu_{B1} + h_c \mu_{C} / \mu_{B2}}$} By the second term in \eqref{opt_type_3_1}, it is more efficient to increase $q_6$ and $q_7$ together with the same rate (if possible). Suppose that $q_6=q_7=0$. If we increase $q_6$ and $q_7$ by sufficiently small $\epsilon_{10}>0$, then the objective function value increases by ($h_b- h_a \mu_{A}/ \mu_{B1} - h_c \mu_{C} / \mu_{B2})\epsilon_{10}>0$. Therefore, we should not increase $q_6$ and $q_7$ at all, and $q_6=q_7=0$ in the optimal solution. Lastly, \eqref{DCP_sln_4} follows by the two equality constraints of the optimization problem \eqref{DCP_12}.

\paragraph{Case 2: $\bm{h_b <h_a \mu_{A}/ \mu_{B1} + h_c \mu_{C} / \mu_{B2}}$} Suppose that $q_6=q_7=0$. If we increase $q_6$ and $q_7$ by sufficiently small $\epsilon_{10}>0$, then the objective function value changes by ($h_b- h_a \mu_{A}/ \mu_{B1} - h_c \mu_{C} / \mu_{B2})\epsilon_{10}<0$. Hence, it is efficient to increase $q_6$ and $q_7$ together until either $q_6= (\mu_{B1}/\mu_{A})(w_5-q_4)$ or $q_7 = w_6 - (\mu_{B2}/\mu_{C})q_9$. After this point, if we increase at least one of the $q_6$ or $q_7$, then the objective function value increases because first or third term in \eqref{opt_type_3_1} is in its lower bound (which is 0) at this point and $h_b \geq h_a \mu_A/  \mu_{B1}$ and  $h_b \geq h_c \mu_C/ \mu_{B2}$. Therefore, $q_6=q_7=\min\{ (\mu_{B1}/\mu_{A})(w_5-q_4),\;w_6 - (\mu_{B2}/\mu_{C})q_9\}$. Then, \eqref{DCP_sln_5} follows by the two equality constraints of the optimization problem \eqref{DCP_12}.

\begin{remark}\label{type_3_intuition}
The solution \eqref{DCP_sln_5} can be explained intuitively in the following way. Note that, we start with the solution $\ot{Q}_6=\ot{Q}_7=0$ in the proof of \eqref{DCP_sln_5}, which is equivalent to giving full priority to type $b$ jobs in servers 5 and 6. Then we increase $\ot{Q}_6$ and $\ot{Q}_7$, which is equivalent to giving some of the priority to type $a$ and $c$ jobs. We do this until the point in which giving more priority to type $a$ and $b$ jobs does not increase the throughput rate of at least one of these job types because more than necessary jobs accumulates in buffers 11 or 14. Therefore, we stop giving priority to type $a$ and $b$ jobs at this point.
\end{remark}

\section{Detailed Simulation Results}\label{detailed_numerical_results}

In this section, we present the detailed results of the simulation experiments.

\begin{landscape}
\begin{table}[ht]\scriptsize
\centering \caption{Detailed results of the simulation experiments: Average queue lengths with their $95\%$ confidence intervals.}
\resizebox{1.4\textwidth}{!}{
\begin{tabular}{ c r r|r r|r r|r r|r r|r r }
\cline{2-13}
 & \multicolumn{2}{ c|}{\textbf{Proposed Policy}} & \multicolumn{2}{|c|}{\textbf{SDP Policy}} & \multicolumn{2}{|c|}{\textbf{Static Priority Policy}} & \multicolumn{2}{|c|}{\textbf{FCFS Policy}} & \multicolumn{2}{|c|}{\textbf{Randomized Policy}} & \multicolumn{2}{|c }{\textbf{Randomized-$\bm{2/3}$ Policy}}\\
\hline
\textbf{Ins.} & $\bm{\ol{Q}^p_3(i)+\ol{Q}_7^p(i)}$ & $\bm{\ol{Q}_6^p(i)+\ol{Q}_{10}^p(i)}$ & $\bm{\ol{Q}^p_3(i)+\ol{Q}_7^p(i)}$ & $\bm{\ol{Q}_6^p(i)+\ol{Q}_{10}^p(i)}$ & $\bm{\ol{Q}^p_3(i)+\ol{Q}_7^p(i)}$ & $\bm{\ol{Q}_6^p(i)+\ol{Q}_{10}^p(i)}$ & $\bm{\ol{Q}^p_3(i)+\ol{Q}_7^p(i)}$ & $\bm{\ol{Q}_6^p(i)+\ol{Q}_{10}^p(i)}$ & $\bm{\ol{Q}^p_3(i)+\ol{Q}_7^p(i)}$ & $\bm{\ol{Q}_6^p(i)+\ol{Q}_{10}^p(i)}$ & $\bm{\ol{Q}^p_3(i)+\ol{Q}_7^p(i)}$ & $\bm{\ol{Q}_6^p(i)+\ol{Q}_{10}^p(i)}$ \\
\hline
1 & 14.01 $\pm$ 0.16 & 14.43 $\pm$ 0.17 & 14.01 $\pm$ 0.16 & 14.43 $\pm$ 0.17 & 13.6 $\pm$ 0.13 & 18.21 $\pm$ 0.20 & 14.72 $\pm$ 0.14 & 14.75 $\pm$ 0.17 & 14.51 $\pm$ 0.14 & 14.72 $\pm$ 0.17 & 13.59 $\pm$ 0.15 & 17.31 $\pm$ 0.14\\

2 & 21.91 $\pm$ 0.32 & 23.9 $\pm$ 0.27 & 21.91 $\pm$ 0.32 & 23.9 $\pm$ 0.27 & 21.55 $\pm$ 0.32 & 29.58 $\pm$ 0.36 & 23.36 $\pm$ 0.24 & 23.33 $\pm$ 0.18 & 23.75 $\pm$ 0.37 & 23.14 $\pm$ 0.22 & 21.49 $\pm$ 0.27 & 28.09 $\pm$ 0.30\\

3 & 45.73 $\pm$ 0.84 & 52.44 $\pm$ 0.82 & 45.73 $\pm$ 0.84 & 52.44 $\pm$ 0.82 & 44.8 $\pm$ 0.88 & 62.75 $\pm$ 0.89 & 50.01 $\pm$ 0.69 & 49.86 $\pm$ 0.70 & 50.02 $\pm$ 0.71 & 49.66 $\pm$ 0.75 & 44.91 $\pm$ 0.79 & 60.59 $\pm$ 0.92\\

4 & 14.09 $\pm$ 0.19 & 6.53 $\pm$ 0.06 & 14.09 $\pm$ 0.19 & 6.53 $\pm$ 0.06 & 13.67 $\pm$ 0.19 & 13.28 $\pm$ 0.11 & 14.7 $\pm$ 0.17 & 7.56 $\pm$ 0.05 & 14.49 $\pm$ 0.13 & 7.71 $\pm$ 0.06 & 13.6 $\pm$ 0.15 & 12.12 $\pm$ 0.14\\

5 & 21.86 $\pm$ 0.22 & 11.63 $\pm$ 0.16 & 21.86 $\pm$ 0.22 & 11.63 $\pm$ 0.16 & 21.52 $\pm$ 0.27 & 21.35 $\pm$ 0.24 & 23.43 $\pm$ 0.28 & 12.17 $\pm$ 0.11 & 23.62 $\pm$ 0.24 & 12.18 $\pm$ 0.10 & 21.79 $\pm$ 0.29 & 19.53 $\pm$ 0.20\\

6 & 45.84 $\pm$ 0.75 & 26.69 $\pm$ 0.54 & 45.84 $\pm$ 0.75 & 26.69 $\pm$ 0.54 & 44.9 $\pm$ 0.84 & 44.9 $\pm$ 0.49 & 49.7 $\pm$ 0.58 & 25.78 $\pm$ 0.26 & 50.57 $\pm$ 0.96 & 25.57 $\pm$ 0.31 & 44.31 $\pm$ 0.74 & 41.86 $\pm$ 0.64\\

7 & 3.67 $\pm$ 0.00 & 17.27 $\pm$ 0.14 & 3.67 $\pm$ 0.00 & 17.27 $\pm$ 0.14 & 3.03 $\pm$ 0.00 & 18.07 $\pm$ 0.14 & 7.59 $\pm$ 0.05 & 14.69 $\pm$ 0.15 & 7.75 $\pm$ 0.07 & 14.73 $\pm$ 0.22 & 3.64 $\pm$ 0.01 & 17.32 $\pm$ 0.16\\

8 & 5.97 $\pm$ 0.01 & 28.03 $\pm$ 0.35 & 5.97 $\pm$ 0.01 & 28.03 $\pm$ 0.35 & 5.12 $\pm$ 0.01 & 29.11 $\pm$ 0.33 & 12.26 $\pm$ 0.11 & 23.64 $\pm$ 0.24 & 12.37 $\pm$ 0.13 & 23.58 $\pm$ 0.28 & 5.9 $\pm$ 0.01 & 27.95 $\pm$ 0.27\\

9 & 12.59 $\pm$ 0.04 & 60.99 $\pm$ 1.00 & 12.59 $\pm$ 0.04 & 60.99 $\pm$ 1.00 & 10.99 $\pm$ 0.03 & 62.74 $\pm$ 0.79 & 25.73 $\pm$ 0.32 & 50.75 $\pm$ 0.86 & 25.58 $\pm$ 0.33 & 50.15 $\pm$ 0.74 & 12.42 $\pm$ 0.04 & 60.49 $\pm$ 0.79\\

10 & 3.67 $\pm$ 0.00 & 11.95 $\pm$ 0.12 & 3.67 $\pm$ 0.00 & 11.95 $\pm$ 0.12 & 3.03 $\pm$ 0.00 & 13.25 $\pm$ 0.11 & 7.57 $\pm$ 0.06 & 7.55 $\pm$ 0.06 & 7.69 $\pm$ 0.06 & 7.75 $\pm$ 0.05 & 3.64 $\pm$ 0.01 & 12.23 $\pm$ 0.14\\

11 & 5.98 $\pm$ 0.01 & 19.13 $\pm$ 0.16 & 5.98 $\pm$ 0.01 & 19.13 $\pm$ 0.16 & 5.12 $\pm$ 0.01 & 21.55 $\pm$ 0.20 & 12.15 $\pm$ 0.10 & 12.13 $\pm$ 0.10 & 12.24 $\pm$ 0.10 & 12.22 $\pm$ 0.10 & 5.9 $\pm$ 0.01 & 19.43 $\pm$ 0.25\\

12 & 12.58 $\pm$ 0.03 & 41.4 $\pm$ 0.54 & 12.58 $\pm$ 0.03 & 41.4 $\pm$ 0.54 & 10.99 $\pm$ 0.03 & 45.94 $\pm$ 0.53 & 25.4 $\pm$ 0.29 & 25.33 $\pm$ 0.29 & 25.84 $\pm$ 0.39 & 25.59 $\pm$ 0.40 & 12.44 $\pm$ 0.04 & 41.74 $\pm$ 0.48\\

13 & 2.4 $\pm$ 0.00 & 17.95 $\pm$ 0.13 & 2.73 $\pm$ 0.00 & 17.83 $\pm$ 0.18 & 2.4 $\pm$ 0.00 & 17.95 $\pm$ 0.13 & 7.49 $\pm$ 0.06 & 14.8 $\pm$ 0.15 & 7.55 $\pm$ 0.06 & 14.58 $\pm$ 0.16 & 3.32 $\pm$ 0.01 & 17.4 $\pm$ 0.18\\

14 & 3.79 $\pm$ 0.01 & 29.18 $\pm$ 0.26 & 4.19 $\pm$ 0.01 & 29.06 $\pm$ 0.27 & 3.79 $\pm$ 0.01 & 29.18 $\pm$ 0.26 & 11.78 $\pm$ 0.10 & 23.62 $\pm$ 0.32 & 11.92 $\pm$ 0.12 & 23.79 $\pm$ 0.31 & 5.09 $\pm$ 0.01 & 28.11 $\pm$ 0.19\\

15 & 7.76 $\pm$ 0.02 & 62.26 $\pm$ 0.84 & 8.4 $\pm$ 0.02 & 62.04 $\pm$ 0.89 & 7.76 $\pm$ 0.02 & 62.26 $\pm$ 0.84 & 24.96 $\pm$ 0.28 & 49.49 $\pm$ 0.67 & 24.91 $\pm$ 0.28 & 50.05 $\pm$ 0.66 & 10.15 $\pm$ 0.03 & 60.9 $\pm$ 0.91\\

16 & 2.41 $\pm$ 0.00 & 13.34 $\pm$ 0.14 & 2.73 $\pm$ 0.00 & 12.95 $\pm$ 0.11 & 2.41 $\pm$ 0.00 & 13.34 $\pm$ 0.14 & 7.45 $\pm$ 0.05 & 7.56 $\pm$ 0.05 & 7.58 $\pm$ 0.05 & 7.66 $\pm$ 0.05 & 3.32 $\pm$ 0.01 & 12.21 $\pm$ 0.11\\

17 & 3.8 $\pm$ 0.01 & 21.44 $\pm$ 0.22 & 4.19 $\pm$ 0.01 & 20.91 $\pm$ 0.23 & 3.8 $\pm$ 0.01 & 21.44 $\pm$ 0.22 & 11.76 $\pm$ 0.11 & 12.05 $\pm$ 0.11 & 11.94 $\pm$ 0.10 & 12.23 $\pm$ 0.11 & 5.08 $\pm$ 0.01 & 19.4 $\pm$ 0.19\\

18 & 7.76 $\pm$ 0.02 & 45.28 $\pm$ 0.63 & 8.41 $\pm$ 0.02 & 45.26 $\pm$ 0.61 & 7.76 $\pm$ 0.02 & 45.28 $\pm$ 0.63 & 25.08 $\pm$ 0.30 & 25.83 $\pm$ 0.29 & 24.78 $\pm$ 0.32 & 25.62 $\pm$ 0.32 & 10.16 $\pm$ 0.02 & 41.69 $\pm$ 0.48\\

19 & 38.22 $\pm$ 0.46 & 24.1 $\pm$ 0.30 & 38.22 $\pm$ 0.46 & 24.1 $\pm$ 0.30 & 38.45 $\pm$ 0.45 & 29.22 $\pm$ 0.37 & 40.1 $\pm$ 0.44 & 23.77 $\pm$ 0.29 & 40.02 $\pm$ 0.49 & 23.56 $\pm$ 0.30 & 38.26 $\pm$ 0.41 & 28.41 $\pm$ 0.33\\

20 & 22.7 $\pm$ 0.21 & 28.23 $\pm$ 0.31 & 22.7 $\pm$ 0.21 & 28.23 $\pm$ 0.31 & 21.7 $\pm$ 0.28 & 29.38 $\pm$ 0.31 & 28.63 $\pm$ 0.33 & 23.27 $\pm$ 0.29 & 29.38 $\pm$ 0.38 & 23.77 $\pm$ 0.30 & 22.46 $\pm$ 0.25 & 28.03 $\pm$ 0.27\\

21 & 20.34 $\pm$ 0.24 & 29.25 $\pm$ 0.28 & 20.86 $\pm$ 0.27 & 28.94 $\pm$ 0.33 & 20.34 $\pm$ 0.24 & 29.25 $\pm$ 0.28 & 28.44 $\pm$ 0.31 & 23.74 $\pm$ 0.30 & 29.45 $\pm$ 0.41 & 23.38 $\pm$ 0.29 & 21.71 $\pm$ 0.36 & 28.26 $\pm$ 0.36\\

22 & 38.69 $\pm$ 0.41 & 11.79 $\pm$ 0.17 & 38.69 $\pm$ 0.41 & 11.79 $\pm$ 0.17 & 37.7 $\pm$ 0.52 & 21.36 $\pm$ 0.21 & 40.21 $\pm$ 0.57 & 12.13 $\pm$ 0.10 & 39.71 $\pm$ 0.55 & 12.27 $\pm$ 0.12 & 37.54 $\pm$ 0.50 & 19.66 $\pm$ 0.17\\

23 & 22.53 $\pm$ 0.32 & 19.18 $\pm$ 0.18 & 22.53 $\pm$ 0.32 & 19.18 $\pm$ 0.18 & 21.72 $\pm$ 0.25 & 21.24 $\pm$ 0.25 & 28.77 $\pm$ 0.29 & 12.12 $\pm$ 0.09 & 29.32 $\pm$ 0.33 & 12.25 $\pm$ 0.14 & 22.64 $\pm$ 0.25 & 19.66 $\pm$ 0.20\\

24 & 20.52 $\pm$ 0.35 & 21.53 $\pm$ 0.23 & 20.75 $\pm$ 0.39 & 20.71 $\pm$ 0.22 & 20.52 $\pm$ 0.35 & 21.53 $\pm$ 0.23 & 28.43 $\pm$ 0.38 & 12.19 $\pm$ 0.10 & 28.63 $\pm$ 0.30 & 12.19 $\pm$ 0.10 & 21.76 $\pm$ 0.30 & 19.54 $\pm$ 0.21\\

25 & 38.45 $\pm$ 0.47 & 41.92 $\pm$ 0.44 & 38.45 $\pm$ 0.47 & 41.92 $\pm$ 0.44 & 38.03 $\pm$ 0.56 & 49.34 $\pm$ 0.49 & 39.97 $\pm$ 0.46 & 39.95 $\pm$ 0.49 & 40.03 $\pm$ 0.42 & 40.24 $\pm$ 0.41 & 38.01 $\pm$ 0.44 & 47.26 $\pm$ 0.42\\

26 & 22.54 $\pm$ 0.28 & 47.98 $\pm$ 0.45 & 22.54 $\pm$ 0.28 & 47.98 $\pm$ 0.45 & 21.58 $\pm$ 0.26 & 48.6 $\pm$ 0.52 & 28.67 $\pm$ 0.28 & 40.13 $\pm$ 0.51 & 28.94 $\pm$ 0.38 & 40.05 $\pm$ 0.50 & 22.56 $\pm$ 0.33 & 47.45 $\pm$ 0.48\\

27 & 20.42 $\pm$ 0.22 & 49.32 $\pm$ 0.47 & 21.15 $\pm$ 0.33 & 48.82 $\pm$ 0.49 & 20.42 $\pm$ 0.22 & 49.32 $\pm$ 0.47 & 28.79 $\pm$ 0.30 & 39.97 $\pm$ 0.50 & 29.09 $\pm$ 0.32 & 40.1 $\pm$ 0.43 & 21.72 $\pm$ 0.34 & 47.28 $\pm$ 0.56\\

28 & 38.42 $\pm$ 0.46 & 31.2 $\pm$ 0.40 & 38.42 $\pm$ 0.46 & 31.2 $\pm$ 0.40 & 37.97 $\pm$ 0.44 & 43.79 $\pm$ 0.50 & 40.09 $\pm$ 0.42 & 28.64 $\pm$ 0.40 & 39.99 $\pm$ 0.40 & 29.1 $\pm$ 0.28 & 38.25 $\pm$ 0.51 & 41.24 $\pm$ 0.56\\

29 & 22.37 $\pm$ 0.30 & 41.59 $\pm$ 0.54 & 22.37 $\pm$ 0.30 & 41.59 $\pm$ 0.54 & 21.52 $\pm$ 0.31 & 44.12 $\pm$ 0.48 & 28.75 $\pm$ 0.35 & 28.51 $\pm$ 0.24 & 29.22 $\pm$ 0.38 & 29.25 $\pm$ 0.30 & 22.38 $\pm$ 0.25 & 41.05 $\pm$ 0.51\\

30 & 20.3 $\pm$ 0.27 & 43.93 $\pm$ 0.53 & 20.81 $\pm$ 0.33 & 43.96 $\pm$ 0.58 & 20.3 $\pm$ 0.27 & 43.93 $\pm$ 0.53 & 28.5 $\pm$ 0.35 & 28.74 $\pm$ 0.27 & 29.34 $\pm$ 0.28 & 29.48 $\pm$ 0.27 & 21.61 $\pm$ 0.26 & 41.1 $\pm$ 0.52\\

31 & 21.88 $\pm$ 0.32 & 41.75 $\pm$ 0.61 & 21.88 $\pm$ 0.32 & 41.75 $\pm$ 0.61 & 21.39 $\pm$ 0.25 & 48.81 $\pm$ 0.43 & 23.6 $\pm$ 0.21 & 40.01 $\pm$ 0.45 & 23.45 $\pm$ 0.32 & 40.24 $\pm$ 0.37 & 21.52 $\pm$ 0.31 & 47.52 $\pm$ 0.49\\

32 & 5.97 $\pm$ 0.01 & 47.48 $\pm$ 0.61 & 5.97 $\pm$ 0.01 & 47.48 $\pm$ 0.61 & 5.11 $\pm$ 0.01 & 49.18 $\pm$ 0.51 & 12.33 $\pm$ 0.10 & 40.12 $\pm$ 0.44 & 12.26 $\pm$ 0.10 & 40.05 $\pm$ 0.44 & 5.91 $\pm$ 0.01 & 47.13 $\pm$ 0.34\\

33 & 3.79 $\pm$ 0.00 & 48.75 $\pm$ 0.43 & 4.19 $\pm$ 0.01 & 48.52 $\pm$ 0.50 & 3.79 $\pm$ 0.00 & 48.75 $\pm$ 0.43 & 11.8 $\pm$ 0.07 & 39.66 $\pm$ 0.49 & 11.96 $\pm$ 0.12 & 39.93 $\pm$ 0.40 & 5.1 $\pm$ 0.01 & 46.74 $\pm$ 0.46\\

34 & 22.1 $\pm$ 0.33 & 30.72 $\pm$ 0.46 & 22.1 $\pm$ 0.33 & 30.72 $\pm$ 0.46 & 21.39 $\pm$ 0.29 & 43.98 $\pm$ 0.37 & 23.4 $\pm$ 0.26 & 28.67 $\pm$ 0.31 & 23.37 $\pm$ 0.30 & 29.26 $\pm$ 0.28 & 21.64 $\pm$ 0.30 & 41.25 $\pm$ 0.41\\

35 & 5.97 $\pm$ 0.01 & 41.51 $\pm$ 0.45 & 5.97 $\pm$ 0.01 & 41.51 $\pm$ 0.45 & 5.11 $\pm$ 0.01 & 43.48 $\pm$ 0.44 & 12.18 $\pm$ 0.08 & 28.72 $\pm$ 0.25 & 12.26 $\pm$ 0.12 & 29.25 $\pm$ 0.31 & 5.91 $\pm$ 0.01 & 41.22 $\pm$ 0.43\\

36 & 3.79 $\pm$ 0.01 & 44.01 $\pm$ 0.52 & 4.19 $\pm$ 0.01 & 43.22 $\pm$ 0.45 & 3.79 $\pm$ 0.01 & 44.01 $\pm$ 0.52 & 11.88 $\pm$ 0.12 & 28.96 $\pm$ 0.26 & 11.89 $\pm$ 0.12 & 29.03 $\pm$ 0.38 & 5.09 $\pm$ 0.01 & 40.94 $\pm$ 0.46\\
\hline
\end{tabular}}\label{detailed_results}
\end{table}
\end{landscape}


\bibliographystyle{plainnat}
\bibliography{fork_join}

\end{document}